%% file: AAGD_MC.tex
\documentclass{article}

\usepackage{microtype}
\usepackage{graphicx}
\usepackage{subfigure}
\usepackage{booktabs} 

\usepackage{array}
\usepackage{textcomp}
\usepackage{subfigure}
\usepackage{multirow}
\usepackage{graphicx}
\usepackage{array}
\usepackage[justification=centering]{caption}
\usepackage{tabularx}
\usepackage{float}
\usepackage{times}
\usepackage{caption}
\usepackage{epstopdf}
\usepackage{booktabs}
\usepackage{enumitem,calc}

\usepackage[margin=1.25in]{geometry}
\usepackage{hyperref}
\usepackage{natbib}
\usepackage{graphicx} 
\usepackage{subfigure} 
\usepackage{longtable}
\usepackage{bm}
\usepackage{algorithm}
\usepackage{algorithmic}
\usepackage{booktabs}
\usepackage{amsmath,amssymb,amsfonts}
\usepackage{multirow}
\newtheorem{theorem}{Theorem}
\newtheorem{proposition}{Proposition}

\newtheorem{lemma}{Lemma}
\newtheorem{corollary}{Corollary}

\newcommand{\ww}{\mathbf{w}}
\newcommand{\p}{\mathbf{p}}
\newcommand{\HH}{\mathbf{H}}

\newcommand{\bdelta}{ \bm{\delta}}

\newcommand{\vv}{\mathbf{v}}
\newcommand{\x}{\mathbf{x}}
\newcommand{\y}{\mathbf{y}}
\newcommand{\z}{\mathbf{z}}

\newcommand{\A}{\mathbf{A}}

\newcommand{\uu}{\mathbf{u}}
\newcommand{\w}{\mathbf{w}}

\newcommand{\bxi}{\bm{\xi}}
\newcommand{\E}{\mathbb{E}}
\newcommand{\tth}{\hat{h}}
\newcommand{\tx}{\tilde{\x}}
\newcommand{\tna}{\tilde\nabla}

\renewcommand{\a}{\mathbf{a}}
\renewcommand{\c}{\mathbf{c}}

\renewcommand{\d}{\mathbf{d}}

\DeclareMathOperator*{\argmin}{argmin}

\usepackage{graphicx} 
\usepackage{subfigure}
\newtheorem{assumption}{Assumption}

\usepackage{multirow}

\begin{document}

\title{Accelerating Asynchronous Algorithms for Convex Optimization by Momentum Compensation}

\author{\\
Cong Fang\\
\texttt{fangcong@pku.edu.cn}\\
Peking University \\
\and \\
Yameng Huang\\
\texttt{huangyameng@pku.edu.cn}\\
Peking University \\
\and \\
Zhouchen Lin\\
\texttt{zlin@pku.edu.cn} \\
Peking University \\
}

%

%
%
%
%
%
\date{Original circulated date: 27th February, 2018.}
\maketitle

\begin{abstract}
	Asynchronous algorithms have attracted much attention recently due to the crucial demands on solving large-scale optimization problems. However, the accelerated versions of asynchronous algorithms are rarely studied. In this paper, we propose the ``momentum compensation'' technique to accelerate asynchronous algorithms for convex problems. Specifically,  we first accelerate the plain Asynchronous Gradient Descent, which achieves a faster $O(1/\sqrt{\epsilon})$ (v.s. $O(1/\epsilon)$) convergence rate for non-strongly convex functions, and $O(\sqrt{\kappa}\log(1/\epsilon))$  (v.s. $ O(\kappa \log(1/\epsilon))$) for strongly convex functions to reach an $\epsilon$- approximate minimizer with the condition number $\kappa$.  We further apply the technique to accelerate modern stochastic asynchronous algorithms such as Asynchronous Stochastic Coordinate Descent and Asynchronous Stochastic Gradient Descent. Both of the resultant practical algorithms are faster than existing ones by order. To the best of our knowledge, we are the first to consider accelerated algorithms that allow updating by delayed gradients and are the first to propose truly accelerated asynchronous algorithms.  Finally, the experimental results on a shared memory system show  acceleration  leads to significant performance gains on ill-conditioned problems.
\end{abstract}

\section{Introduction}\label{introduction}
With the popularity of multi-core computers and the crucial demands for handling large-scale data in machine learning, designing parallel algorithms have attracted lots of interests in recent years. A straightforward way to implement parallelization is through synchronous update. Since each thread has to wait the precedent one to finish computing, a limited speed up caused by  serious overhead can be observed from synchronous algorithms, especially when the computation costs for each thread are different, or a large load imbalance exists. To avoid the frequent usage of synchronization operation, asynchronous algorithms are designed as a more sophisticated way for parallelization.

The main difference between asynchronous and synchronous algorithms lies  in the state of the parameters for computing the gradient. For synchronous algorithms, their results are essentially identical to the serial one with variants only on  implementation. Asynchronous algorithms are different,  because when one thread is computing the gradient, other threads might have updated the parameters. Take Asynchronous Gradient Decent as an example,  if we assign a global counter $k$ to indicate each update from any thread, the iteration can be formulated as:
\begin{eqnarray}\label{iter3}
\x^{k+1} =  \x^k - \gamma \nabla f(\x^{j(k)}),
\end{eqnarray}
where $\gamma$ is the step size and $\x^{j(k)}$ is the state of $\x$ at the reading time. Typically, $\x^{j(k)}$ can be any of $\{\x^1,\cdots,\x^k  \}$ when the parameters are updated with locks (see Section \ref{prelimi}). So for asynchronous algorithms, the gradient might be delayed.

\begin{table}[t]
	\centering
\scriptsize
	\caption{Convergence rates of asynchronous algorithms and their corresponding serial algorithms for convex optimization. ( `P' is short for optimization  problem. `T' is short for type, `S' is short for serial, `A' is short for asynchronous, and `B' is short for bounded delay assumption. $\mu$ is the strong convexity modulus, $L$ and $L_c$ are the Lipschitz  and coordinate  Lipschitz constants in Eq.~\eqref{l1} and Eq.~\eqref{lc}, respectively.)}	
	\begin{tabular} {|c|c |c |c|c|c|}\hline
		\multirow{1}*{P}&		\multirow{1}*{T} &Algorithm & Convergence Rate for NC& Convergence Rate for SC &Assumptions \\\hline\hline
		\multirow{2}*{\eqref{pro1}}&	\multirow{1}*{S}&AGD \small{\citep{nesterov1983method}} &$\sqrt{L/\epsilon}$ & $\sqrt{\frac{L}{\mu}} \log(1/\epsilon)$& \\\cline{2-6}
		&\multirow{1}*{A}&AAGD \small{(ours)} &$\tau^2\sqrt{L/\epsilon}$ & $\tau^2\sqrt{\frac{L}{\mu}}\log(1/\epsilon)$ & B \\\hline\hline
		\multirow{5}*{\eqref{AASCD1}}		&\multirow{1}*{S}&APCG \small{\citep{APCG}} &$ (n\sqrt{L}+n\sqrt{L_c})\sqrt{1/\epsilon }$ & $ n\sqrt{\frac{L_c}{\mu}}  \log(1/\epsilon)$& \\\cline{2-6}
		&\multirow{4}*{A}&\multirow{2}*{ASCD \small{\citep{ASCD1}}}&\multirow{2}*{$nL_c/\epsilon$} & \multirow{2}*{$n\frac{L_c}{\mu}\log(1/\epsilon) $} &  B, $\tau\leq \frac{\sqrt{n}L_c}{L},$\\
		&& && & NC: $\| \x\|^2\leq R$\\\cline{3-6}
		&&AROCK \small{\citep{ASCD2}} & $\tau n  L_c /\epsilon  $ & $\tau n \frac{L_c}{\mu}\log(1/\epsilon)  $ &  B\\\cline{3-6}
		&&AASCD \small{(ours)} &$ (n\sqrt{L}+n\tau\sqrt{L_c})\sqrt{1/\epsilon }$ & $  n\tau\sqrt{\frac{L_c}{\mu}}  \log(1/\epsilon)$& B \\\hline\hline
		\multirow{8}*{\eqref{problem}}	&\multirow{1}*{S}&Katyusha \small{\citep{kat}} &$n+(n+\sqrt{n})\sqrt{L/\epsilon }$& $n+\sqrt{\frac{nL}{\mu}}  \log(1/\epsilon)$ &  \\\cline{2-6}
		&\multirow{7}*{A}&Hogwild \small{\citep{HOGWILD}}& Not analysis & $\frac{\tau^2 \triangle L}{\epsilon} $ & B, Sparse, Smooth \\\cline{3-6}
		&&ASGD \small{\citep{Duch}}&$\frac{\tau L}{\epsilon} +\frac{L}{\epsilon^2\tau^2}$ & Not analysis & B, Smooth \\\cline{3-6}
		&&ASVRG \small{\citep{asvrg}}&  Not analysis & $ n+(1+\triangle \tau^2) \frac{nL}{\mu}\log(1/\epsilon) $ &B, Sparse, Smooth  \\\cline{3-6}
		&&AASGD \small{\citep{meng2016asynchronous}}&   Not analysis  & $n+\tau\frac{nL}{\mu}\log(1/\epsilon) $ & B, Sparsity, Smooth\\\cline{3-6}
		&	&ASVRG \small{\citep{asvrg3}}&  Not analysis  & $n+\tau^2 \frac{nL}{\mu}\log(1/\epsilon) $ & B \\\cline{3-6}
		&	&ASVRG \small{(ours)}&  $n+\tau^2Ln/\epsilon$  &  Not analysis & B \\\cline{3-6}
		&&AASVRG \small{(ours)} &$n+(n+\tau\sqrt{n})\sqrt{L/\epsilon }$ & $n+\tau \sqrt{\frac{nL}{\mu}}  \log(1/\epsilon)$  & B\\\hline
	\end{tabular}
	\label{convergence}
\end{table}

Up to now, lots of plain asynchronous algorithms are designed. For example, \citet{HOGWILD} and \citet{Duch} propose Asynchronous Stochastic Gradient (ASGD), which achieves $O(1/\epsilon)$ convergence rate for strongly convex  (SC) functions, where $\epsilon$ is the approximate error satisfying $F(\x)-F(\x^*)\leq \epsilon$. Some Variance Reduction (VR) based asynchronous algorithms \citep{asvrg,asvrg1,asvrg2} are also designed later. For Asynchronous Stochastic Coordinate Descent (ASCD)~\citep{ASCD1, ASCD2}, the provable convergence rate is $O(1/\epsilon)$ for non-strongly convex  (NC) functions and $\kappa \log (1/\epsilon)$ for SC, where $\kappa $ is the condition number.  A more detailed comparison for convergence results\footnote{ To the best of our knowledge, there is still no analysis on asynchronous VR algorithms for NC. As a byproduct of our analysis, the convergence rate is $O(n+\tau^2Ln/\epsilon)$, as shown in Table ~\ref{convergence}. The proof is shown in Supplementary Material.} of  asynchronous algorithms for convex problems is shown in Table~\ref{convergence}.

On the other hand, \citet{nesterov1983method,nesterov1988approach} has proposed a well-known accelerated version of gradient descent (AGD) for $L$-smooth convex functions.  AGD achieves $O(1/\sqrt{\epsilon})$ rate for NC and $O(\sqrt{\kappa}\log(1/\epsilon))$ for SC, which provably meets the lower bound (ignoring the constant) and is also observed to be faster than existing ones. After that, many accelerated algorithms have been designed to achieve faster convergence rates. For example, FISTA~\citep{beck2009fast} is a proximal version of AGD. APCG~\citep{fercoq2015accelerated, APCG} is a proximal and accelerated version of Stochastic Coordinate Descent~(SCD). Acc-SDCA~\citep{shalev2014accelerated} uses the black-box technique to accelerate the Stochastic Dual Coordinate Ascent. Katyusha~\citep{kat} is an accelerated version of VR methods.

\begin{table*}[t]
	\centering
	\caption{Notations for different algorithms in this paper}
	\small
	\begin{tabular} {|c |c |c|c|}\hline
		
		\multirow{2}*{AGD} & objective function& $f(\x)+h(\x)$&  $f(\x)$: $L$-Lipschitz continuous gradient\\\cline{2-4}
		& superscript $k$  and $j(k)$& $\x^k, \x^{j(k)}$&  $\x$ in $k$-th, $j(k)$-th iteration, respectively\\\cline{1-4}
		\multirow{5}*{AASCD} & \multirow{2}*{objective function and $n$} &  \multirow{2}*{$f(\x)+h(\x)$} & $f(\x)$: Lipschitz coordinate continuous gradient;  \\
		& & &  $h(\x) = \sum_{i =1}^n h_i(\x_i)$,  $\x \in \mathcal{R}^n$.\\\cline{2-4}
		& superscript $k$  and $j(k)$& $\x^k, \x^{j(k)}$&  $\x$ in $k$-th$, j(k)$-th iteration, respectively,\\\cline{2-4}
		& subscript $i $& $\x_i, \nabla_i f (\x)$&  the $i$-th coordinate of $\x$ and $\nabla f(\x)$, respecitvely.\\\cline{2-4}
		&  $i_k$ &\multicolumn{2}{c|}{The  index of the coordinate randomly chosen at iteration $k$.}\\\cline{1-4}
		\multirow{5}*{AASVRG} & \multirow{2}*{objective function and $n$} &  \multirow{2}*{$f(\x)+h(\x)$} & $f(\x)=\frac{1}{n}\sum_{i=1}^n f_i(\x)$ ;  \\
		& & &  with $f_i$: $L$-Lipschitz  continuous gradient, \\\cline{2-4}
		& superscript $s$, subscript $k$, $j(k)$& $\x^s_k, \x^s_{j(k)}$&  $\x$ in $k$-th, $j(k)$-th iteration, at $s$-th epoch, respectively,\\\cline{2-4}
		&  $i^s_k$ &\multicolumn{2}{c|}{The random index of the function chosen at iteration $k$, epoch $s$.}\\\cline{2-4}
		& $\tx$, $\tna$ &\multicolumn{2}{c|}{  snapshot vector and VR gradient followed by \citep{SVRG} .}\\\cline{1-4}
		
	\end{tabular}
	\label{notation}
\end{table*}

Comparing those plain asynchronous algorithms with serial ones, there is a gap in convergence rate. It is an open problem to fill in the gap by proposing  accelerated asynchronous algorithms. We find that \citet{meng2016asynchronous} integrates momentum, VR, tricks, coordinate sampling to accelerate ASGD, named AASGD.  But the convergence rate is still $O(\kappa \log(1/\epsilon))$ for SC functions. There is no improvement in convergence rate comparing with ASVRG.  Designing an asynchronous accelerated algorithms is not easy. The reason are two-folded:
\begin{itemize}
	\item  In serial accelerated schemes, the extrapolation point are subtly and  strictly connected with $\x^k$ and $\x^{k-1}$,  i.e. $\y^k = \x^k +\frac{\theta^k(1-\theta^{k-1})}{\theta^{k-1}}\left( \x^{k}-\x^{k-1}\right)$. However, such information might not be available for asynchronous algorithms because there are unknown delays in updating the parameters.
	\item  Since $\x^{k+1}$ is updated based on $\y^k$, i.e. $\x^{k+1} = \y^k - \frac{1}{L}\nabla f(\y^{k})$,  $\x^{k+1}$ is related to  past updates (to generate $\y^k$). This is different from unaccelerated algorithms. For example, in gradient descent, $\x^{k+1} = \x^k - \frac{1}{L}\nabla f(\x^{k})$, so $\x^{k+1}$ only depends on $\x^k$.
\end{itemize}
In this paper, we attempt to fill in the gap to some degree.  We propose a technique called ``momentum compensation'' to accelerate asynchronous algorithms for convex problems. We first consider accelerating plain Asynchronous Gradient Descent.  We demonstrate that doing only one original step of momentum prevents us from bounding  the distance between delayed gradient and the latest one.  Instead, by ``momentum compensation'' our algorithm is able to achieve a faster rate, i.e. $O(\tau^2/\sqrt{\epsilon})$ for NC functions and  $O(\tau^2\sqrt{\kappa}\log(1/\epsilon))$ for SC ones, where $\tau$ is the upper bound of delay. We then show that this technique can be further applied to modern stochastic algorithms, by designing Accelerated Asynchronous Stochastic Coordinate Descent (AASCD) and Asynchronous Stochastic Gradient Descent (AASVRG). Both of the resultant algorithms are faster than existing ones by order  and even with less order of $\tau$ comparing with  AAGD.   We also  show that for sparse datasets, the delay  will be largely reduced and  linear speed up is achievable for our algorithms under certain conditions.
Finally, we conduct lots of experiments on a shared memory system to demonstrate the fast convergence of our algorithms. To summarize, we list the contributions of our work as follows:
\begin{itemize}
	\item[1.] We propose   the `` momentum compensation'' technique to accelerate asynchronous algorithms for convex problems.  To the best of our knowledge, we are the first to consider accelerated algorithms for delayed gradients. Our results are strong (improve the rate), general (includes analysis for proximal version and NC), and easy to combine with other techniques (see 2).
	\item[2.] We show that our technique can be applied to modern stochastic asynchronous algorithms.  The resultant algorithms, i.e. AACD and AASVRG, are also faster than existing ones by order.
	\item[3.] We perform lots of numerical experiments on a shared memory system to demonstrate that acceleration can lead to significant performance improvements. We will put our C++ implementation with POSIX threads on website once our paper is accepted.
\end{itemize}

\section{Preliminaries and Notations }\label{prelimi}
In most asynchronous parallelism, there are typically two schemes:
\begin{itemize}
	\item Atom~(consistent read) scheme: The parameter $\x$ is updated as an atom. When $\x$ is read or updated in the central node, it will be locked. So $\x^{j(k)}\in \{\x^0,\x^1, \cdots,\x^k \}$.
	\item Wild~(inconsistent read) scheme: To further reduce the system overhead, there is no lock in implementation.  All the threads may perform modifications on $\x$ at the same time~\citep{HOGWILD}. Obviously, analysis becomes more complicated in this situation.
\end{itemize}
In this work, our analysis focuses on the atom scheme. However, we implement our algorithm in the wild scheme. We leave the analysis of the wild scheme as a future work. For more details, please refer to \citep{HOGWILD,ASGD2}.

Since the gradients are delayed for asynchronous algorithms, lots of algorithms are not guaranteed to converge without additional assumption on delay, let alone acceleration.
In this paper, we follow \citep{HOGWILD, asvrg, ASCD2} to assume a bounded delay.
\begin{assumption}\label{fir}
	We assume that all the updates before $(k -\tau -1)$-th iteration are completed before the ``read'' step of the $k-$th iteration. So for the atom scheme, we have
	\begin{eqnarray}\label{123123}
	{j(k)}\in \{k-\tau, k-\tau+1, \cdots,k \}.
	\end{eqnarray}
\end{assumption}
The parameter $\tau$ has expressed the degree of delay. When there are more threads, the delay accumulates and results in larger $\tau$.

The notations for different algorithms are shown in Table~\ref{notation}. The three algorithms are independent without confusion, so by a little abuse of notation we make it easier for understanding our technique.  For all algorithms, we use $j(k)$ to denote the delayed state which satisfies Eq.~\eqref{123123}.  We say the function $f$ has $L$-Lipschitz continuous gradient ($L$-smooth) if
\begin{eqnarray}\label{l1}
\| \nabla f(\x) -\nabla f(\y)\|\leq L\|\y-\x\|.
\end{eqnarray}
For AASCD, we say the function $f$ has $L_c$ coordinate continuous gradient, where $L_c = \max_{i=1}^n {L_i}$, and $L_i$ is the coordinate Lipschitz constant for $\nabla f$ in the $i$-th coordinate direction:
\begin{eqnarray}\label{lc}
|\nabla_i f(\x) -  \nabla_i f(\y)  |\leq  L_i|\x_i -\y_i  |,\quad 1\leq i\leq n,
\end{eqnarray}
in which $\x_i$ and $\y_i$ denote the $i$-th coordinate of $\x$ and $\y$, respectively. $\nabla_i$ denotes the $i$-th partial gradient as shown in Table~\ref{notation}.

\begin{algorithm}[tb]
	\caption{AGD~\citep{nesterov1983method}}
	\label{AGD-ori}
	\begin{algorithmic}
		\STATE  $\mathbf{Input}$ $\theta^k$, step size $\gamma$, $\x^0= \mathbf{0} $, and $\z^0= \mathbf{0}$ .
		\STATE  $  \mathbf{for}$ $k=0$ $\mathbf{to}$ $K$ $\mathbf{do}$\\
		\STATE1  $\quad \y^k= (1-\theta^k)\x^k+\theta^k\z^k$.
		\STATE2 $\quad \bdelta^k = \argmin_{\bdelta} h(\z^k+\bdelta) +\langle \nabla f(\y^k), \bdelta \rangle +\frac{\theta^k }{2\gamma}\| \bdelta\|^2  $.
		\STATE3  $\quad \z^{k+1} =  \z^k +\bdelta^k  $.
		\STATE4  $\quad  \x^{k+1} =\theta^k \z^{k+1} +(1-\theta^k)\x^k  $.
		\STATE $\mathbf{end \ for}$
		\STATE  $\mathbf{Output}$ $\x^{K+1}$.
	\end{algorithmic}
\end{algorithm}

\section{AAGD}
We first illustrate  our momentum  compensation technique for plain AGD algorithms. The objective function is:
\begin{eqnarray}\label{pro1}
\min_{\x} f(\x)+h(\x),
\end{eqnarray}
where $f(\x)$ has $L$-Lipschitz continuous gradient and both $f(\x)$ and $h(\x)$ are convex.

\subsection{Momentum Compensation}
Recall the  serial  Accelerated Gradient Descent~\citep{nesterov1983method},  shown in Algorithm~\ref{AGD-ori}. If we directly implement AGD~\citep{nesterov1983method} asynchronously, we can only get the gradient $ \nabla f(\y^{j(k)})$ at Step 1 due to the delay. Now we need  to measure the distance between $\y^{j(k)}$ and $\y^{k}$. With some algebraic transformation, we have
\begin{eqnarray}\label{111}
\y^{k+1} = \x^k + \frac{\theta^k(1-\theta^k)}{\theta^k} (\x^k -\x^{k-1}),
\end{eqnarray}
which is known as extrapolation. Set $a^k =  \frac{\theta^k(1-\theta^k)}{\theta^k}$, and $b(l,k) =\prod_{i=l}^{k} a^{i}$, where $l\leq k $. Then by applying  Eq.~\eqref{111} recursively,  for $k\geq j(k)\geq 0$, we have,
\begin{eqnarray}\label{one}
\y^k&=&\y^{k-1}+\sum_{i=j(k)+1}^{k}\left( 1+b(i,k)(\x^i-\y^{i-1})\right) \notag\\
&&+b(j(k),k)(\x^{j(k)}-\x^{j(k)-1}).
\end{eqnarray}
Summing Eq.~\eqref{one} with superscript from $j(k)$ to $k$, we obtain the relation between $\y^{j(k)}$ and $\y^{k}$:
\begin{eqnarray}
\y^k&=& \y^{j(k)}+\sum_{i=j(k)+1}^{k}\left( 1+\sum_{l=i}^{k} b(i,l)\right)(\x^i-\y^{i-1})\notag\\
&&+\left(\sum_{i=j(k)+1}^{k}b(j(k),i)\right)(\x^{j(k)}-\x^{j(k)-1}).
\end{eqnarray}
We find that $\x^{j(k)}-\x^{j(k)-1} $ is related to all the past updates before $j(k)$. If we directly implement AGD asynchronously like most asynchronous algorithms,
then $j(k) <k$ (due to delay), so $\sum_{i=j(k)+1}^{k}b(j(k),i)> 0$.  Since $\x^{j(k)}-\x^{j(k)-1}$ is hard to bound,  it causes difficulty to obtain the accelerated convergence rate.

Instead, we compensate the momentum term and introduce a new extrapolation point $\w^{j(k)}$, such that
\begin{eqnarray}\label{ww}
\w^{j(k)} =  \x^{j(k)}\!\!+\! \left(\! \sum_{i=j(k)}^k b(j(k),i)  \!\right)\!(\x^{j(k)} -\x^{j(k)-1}).
\end{eqnarray}
One can find these are actually several steps of momentum. Then the difference between
$\y^{k}$ and $\w^{j(k)}$ can be directly bounded by the norm of  several latest updates, namely $\| \sum_{i=j(k)+1}^{k}\left( 1+\sum_{l=i}^{k} b(i,l)\right)(\x^i-\y^{i-1}) \|^2$. So we are able to obtain the accelerated rate. The Algorithm is shown in Algorithm \ref{AAGD}.
\begin{algorithm}[t]
	\caption{AAGD}
	\label{AAGD}
	\begin{algorithmic}
		\STATE  $\mathbf{Input}$ $\theta^k$, step size $\gamma$, $\x^0= \mathbf{0} $ and $\z^0= \mathbf{0}$.
		\STATE  $  \mathbf{for}$ $k=0$ $\mathbf{to}$ $K$ $\mathbf{do}$\\
		\STATE1  $\ \ \w^{j(k)} =  \x^{j(k)}\!\!+\! \left(\! \sum_{i=j(k)}^k b(j(k),i)  \!\right)\!(\x^{j(k)} -\x^{j(k)-1})$.
		\STATE2 $\ \ \bdelta^k = \argmin_{\bdelta} h(\z^k+\bdelta) +\langle \nabla f(\w^{j(k)}), \bdelta \rangle +\frac{\theta^k }{2\gamma}\| \bdelta\|^2  $.
		\STATE3  $\ \ \z^{k+1} =  \z^k +\bdelta^k  $.
		\STATE4  $\ \  \x^{k+1} =\theta^k \z^{k+1} +(1-\theta^k)\x^k  $.
		\STATE $\mathbf{end \ for}$
		\STATE  $\mathbf{Output}$ $\x^{K+1}$.
	\end{algorithmic}
\end{algorithm}

\subsection{Convergence Results}

After introducing $\w^{j(k)}$, we separately analysis $f(\x^{k+1})-f(\x^*)$  and $\|\z^{k+1}-\z^*\|^2$ like the Lyapunov technique~\citep{asvrg,asvrg1}, and bound them through the existing terms in serial AGD~\citep{nesterov1983method} and additional $\|\w^{j(k)}-\y^k\|^2$. Then we  choose a proper step size to obtain a faster convergence rate. We directly give the convergence results of AAGD. All the proofs can be found in Supplementary Material.

\begin{algorithm}[tb]
	\caption{AAGD-implementation}
	\label{AAGD-implementation}
	\begin{algorithmic}
		\STATE  $\mathbf{Input}$  $\theta^k$, step size $\gamma$, $\uu^0= \mathbf{0} $, $\vv^0= \mathbf{0}$, and $d^0 =1$.
		\STATE  $  \mathbf{for}$ $k=0$ $\mathbf{to}$ $K$ $\mathbf{do}$\\
		\STATE1  $\quad  \!\!\!d^{k+1}= d^{k}(1-\theta^k)$.
		\STATE2  $\quad\!\!\!\w^{j(k)} =  \uu^{j(k)} + d^{k+1} \vv^{j(k)} .$
		\STATE3 $\quad  \!\!\! \bdelta^k = \argmin_{\bdelta} h(\z^k+\bdelta) +\langle \nabla f(\w^{j(k)}), \bdelta \rangle +\frac{\theta^k }{2\gamma}\| \bdelta\|^2  $.
		\STATE4  $\quad\!\!\! \uu^{k+1} =  \uu^{k}+ \bdelta^k  $.
		\STATE5  $\quad \!\!\!\vv^{k+1} =  \vv^{k}- \frac{\bdelta^k}{d^k} $.
		\STATE $\mathbf{end \ for}$
		\STATE  $\mathbf{Output}$ $\x^{K+1}=\uu^{K+1} +  d^{K+1} \vv^{K+1}  $.
	\end{algorithmic}
\end{algorithm}

\begin{theorem}\label{th1}
	Under Assumption 1, for Algorithm~\ref{AAGD}, for for
	non-strongly convex case, if the step size satisfies $ 2\gamma L  + 3\gamma L(\tau^2+3\tau)^2\leq 1  $, $\theta^k=\frac{2}{k+2}$, and the first $\tau$ iterations are updated in serial\footnote{We use this assumption only for simplicity. This assumption is removed in the analysis of AASCD and AASVRG.}, we have
	\begin{eqnarray}
	F(\x^{K+1}) -F(\x^*) \leq  (\theta^k)^2 \left(\frac{1}{2\gamma} \|  \z^{0} - \x^*\|^2\right).
	\end{eqnarray}
	When $h(\x)$ is strongly convex with modulus $\mu \leq L$, the step size satisfies $ \frac{5}{2}\gamma L  + \gamma L(\tau^2+3\tau)^2\leq 1$, and $\theta^k = \frac{-\gamma\mu+\sqrt{\gamma\mu^2+4\gamma\mu}}{2}$ is denoted as $\theta$ instead, we have
	\begin{eqnarray}
	F(\x^{K+1}) -F(\x^*) \leq (1-\theta)^{K+1}\left( F(\x^{0}) -F(\x^*)\right) \notag\\
	+  (1-\theta)^{K+1}\left(\left(\frac{\theta^2}{2\gamma}+\frac{\mu\theta}{2} \right)  \|  \z^{0} - \x^*\|^2\right).
	\end{eqnarray}
\end{theorem}

\begin{corollary}
	For Algorithm~\ref{AAGD}, under the assumption of Theorem \ref{th1}, the Iteration First-Order (IFO) calls are $O(\tau^2 \sqrt{L/\epsilon})$ for NC and $O(\tau^2 \sqrt{L/\mu}\log(1/\epsilon))$ for SC.
\end{corollary}
The order of $\tau$ is large for AAGD, we will show that for stochastic  asynchronous algorithms, the order of $\tau$ will be largely reduced.
\subsection{AAGD in Implementation}
In Eq.~\eqref{ww}, we need to compute $\sum_{i=j(k)}^k b(j(k),i) $, which is a little complicated. To make our algorithm clearer, inspired by \citep{fercoq2015accelerated, APCG}, we can change variable as follows:  $\z^k= \uu^k$, $\x^k = \uu^k +a^k\vv^k$, and $\y^k = \uu^k +a^{k+1}\vv^k$. The algorithm is shown in Algorithm~\ref{AAGD-implementation}. The equivalent of Algorithm \ref{AAGD}  and \ref{AAGD-implementation}  is shown Supplementary Material. Another advantage for Algorithm~\ref{AAGD-implementation} is the ability to  sparse update for the sparse dataset.

\section{Practical  Asynchronous Algorithms}
\begin{algorithm}[t]
	\caption{AASCD}
	\label{AASCD-ori}
	\begin{algorithmic}
		\STATE  $\mathbf{Input}$ $\theta^k$, step size $\gamma$, $\x^0= \mathbf{0} $ and $\z^0= \mathbf{0}$.
		\STATE 	Define $a^k=\frac{\theta^{k}(1-\theta^k)}{\theta^{k-1}} $, $b(l,k)=\prod_{i=l}^{k} a^{l}$.
		\STATE  $  \mathbf{for}$ $k=0$ $\mathbf{to}$ $K$ $\mathbf{do}$\\
		\STATE1  $\quad\w^{j(k)} =  \y^{j(k)} + \sum_{i=j(k)+1}^k b(i,k)   (\y^{j(k)} -\x^{j(k)-1}).$
		\STATE2 $\quad$Randomly choose an index  $i_k$ form $[1,2,\cdots,n]$.
		\STATE3 $\quad \bdelta^k \!=\! \argmin_{\bdelta} h_{i_k}(\z^k+\bdelta) +\langle \nabla_{i_k} f(\w^{j(k)}), \bdelta \rangle +\frac{\theta^k }{2\gamma}\| \bdelta\|^2.$
		\STATE4  $\quad \z_{i_k}^{k+1} =  \z_{i_k}^k +\bdelta^k  $ with other coordinates unchanged.
		\STATE5 	$\quad  \y^{k} = (1-\theta^k)\x^k +\theta^k\z^k $.
		\STATE6  $\quad  \x^{k+1} = (1-\theta^k)\x^k +n\theta^k\z^{k+1}-(n-1)\theta^k\z^k  $.
		\STATE $\mathbf{end \ for}$
		\STATE  $\mathbf{Output}$ $\x^{K+1}$.
	\end{algorithmic}
\end{algorithm}


To meet the large-scale of machine learning,  most asynchronous algorithms are designed in a stochastic fashion. We are now to demonstrate that our technique can further be applied to accelerate modern state-of-the-art stochastic asynchronous algorithms, such as ASCD~\citep{ASCD1} and ASVRG~\citep{asvrg, asvrg3}.  The proofs of our AASCD and AASVRG are similar to that of AAGD, but are much  involved.  It needs to  further fuse other techniques, such as  Estimate Sequence  technique in~\citep{fercoq2015accelerated} for AASCD and the  negative momentum technique~\citep{kat} for AASVRG. Like AAGD, the two algorithms also be changed variables to  be clearer and able to sparse update. We directly demonstrate the algorithms and the convergence results. All the proofs can also be found in Supplementary Material.
\subsection{AASCD}
(Asynchronous) Stochastic Coordinate Descent algorithms mainly  solves the following problem:
\begin{eqnarray}\label{AASCD1}
\min_{\x\in \mathcal{R}^n} f(\x)+h(\x),
\end{eqnarray}
where $f(\x)$ has $L_c$-Lipschitz coordinate continuous gradient, $h(\x)$ has coordinate separable structure, i.e. $h(\x) = \sum_{i=1}^n h_i(\x_i)$, and  $f(\x)$ and $h_i(\x)$ are convex.
At each iteration, the algorithms choose one coordinate $\x_i$ to sufficiently reduce the objective value while keeping other coordinates fixed which reduces the per-iteration cost. In more detail, in each iteration the following types of proximal subproblem is solved:
\begin{eqnarray}
\bdelta^k  = \argmin_{\bdelta} h_{i_k}(\x^k+\bdelta) +\langle \nabla_{i_k} f(\x^k), \bdelta \rangle +\frac{\theta^k }{2\gamma}\| \bdelta\|^2,
\end{eqnarray}
where $\nabla_{i_k} f(\x)$ denotes the partial gradient of $f$ with respect to $\x_i$.

For asynchronous algorithms, the partial gradient will be delayed, and at iteration $k$ we could only obtain $\nabla_{i_k} f(\x^{j(k)})$ instead of $\nabla_{i_k} f(\x^{k})$.

Now we propose our accelerated algorithm. For simplicity, we assume that each coordinate Lipschitz constant  $L_i$ are the same,  then $L_c=L_i$, $i = 1,2,\cdots,n$\footnote{This is the case that the data are normalized. When $L_i$ are different,
	$ n^2L_c$ can be extended to $(n\sum_{i=1}^n(L_i)^2)$, also we can  fuse  the non-uniform sampling \citep{allen2016even} technique and replace it with   smaller $(\sum_{i=1}^n(\sqrt{L_i})^2$ in convergence rate.}.  By  judging the distance between the delayed extrapolation points  and the newest noes and compensating the  ''lost'' momentum term, we obtain Algorithm~\ref{AASCD-ori}.  We have the following theorem:
\begin{theorem}\label{th2}
	Under  Assumption 1, and $\tau\leq \sqrt{n}$, for Algorithm~\ref{AASCD-ori},  if the step size satisfies $2\gamma L_c  + (1+\frac{1}{n})\gamma L_c \left(\frac{\tau^2+\tau}{n} +2\tau \right)^2\leq 1$, and $\theta^k=\frac{2}{2n+k}$, we have
	\begin{eqnarray}
	&&\frac{\E [F(\x^{K+1})] -F(\x^*)}{(\theta^K)^2} +\frac{n^2}{2\gamma}\E\| \z^{K+1} -\x^* \|^2 \notag\\
	&\leq&\frac{F(\x^0)-F(\x^*)}{(\theta^{-1})^2} +\frac{n^2}{2\gamma}\|\z^0 -\x^*\|^2.
	\end{eqnarray}
	When $h(\x)$ is strongly convex with modulus $\mu \leq L_c$, the step size satisfies $2\gamma L_c  + (\frac{3}{4}+\frac{3}{8n})\gamma L_c\left((\tau^2+\tau)/n +2\tau\right)^2\leq 1$, and  $\theta^k = \frac{-\gamma\mu+\sqrt{\gamma\mu^2+4\gamma\mu}}{2n}$ is denoted as $\theta$ instead, we have
	\begin{eqnarray}
	&&\!\!\!\!\!\!\!\!\E [F(\x^{K+1})] -F(\x^*)\leq (1-\theta)^{K+1} \left(F(\x^0) -F(\x^*)\right)\notag\\
	&& +(1-\theta)^{K+1} \left(\frac{n^2(\theta)^2+n\theta\mu\gamma}{2\gamma}\| \z^{0} -\x^* \|^2 \right).
	\end{eqnarray}
\end{theorem}
\begin{corollary}
	For Algorithm~\ref{AASCD-ori}, under the assumption of Theorem \ref{th2}, the IFO calls are $ O\left((n\sqrt{L}+n\tau\sqrt{L_c})\sqrt{1/\epsilon }\right)$ for NC and $O\left(  n\tau\sqrt{\frac{L_c}{\mu}}  \log(1/\epsilon)\right)$ for SC.
\end{corollary}
We can find that the order of $\tau$ are reduced comparing with AAGD due to the stochastic effect.

\begin{algorithm}[tb]
	\caption{AASVRG}
	\label{AASVRG-ori}
	\begin{algorithmic}
		\STATE  $\mathbf{Input}$ $\theta^s_1$, step size $\gamma$, $\x^0_0= \mathbf{0}$, $\tx^0= \mathbf{0}$,  and $\z^0_0= \mathbf{0}$, $\theta_2 = \frac{1}{2}$,  $m=n$,  and $a^s=1-\theta_2 -\theta^s_{1} $.
		\STATE  $  \mathbf{for}$ $s=0$ $\mathbf{to}$ $S$ $\mathbf{do}$ \\
		\STATE  $ \ \ \ \mathbf{for}$ $k=0$ $\mathbf{to}$ $m-1$ $\mathbf{do}$$\quad\quad\circ$ start asynchronous update\\
		\STATE1  $\ \ \ \  \w^s_{j(k)} =  \x^s_{j(k)}+ \frac{a^s\left(1-(a^s)^{k-j(k)+1}\right)}{1-a^s} (\x^s_{j(k)} -\x^s_{j(k)-1}).$
		\STATE2 $\ \ \ \  $Randomly selected an sample with index $i^s_k$.
		\STATE3 $\ \ \ \  $$\tna^s_k =  \nabla f_{i^s_k}(\w^s_{j(k)})-\nabla f_{i^s_k}(\tx)+ \nabla f(\tx)$.
		\STATE4 $\ \ \ \  \bdelta^s_k = \argmin_{\bdelta} h(\z^s_k+\bdelta) +\langle \tna^s_k, \bdelta \rangle +\frac{\theta^s_1 }{2\gamma}\| \bdelta\|^2  $.
		\STATE5  $\ \ \ \  \z^s_{k+1} =  \z^s_k +\bdelta^s_k  $.
		\STATE6  $\ \ \ \  \x^s_{k+1} =\theta_1^s \z^s_{k+1} +\theta_2\tilde{\x} +a^s \x^s_k  $.
		\STATE $ \ \ \ \mathbf{end \ for}$ $k$.$\quad\quad\quad\quad\quad\quad\quad\circ$ synchronization
		\STATE $ \quad \x^{s+1}_{0} = \x^s_m$, $\z^{s+1}_{0} = \z^s_m$.
		\STATE $ \quad$For NC: $\tx^{s+1} =\frac{1}{m}\sum_{k=0}^{m-1}\x^{s}_k$,
		\STATE $ \quad$For SC:
		\STATE $ \quad$ $\tx^{s+1} =\left(\sum_{k=0}^{m-1} (1+\theta^s_1)^{i}\right)^{-1}\sum_{k=0}^{m-1} (1+\theta^s_1)^{i}\x^{s}_k$.
		\STATE $\mathbf{end \ for}$ $s$.		
		\STATE  $\mathbf{Output}$ $\x^{S+1}_0$.
	\end{algorithmic}
\end{algorithm}

\subsection{AASVRG}
We  consider the following  composite finite-sum convex optimization problem:
\begin{eqnarray}\label{problem}
F(\x) = h(\x) + \frac{1}{n}\sum_{i=1}^n f_i(\x).
\end{eqnarray}
where $f_i(\x)$'s , $i= 1,2,\cdots,n$, are convex and have Lipschitz continuous gradients, and $h(\x)$ is also convex. We denote $f(\x)=\frac{1}{n}\sum_{i=1}^n f_i(\x)$. To solve Eq.~\eqref{problem}, stochastic methods computes a gradient estimator from one or several $f_i(\x)$ to reduce the computation cost.   For asynchronous algorithms, the VR based asynchronous algorithms are proposed as state-of-the-art methods to solve Eq~\eqref{problem}. We show that our technique can further accelerate these algorithms. Like ASVRG, we adopt asynchronous update in each inner loop. There will be synchronization operation after each epoch.  Since $n$ is large, the cost for synchronization is small comparing with the cost for computation. The algorithm is shown in Algorithm~\ref{AASVRG-ori}.  We have the following theorem:
\begin{theorem}\label{th3}
	Under the Assumption 1, for Algorithm~\ref{AASCD-ori}, if the step size satisfies $ 5\gamma L +10 \gamma L\tau^2 \leq 1$, and $\theta^s_1=\frac{2}{s+4}$, we have
	\begin{eqnarray}
	&& \!\!\!\!\!\!\!\!\!\!\!\!\!\!\!\!\E\left( F(\x^S_{n})-F(\x^*)\right)+(\theta_{2}+\theta^S_1)\sum_{k=1}^{n-1}\E\left(F(\x^S_{k})-F(\x^*) \right)\notag\\
	&\leq&2n (\theta^S_1)^2 (F(\x^0_0)-F(\x^*))+\frac{(\theta^S_1)^2}{2\gamma}\|\z^0_0-\x^*\|^2.
	\end{eqnarray}
	When $h(\x)$ is strongly convex with modulus $\mu \leq \frac{L\tau^2}{4n}$, the step size satisfies $ 5\gamma L+\frac{95}{8}\tau^2\gamma L\leq 1$, $\theta^s_1= \frac{1}{\tau}\sqrt{\frac{n\mu}{L}}$, and $\theta_3 = 1+\frac{\mu\gamma}{\theta_1^s}$, we have
	\begin{eqnarray}
	&&\!\!\!\!\!\!\!\!\!\!\!\left(  F(\tx^{S+1})-F(\x^*)   \right)\leq (\theta_3)^{-Sn}\left(\frac{\gamma}{4n}\|\z^0_0-\x^* \|^2 \right)\notag\\
	&& +  (\theta_3)^{-Sn} \left(    (1+\frac{1}{n}) \left(F(\x^0_0)-F(\x^*) \right) \right).
	\end{eqnarray}
\end{theorem}
\begin{corollary}
	For Algorithm~\ref{AASCD-ori}, under the assumption of Theorem \ref{th3}, the IFO calls are $O\left(n+(n+\tau\sqrt{n})\sqrt{L/\epsilon }\right)$ for NC and $O\left(n+\tau \sqrt{\frac{nL}{\mu}}  \log(1/\epsilon)\right)$ for SC.
\end{corollary}
We can find that the order of $\tau$ is also lower than AAGD.

\section{Applications}
We focus on solving  Empirical Risk Minimization problems:
\begin{eqnarray}\label{13}
\min_{\x \in \mathcal{R}^d} P(\x) = \frac{1}{n} \sum_{i = 1}^n \phi_i(\A^T_i \x) +\lambda g(\x),
\end{eqnarray}
where $\lambda >0$, $g(x)$ is typical a regular terms, and $\sum_{i = 1}^n \phi_i(\A^T_i \x)$  are loss functions over training samples. Lots of machine learning problem can be formulated into Eq~\eqref{13}, such as linear SVM, Ridge Regression, and Logistic Regression. For AASVRG,  solving Eq.\eqref{13} is equivalent to Eq.~\eqref{problem}.
For AASCD, we consider solve Eq.~\eqref{13} in dual.
When $g(\x) = \|\x\|^2$. 	 The dual formulation of  Eq.~\eqref{13} is :
\begin{eqnarray}\label{dual}
\min_{\a \in \mathcal{R}^n} D(\a) = \frac{1}{n} \sum_{i = 1}^n \phi^*_i(-\a_i) + \frac{\lambda}{2}\|\frac{1}{\lambda n}\A\a\|^2.
\end{eqnarray}
Through the technique of \citep{APCG}, for SC, we can obtain the convergence rate for AASCD on primal:
\begin{theorem}
	Assume that each function $\phi_i$ is $L_2$-smooth, $g(\cdot)$ has a unit convexity modulus $1$, and $\| \A_i\|\leq R$, for all $i = 1,\cdots, n$.  Then the IFO calls to reach both the dual  optimality gap ($\E [D^*-D(\a^k)]\leq  \epsilon$) and the primal one ($\E [P(\x)-P^*]\leq  \epsilon$)   through Algorithm~\ref{AASCD-ori} are $O\left(n+\tau \sqrt{\frac{nR^2 L_2}{\lambda}}\log(1/\epsilon)\right)$.	
\end{theorem}

\begin{table}
	\centering
	\caption{Details of sparse datasets.}
	\begin{tabular}{crrr}
		\toprule
		Datasets   & \#samples  & \#features & \#nonzeros \\
		\midrule
		real-sim   & 72,309  & 20,958        & 3,709,083  \\
		news20     & 19,996  & 1,355,191     & 9,097,916  \\
		rcv1       & 20,242  & 47,236        & 49,556,258  \\
		url        & 2,396,130	& 3,231,961  & 277,058,644 \\			
		\bottomrule
	\end{tabular}
	\label{tab:dataset}
\end{table}

\begin{figure*}[t]
	\center
	\hspace{-0.15in}
	\subfigure[rcv1\_ASCD]{\includegraphics[width=1.55in]{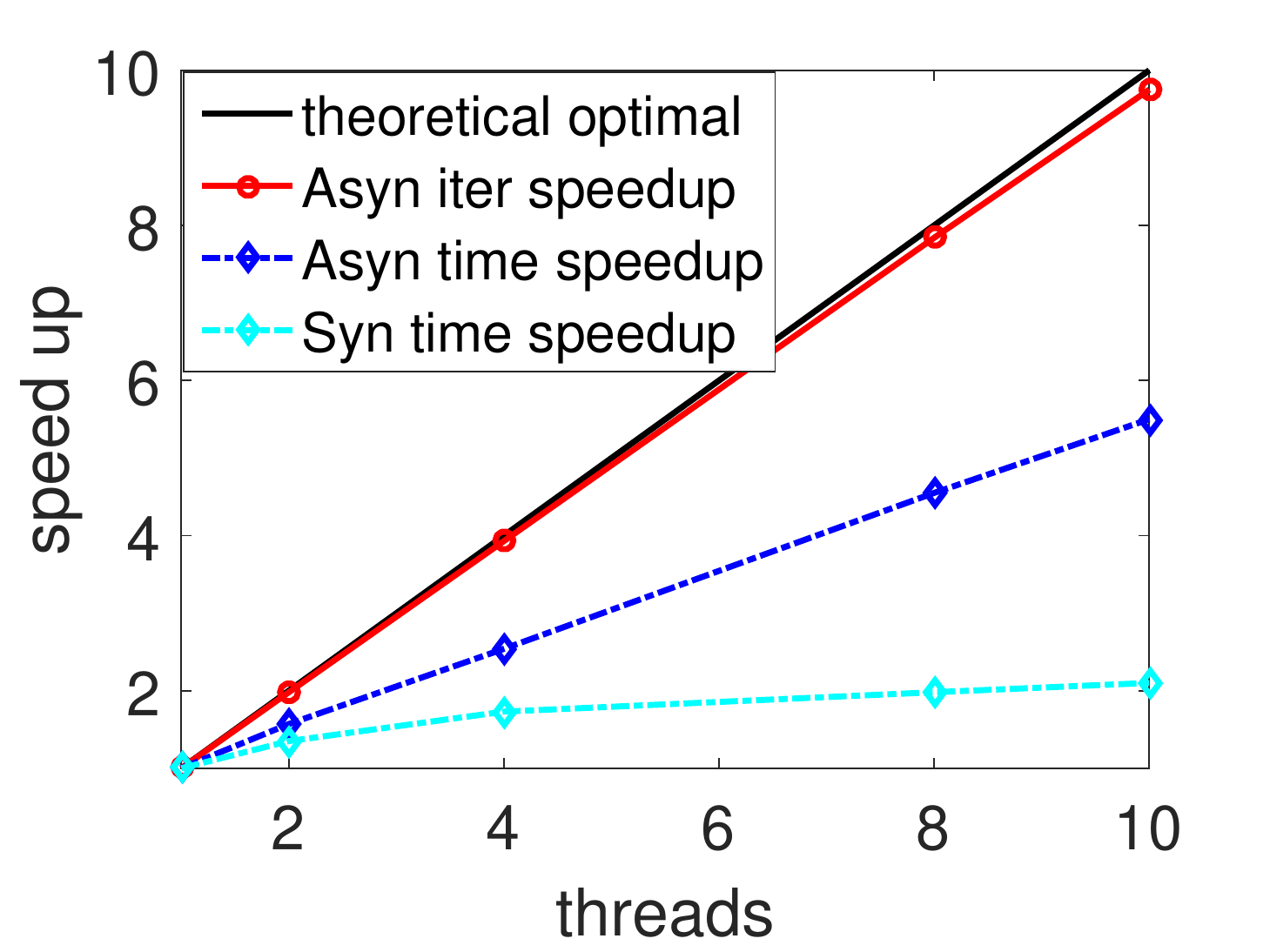}}
	\hspace{-0.15in}
	\subfigure[rcv1\_AASVRG]{\includegraphics[width=1.55in]{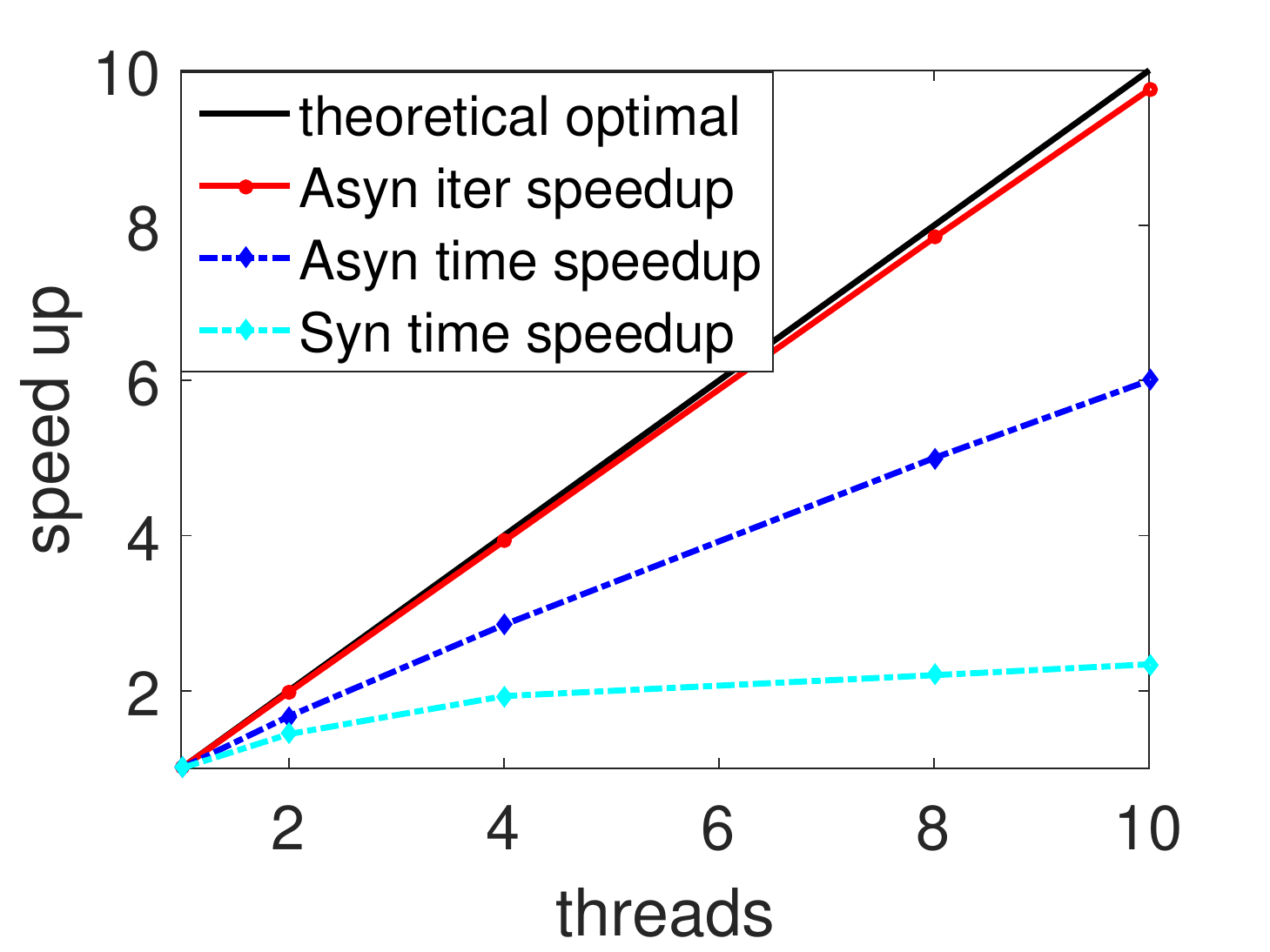}}
	\hspace{-0.15in}
	\subfigure[real-sim\_ASCD]{\includegraphics[width=1.55in]{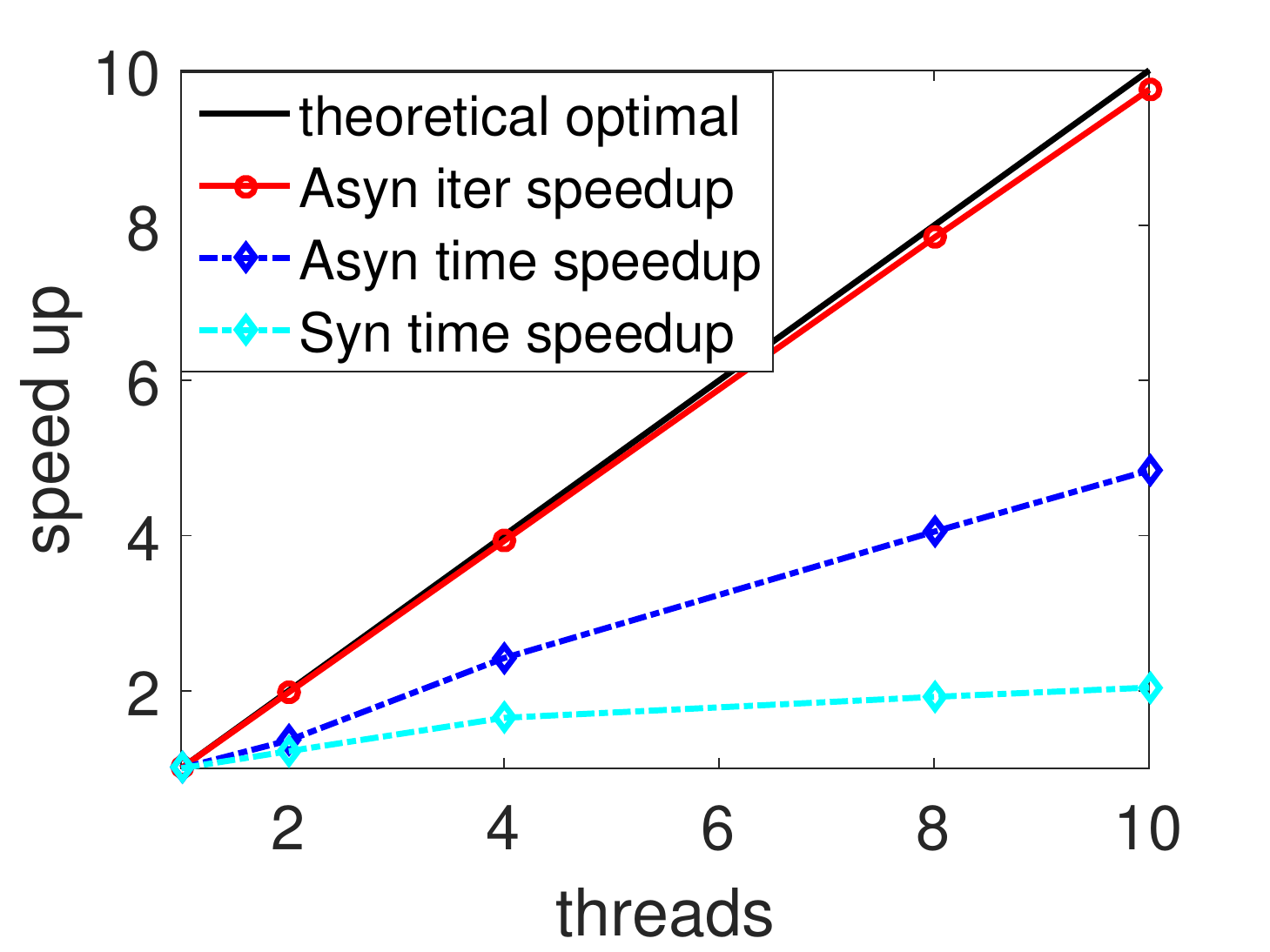}}
	\hspace{-0.15in}
	\subfigure[real-sim\_AASVRG]{\includegraphics[width=1.55in]{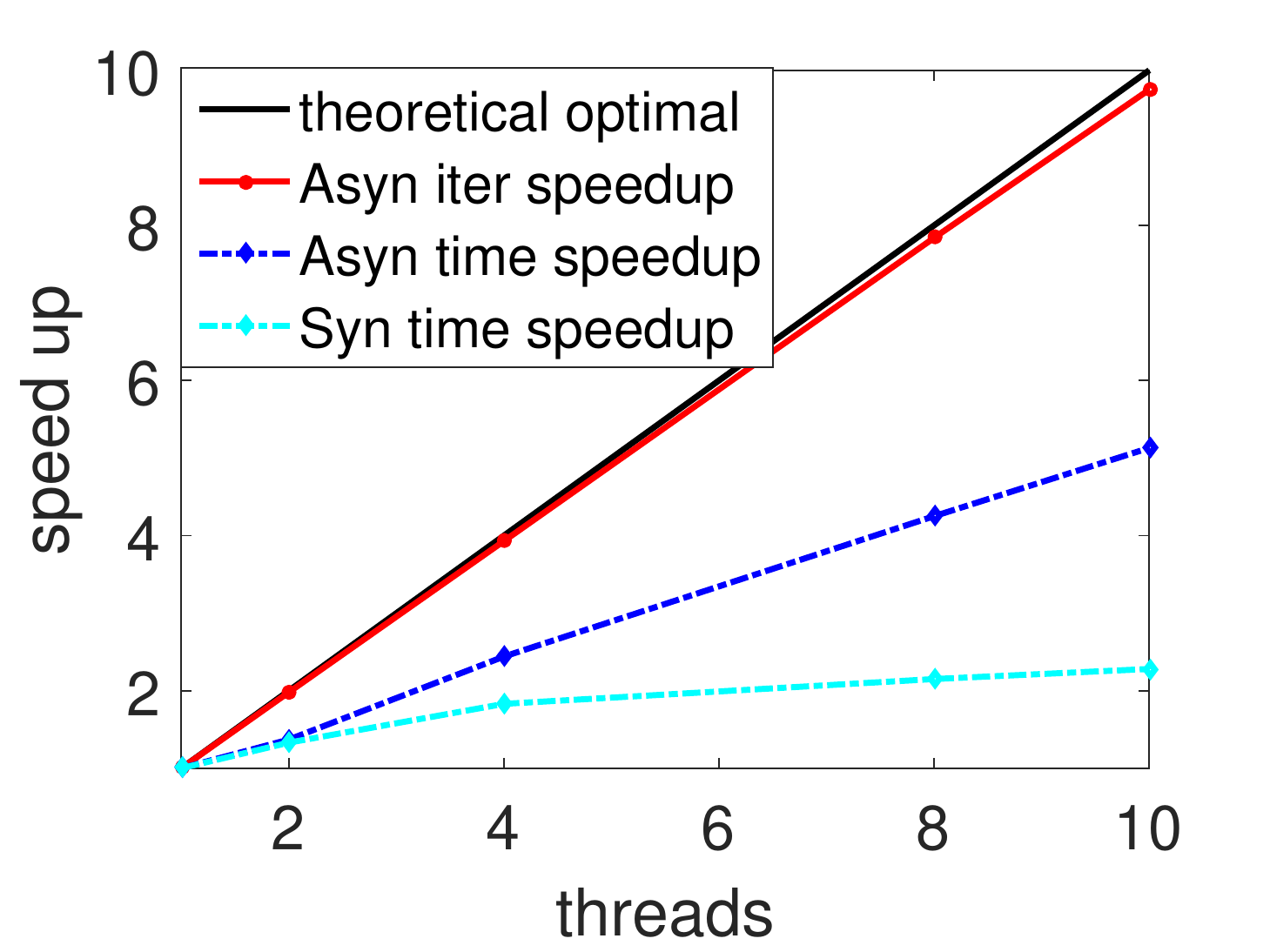}}
	\vspace{-0.3cm}	
	\caption{Experimental results on speed up for asynchronous implementation versus synchronous implementation.}
	\label{fig:speedup}
	\vspace{-0.3cm}
\end{figure*}
\begin{figure*}[t]
	\center
	\hspace{-0.15in}
	\subfigure[rcv1]{\includegraphics[width=1.55in]{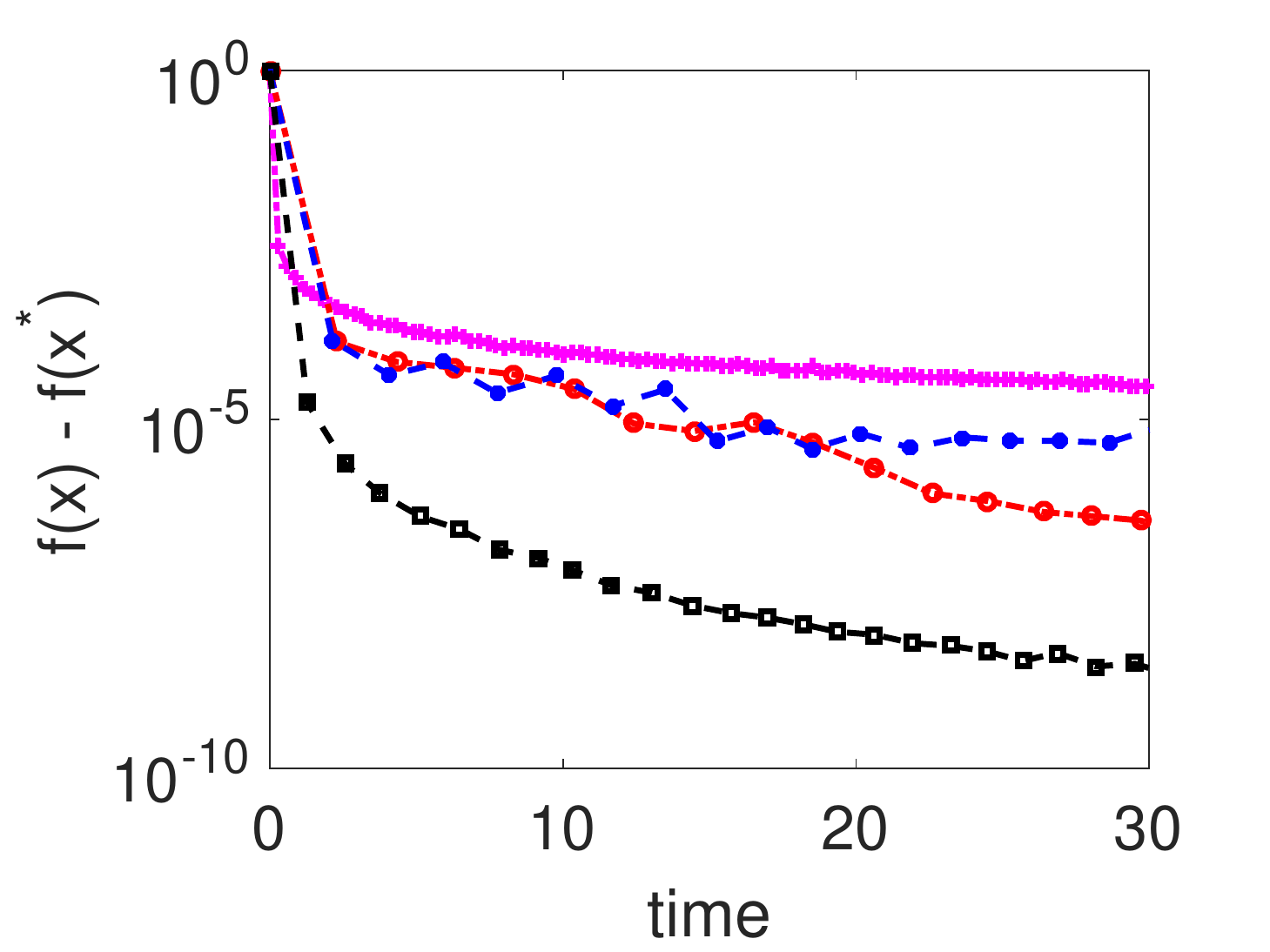}}
	\hspace{-0.15in}
	\subfigure[real-sim]{\includegraphics[width=1.55in]{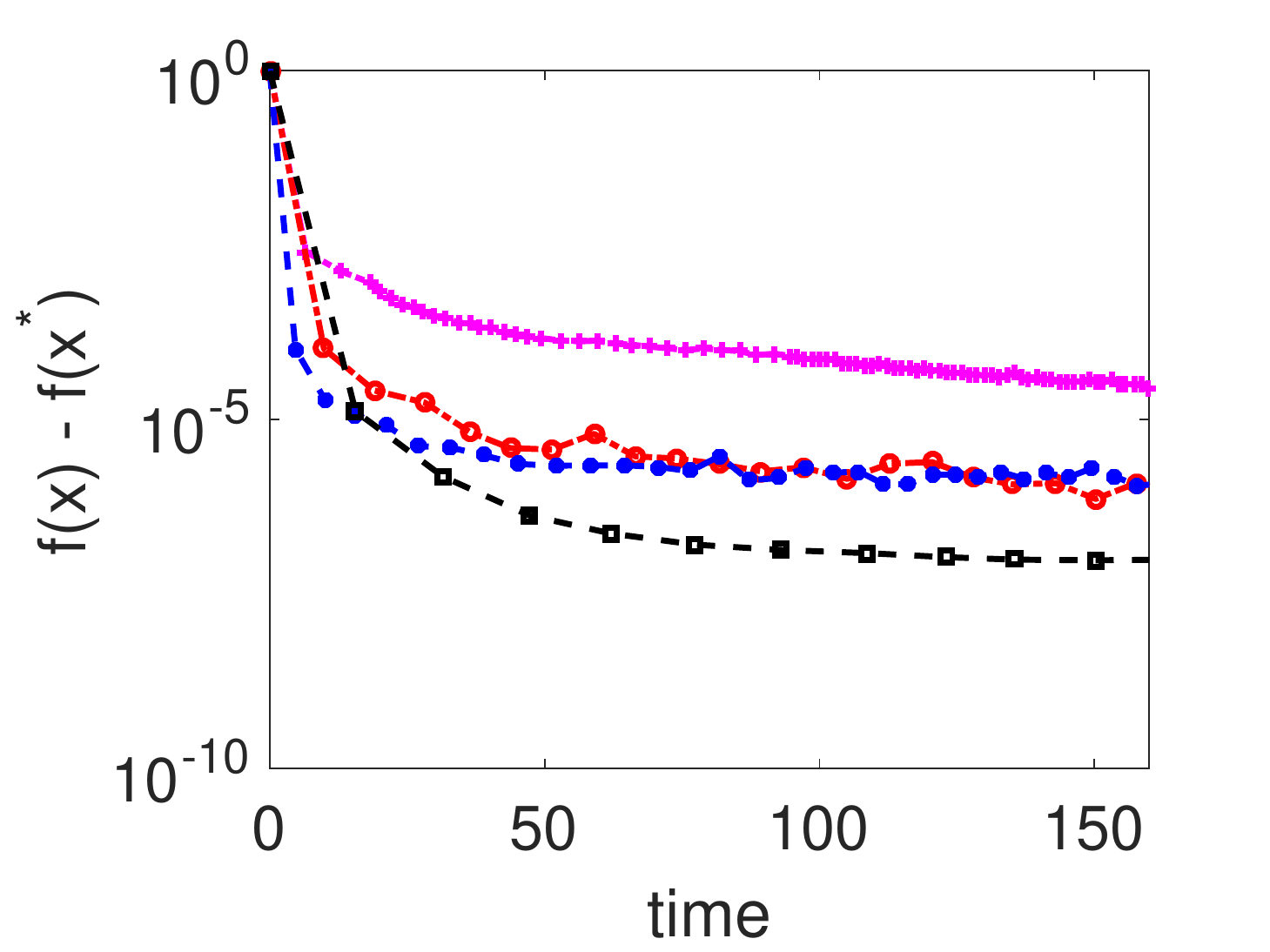}}
	\hspace{-0.15in}
	\subfigure[new20]{\includegraphics[width=1.55in]{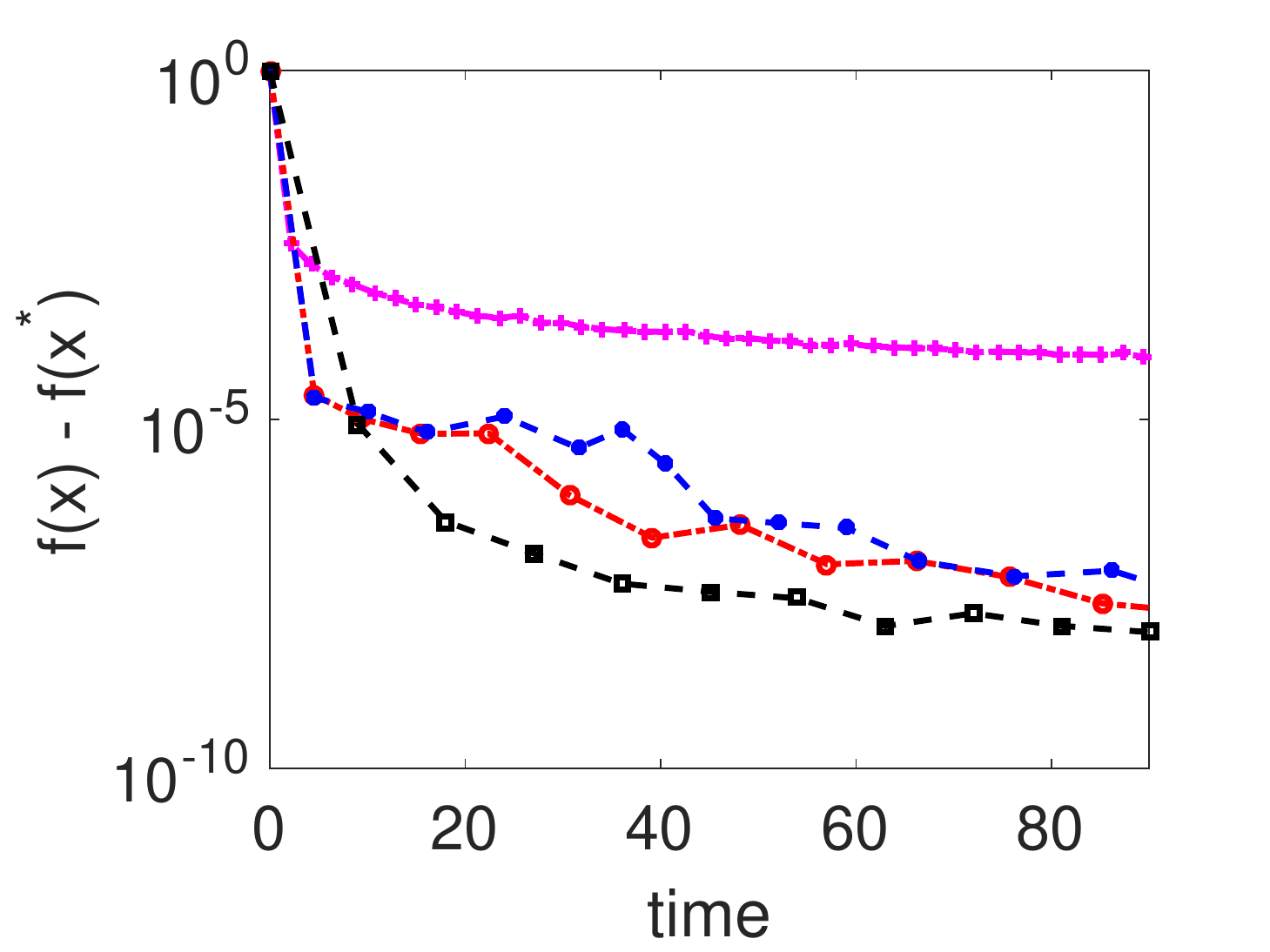}}
	\hspace{-0.15in}
	\subfigure[url]{\includegraphics[width=1.55in]{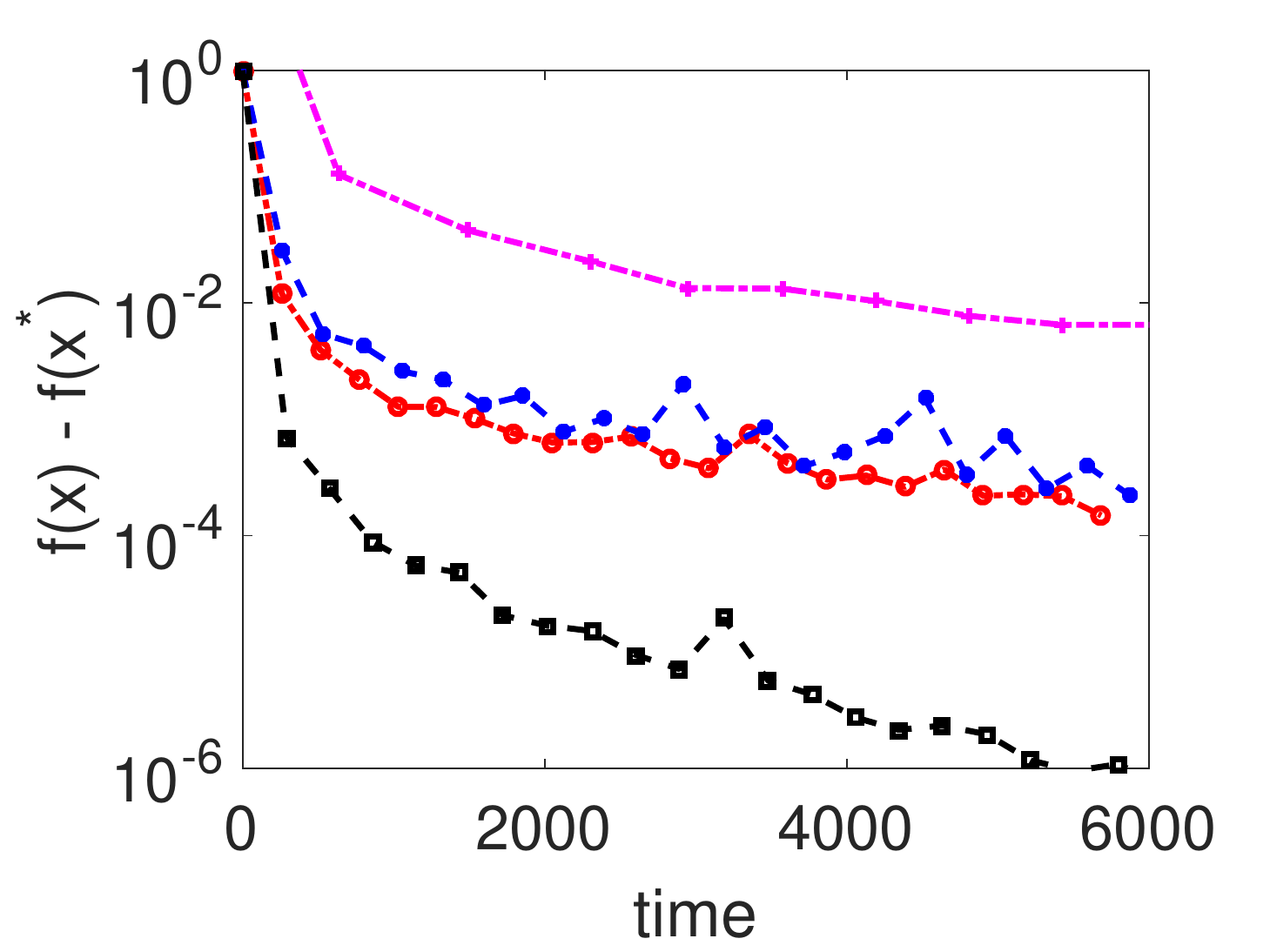}}
	\vspace{-0.1cm}
	
	\includegraphics[height=0.17in]{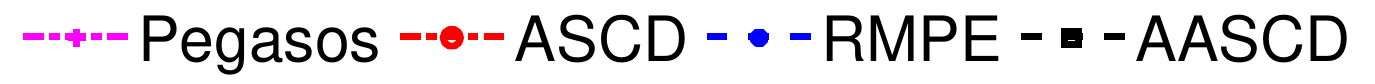}
	\vspace{-0.3cm}
	\caption{Residuals vs CPU training time (s) for solving the Linear SVM problem on four test datasets.}
	\label{fig:dual}
	\vspace{-0.3cm}
\end{figure*}

\subsection{Sparse Dataset}\label{sparse}
One crucial application for asynchronous algorithms in shared memory systems is to solve sparse data. The main reasons are two folded: (1) as the non-zeros coordinates for samples are varying, the computation cost is different for each thread. In this case, asynchronous algorithms are more practical than synchronous ones because threads do not need to wait for synchronization; (2) the data matrix are sparse and ``disjoint'', so the delay effect will be largely reduced. We formulate this fact in the following proposition:
\begin{proposition}
	For a given dataset, if each example is generated  i.i.d and has non-zero component $i$ with probability $\beta_i$, then $\E(m_i(k, j(k)))= \beta_i (k-j(k))$, where $m_i(k_i,k_j)$ is the total number of nonzero updates in component i from iteration $k_j$ to iteration $k_i$.
\end{proposition}
Since $\beta \ll 1$, the delay effect is largely reduced. \citet{asvrg} proposes the $ \triangle$-assumption to judge the sparsity (see Supplementary Material). Under this assumption, our algorithms are able to achieve linear speed up. For example, for AASVRG we have the following property:

\begin{proposition}
	Under the $ \triangle$-assumption ($\triangle \ll 1$) proposed by \citep{asvrg}, for Algorithm~\ref{AASVRG-ori}, the IFO calls is $O\left( n +  (n+(1+\triangle\tau)\sqrt{n})\sqrt{L/\epsilon}\right)$ for NC and $O\left(n+(1+\triangle\tau) \sqrt{nL/\mu}\log(1/\epsilon)\right)$ for SC, respectively. Thus linear speedup is achievable.
\end{proposition}

%
%
%
%

\begin{figure*}[t]
	\center
	\vspace{-0.2cm}
	\hspace{-0.15in}
	\subfigure[rcv1~(iters)]    {\includegraphics[width=1.55in]{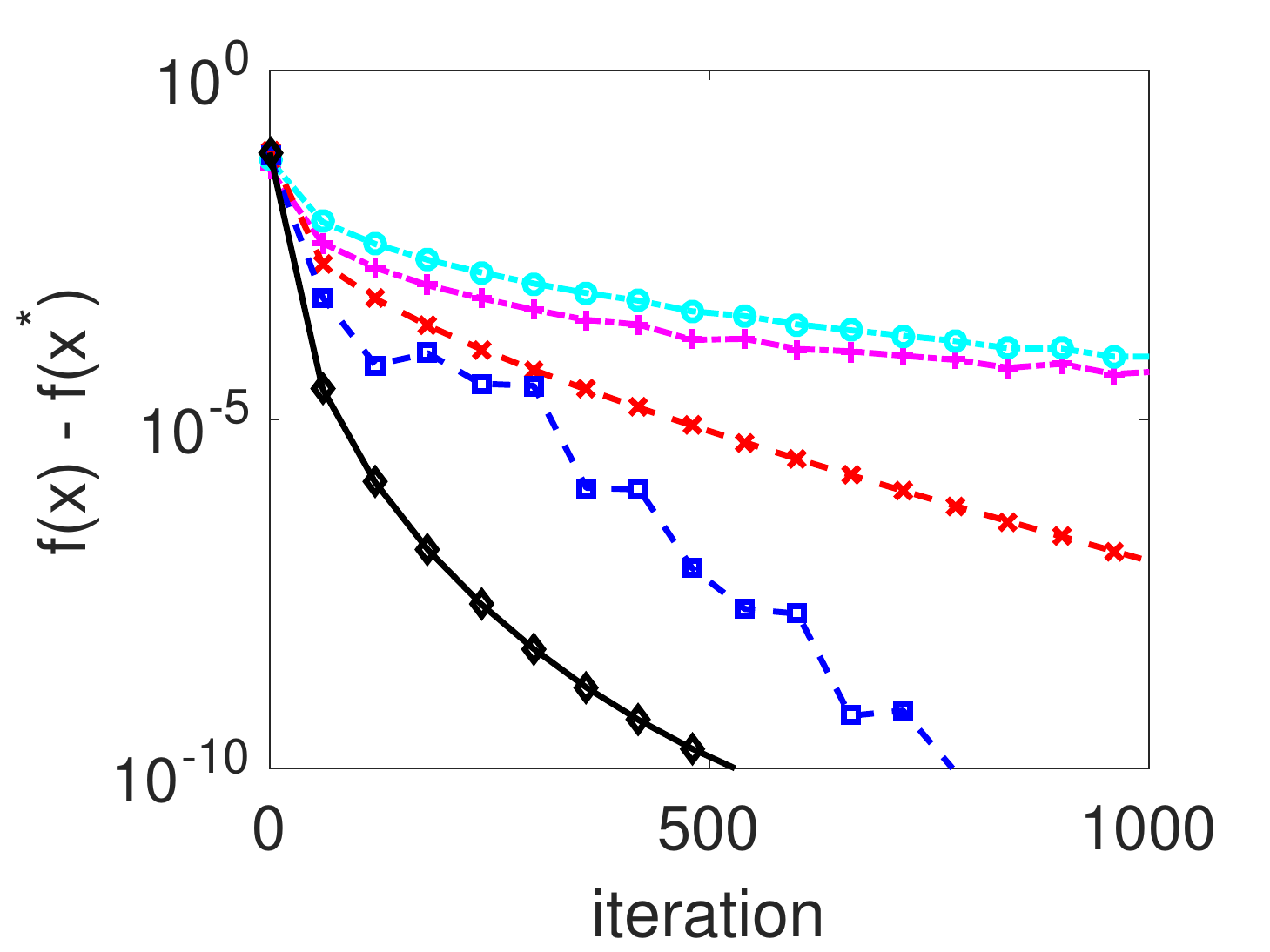}}
	\hspace{-0.15in}
	\subfigure[real-sim~(iters)]{\includegraphics[width=1.55in]{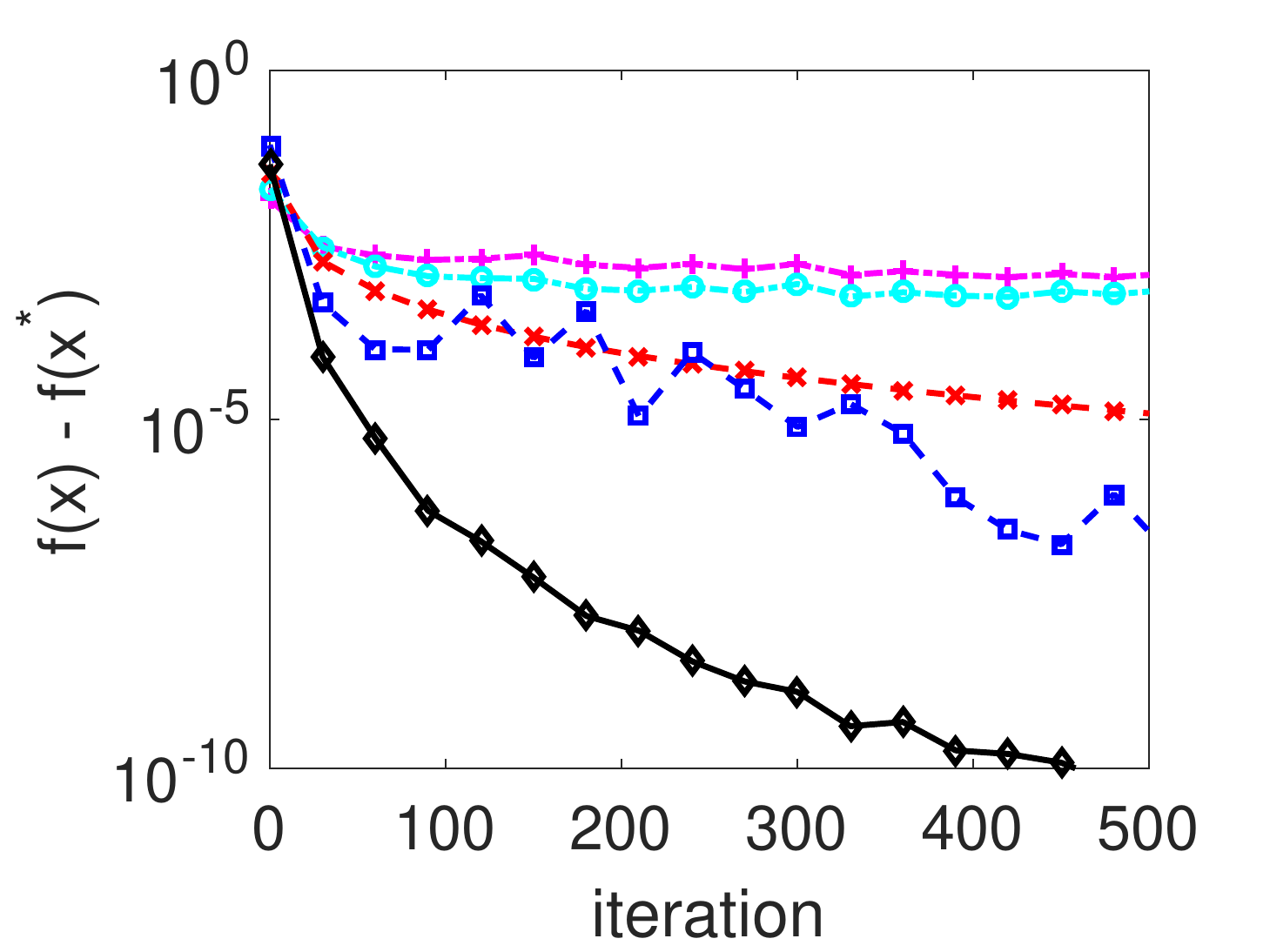}}
	\hspace{-0.15in}
	\subfigure[mnist~(iters)]   {\includegraphics[width=1.55in]{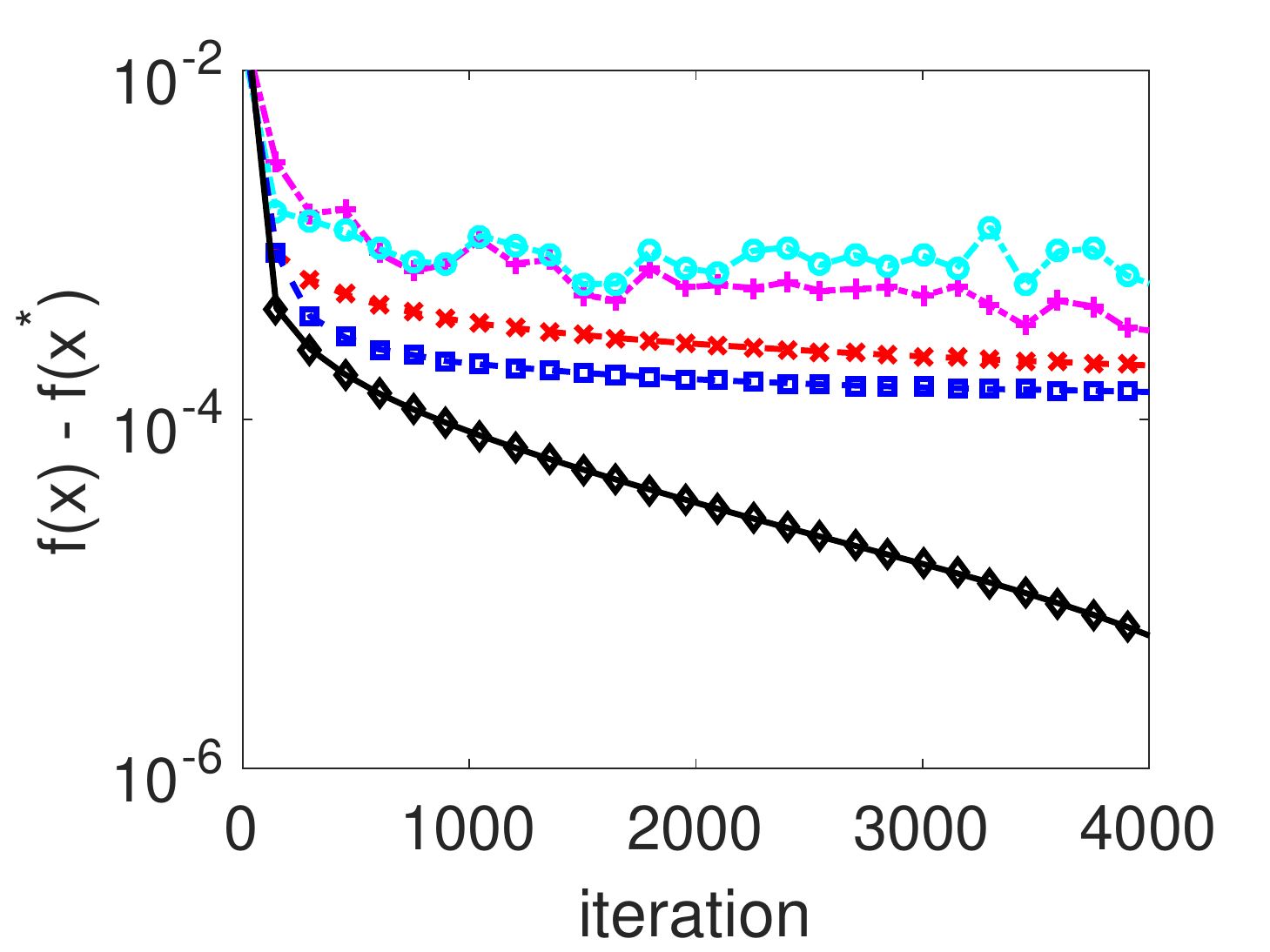}}
	\hspace{-0.15in}
	\subfigure[epsilon~(iters)] {\includegraphics[width=1.55in]{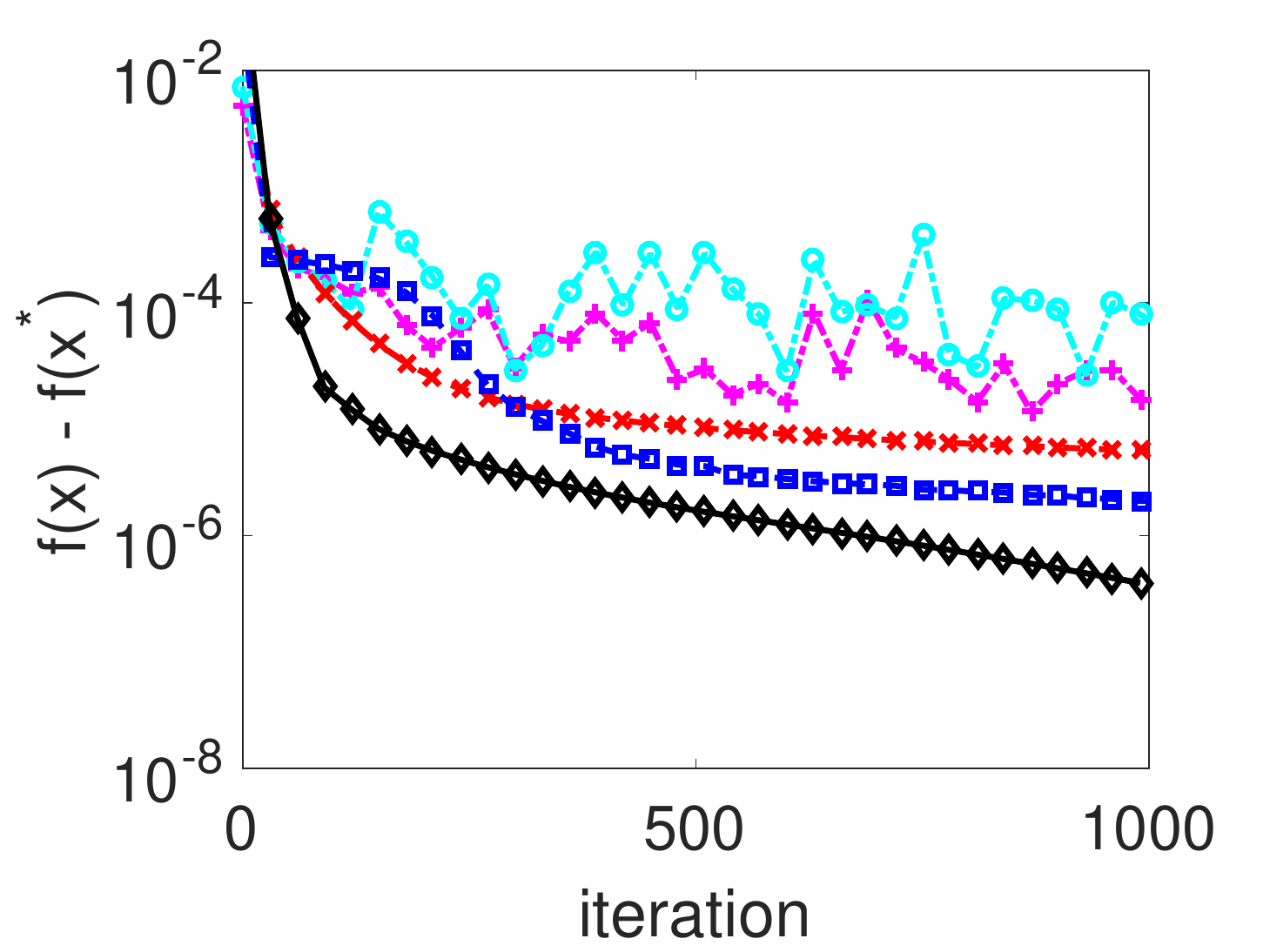}}
	\vspace{-0.3cm}
	
	\hspace{-0.15in}
	\subfigure[rcv1~(time)]    {\includegraphics[width=1.55in]{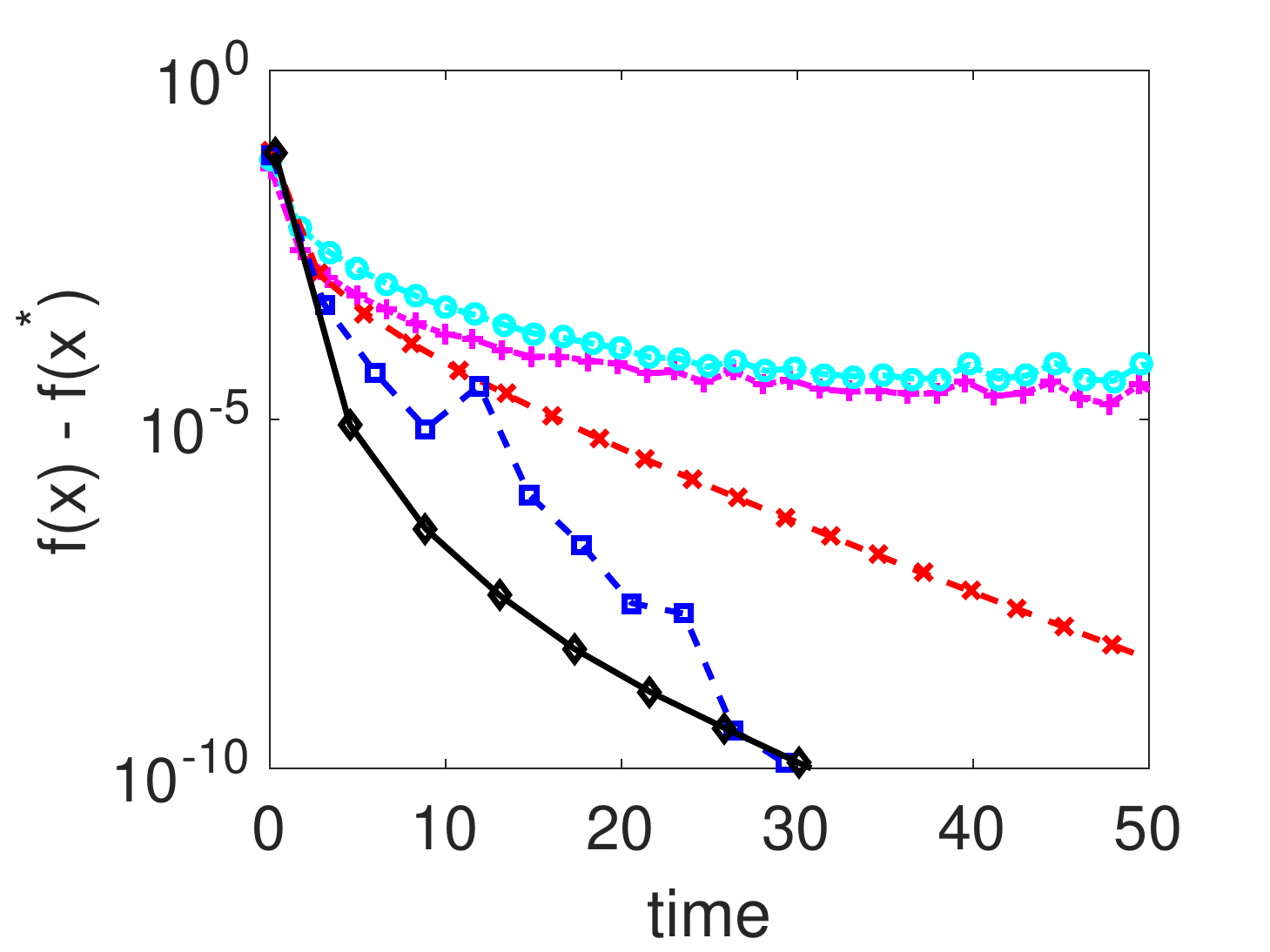}}
	\hspace{-0.15in}
	\subfigure[real-sim~(time)]{\includegraphics[width=1.55in]{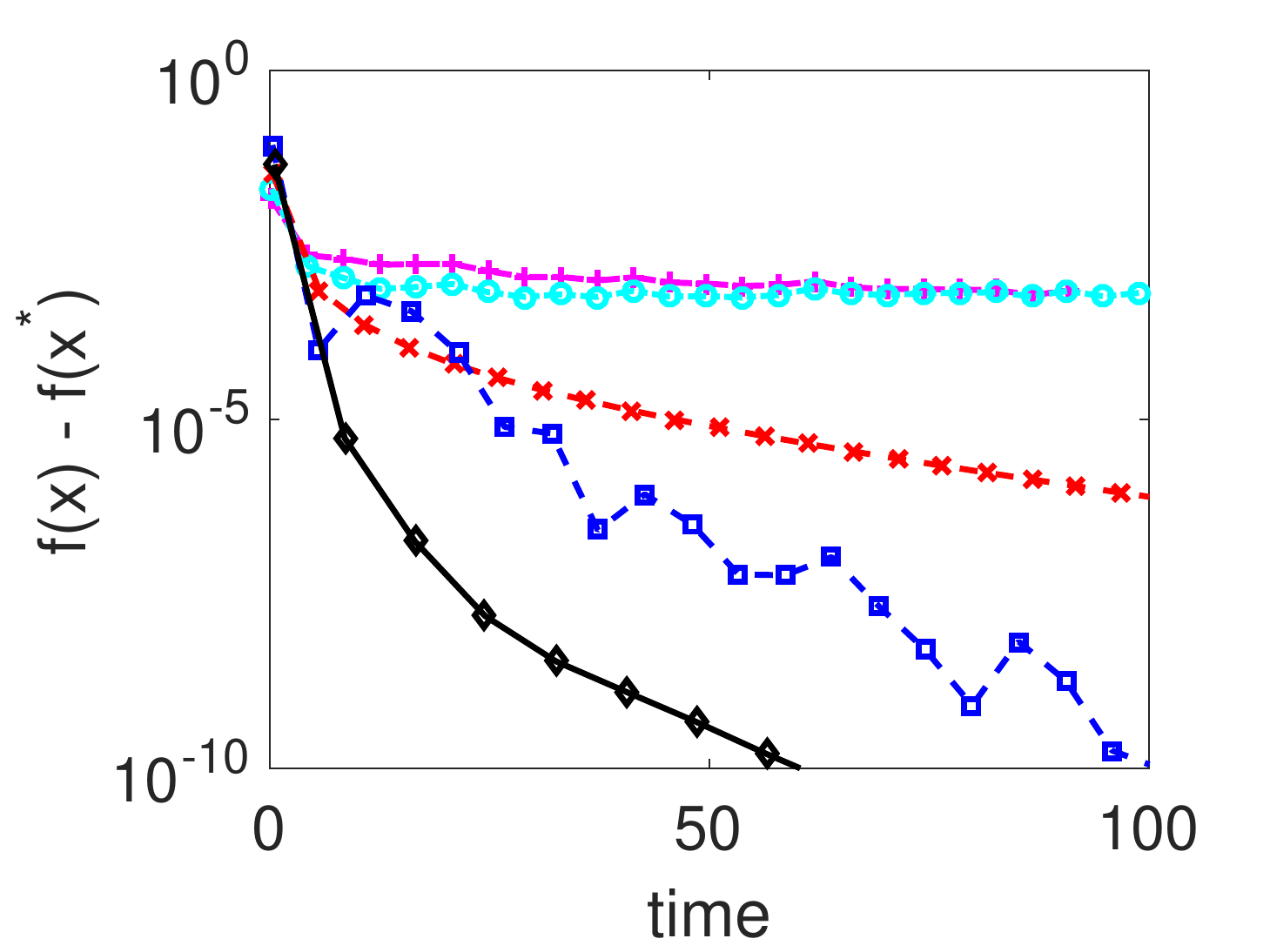}}
	\hspace{-0.15in}
	\subfigure[mnist~(time)]   {\includegraphics[width=1.55in]{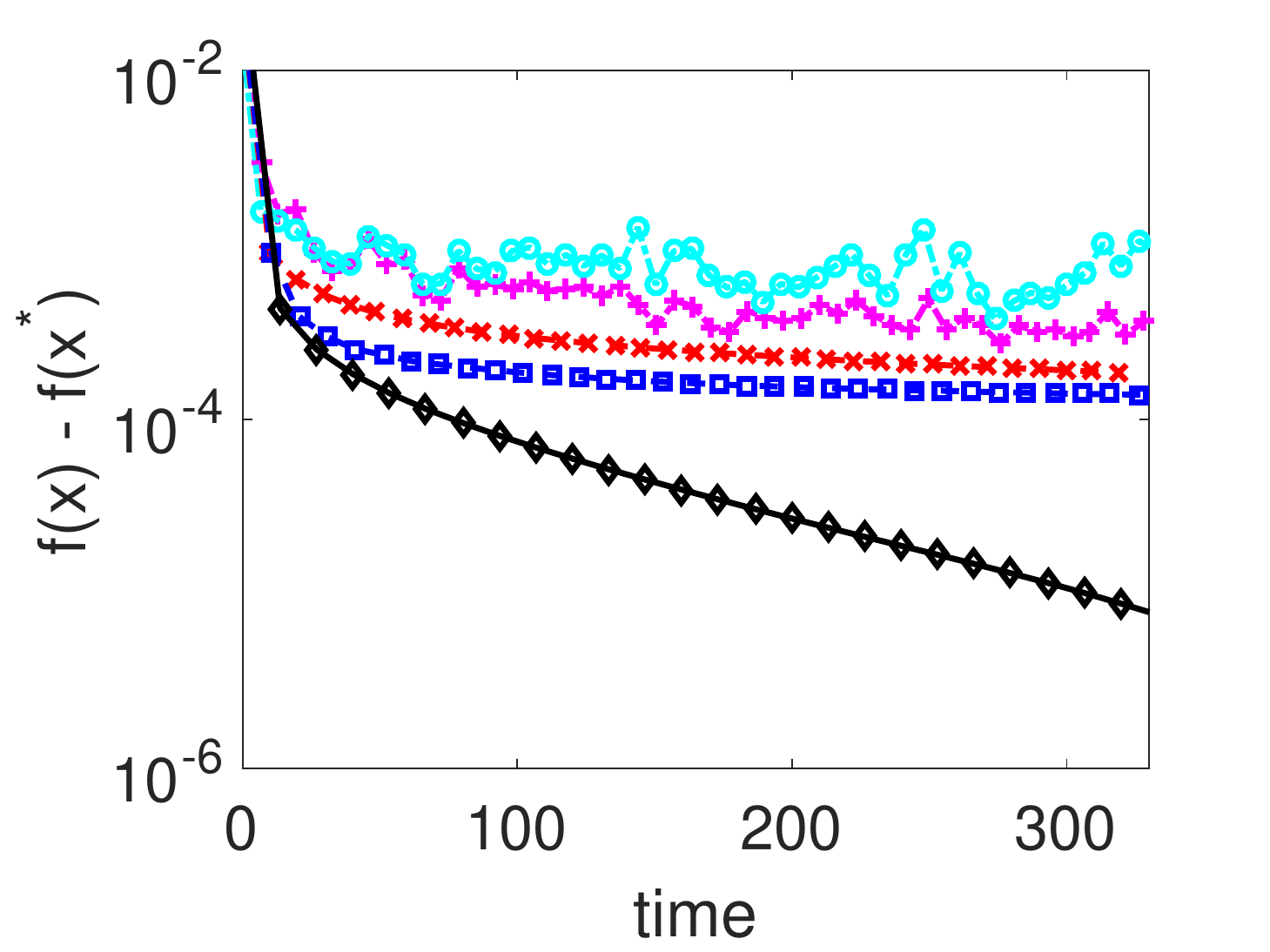}}
	\hspace{-0.15in}
	\subfigure[epsilon~(time)] {\includegraphics[width=1.55in]{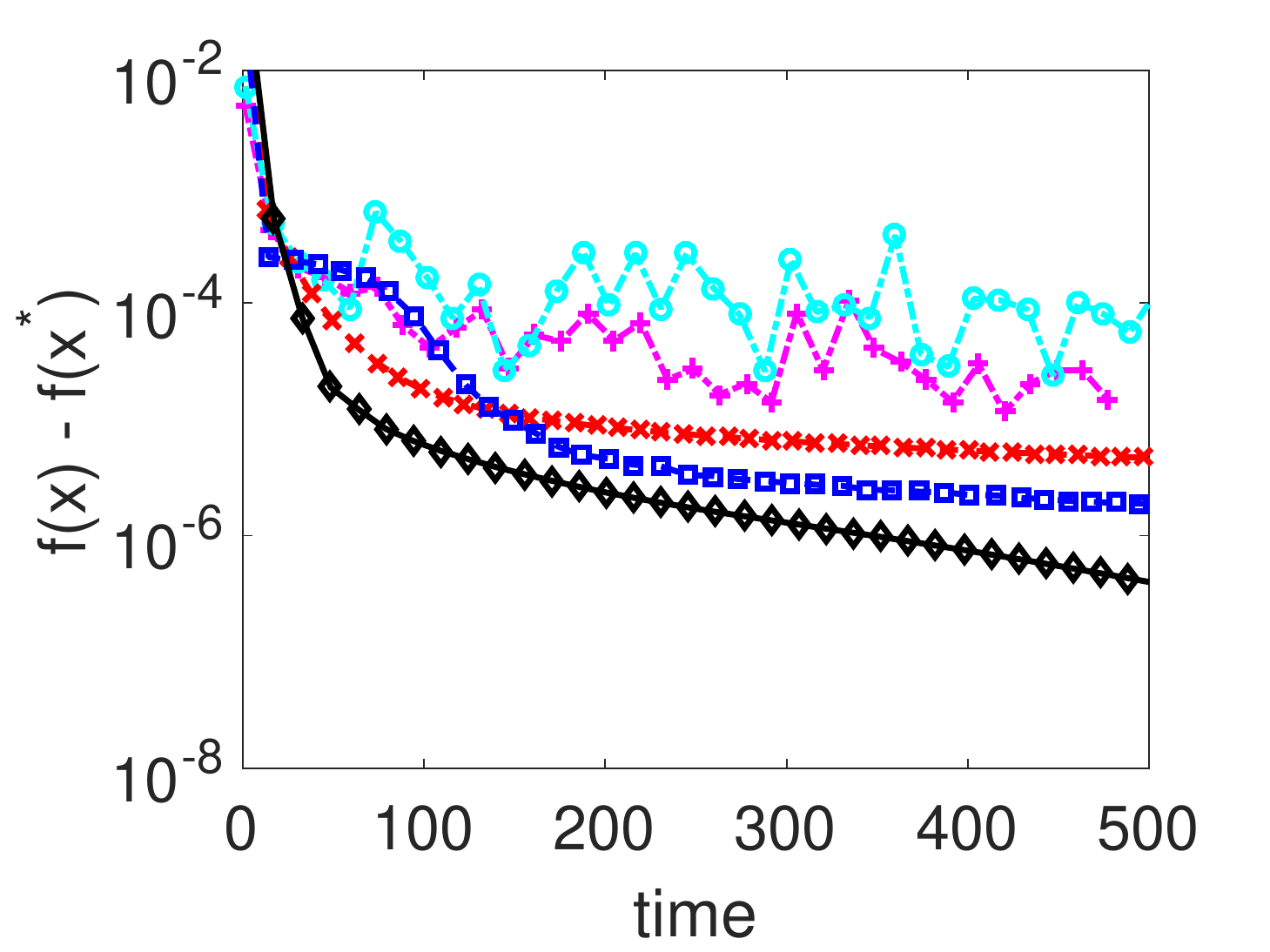}}	
	
	\includegraphics[height=0.17in]{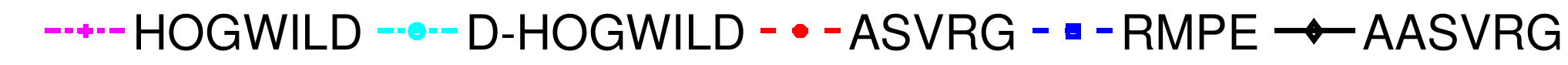}
	\vspace{-0.3cm}
	\caption{Residuals vs Iterations and CPU training time (s) for solving the Ridge Regression problem on four test datasets.}
	\label{fig:primal}
	\vspace{-0.4cm}
\end{figure*}
\section{Experiments}
We have conducted extensive experiments to demonstrate the effectiveness of our method.  We study the problem of Linear SVM for AASCD and Ridge Regression for AASVRG.  We have performed experiments on lots of datasets. For sparse datasets, we choose  four benchmark sparse datasets rcv1, real-sim, news20, and url\footnote{These datasets can be downloaded from \url{https://www.csie.ntu.edu.tw/\~cjlin/libsvmtools/datasets/}.}.  The details of the datasets are shown in  Table~\ref{tab:dataset}. Similar to ~\citep{asvrg}, we have a careful implementation for sparse gradient and computation. We mainly focus on ill-condition problems, so we set the regularizer weight to be $1/(100n)$ in all experiments, and we tune the step size to give the best convergence results. All experiments are done on an Intel multi-core $4$-socket machine with each one contains $8$ cores.

\subsection{Experiments for AASCD}
We compare AASCD with the following methods: 1) Pegasos~\citep{shalev2011pegasos}, which can be considered as one of the best single thread implementation to solve Linear SVM; 2) ASCD~\citep{ASCD}; we also compare some accelerated algorithms, though they have no theoretic guarantees. We compare with RMPE~\citep{scieur2016regularized}, which a regularized nonlinear acceleration algorithm, we also implement   AASGD~\citep{meng2016asynchronous}. However, we find it hard to converge on the sparse data.

\subsection{Experiments for AASVRG}
We compare AASVRG with the following methods: 1) HOGWILD, a lock-free asynchronous variant of  SGD and D-HOGWILD, which  chooses decaying step size as $\eta_0 \sqrt{\sigma_0 / (t + \sigma_0)}$); 2) ASVRG, the lock-free asynchronous variant of SVRG; and 3) RMPE~\citep{scieur2016regularized}, a regularized nonlinear acceleration algorithm.

\subsection{Experiments for Dense Dataset}
We also do experiments on Dense Dataset for AASVRG which  performs on USPS, SENSIT, MNIST, and EPSILON with the details shown in Supplementary Material.

\subsection{Results}
We first  measure the speedup achieved by our algorithms on the sparse dataset. The time speed up is defined as the ratio of the runtime to achieve a given precision with the serial implementation to the runtime with the asynchronous implementation with $P$ threads, and the iteration speedup is defined as
\begin{equation}
\text{iteration speedup} = \frac{\# \text{of iters for seri. algorithms}}{\# \text{of iters for asyn. algorithms}} \times P. \notag
\end{equation}
The results are shown in Fig~\ref{fig:speedup}. There is a linear speedup for iteration, and a nearly linear speedup for time, which verifies our theorem on the sparse dataset. Asynchronous algorithm achieves higher speedup than synchronous one.

To compare these algorithms, we consider the training loss residual versus CPU time. For AASCD, the results is shown in Fig.~\ref{fig:dual}. The experiments are conducted on $10$ cores. It is clear that our algorithm  converges fastest in all four datasets among other the algorithms.

The experiment results for AASVRG is shown in Fig.~\ref{fig:primal}. To demonstrate that our algorithm has a faster speed, we also reports the training loss residual versus iteration. It is also  clear that our algorithm are much faster.

Due to space limit,  implementation details and more experimental results, e.g. variant regularizer weight terms are shown in Supplementary Material.

\input{appendix}

\bibliographystyle{custom}
\setlength{\bibsep}{3pt}
\bibliography{AAGD_MC}

\end{document}

%% file: appendix.tex
\section{Supplementary Materials}

The Supplementary Material is structured as follows: in Section~\ref{AAGD-sec}, we give the proof for AAGD;  in Section~\ref{AASCD-sec}, we give the proof for AASCD; in Section~\ref{AASVRG-sec}, we give the proof for AASVRG; Also an outline of the proof is at the beginning of each Section.   In Section~\ref{ASVRG-sec}, we give the proof for ASVRG.  In Section~\ref{Other Material-sec},  we show  some implementation details and more experimental results.

\subsection{AAGD}\label{AAGD-sec}
We set 
\begin{eqnarray}\label{b1}
\y^k = (1-\theta^k)\x^k +\theta^k \z^k.
\end{eqnarray}
The through the step 4 in Algorithm 2 in the paper, we have
\begin{eqnarray}\label{b2}
\x^{k+1} = \y^k +\theta^k\bdelta^k.
\end{eqnarray}
\bfseries Outline of the Proof: \mdseries \\
Step 1: Through the update rule, we have that
\begin{eqnarray}\label{step1}
\y^{k}- \w^{j(k)} = \sum_{i = j(k)+1}^k \left( 1+ \sum_{l = i}^k b(i,l) \right) (\x^i -\y^{i-1}). 
\end{eqnarray}
Step 2: By analyzing the function value,  we have
\begin{eqnarray}\label{im11}
f(\x^{k+1}) &\leq& f(\y^k)-\gamma (1-\frac{\gamma L}{2})\left\| \frac{\x^{k+1}-\y^k}{\gamma} \right\|^2-\langle \bxi^k, \x^{k+1}-\y^k\rangle\notag\\
&& + \frac{\gamma L^2}{2C_1}\left\|  \w^{j(k)}-\y^k \right\|^2 + \frac{\gamma C_1}{2}\left\| \frac{\x^{k+1}-\y^k}{\gamma} \right\|^2. 
\end{eqnarray}
Step 3: By analyzing the $\left\|\z^{k+1}-\x^*\right\|^2$,  we have
\begin{eqnarray}\label{im2}
&&\frac{1}{2\gamma}\left\|\theta^k \z^{k+1} -\theta^k \x^*  \right\|^2\\
&=&\frac{1}{2\gamma}\left\| \theta^k \z^{k} -\theta^k \x^*\right\|^2 +\frac{1}{2\gamma}\left\|\x^{k+1}-\y^k \right\|^2-\langle  \bxi^k,  \theta^k\z^{k}-\theta^k\x^*  \rangle\notag\\
&&+(1-\theta^k)f(\x^k)+\theta^k f(\x^*) -f(\y^k)+\langle \nabla f(\y^k)-\nabla f(\w^{j(k)}), \y^k- \w^{j(k)}\rangle\notag.
\end{eqnarray}
Step 4: By adding Eq.~\eqref{im11} and Eq.~\eqref{im2},  we have
\begin{eqnarray}\label{im3}
F(\x^{k+1})&\leq&(1-\theta^k)F(\x^k)+\theta^k F(\x^*)-\gamma (\frac{1}{2}-\frac{\gamma L}{2})\left\| \frac{\x^{k+1}-\y^k}{\gamma} \right\|^2\\
&&+\gamma\left(\frac{\gamma^2 L^2}{2C_1}+\gamma L\right)\left\| \frac{\w^{j(k)} - \y^k}{\gamma} \right\|^2 + \frac{\gamma C_1}{2}\left\|\frac{\x^{k+1}-\y^k }{\gamma}\right\|^2 \notag\\
&&+\frac{1}{2\gamma} \left\| \theta^k \z^{k} -\theta^k \x^*\right\|^2-\left(\frac{1}{2\gamma}+\frac{\mu}{2\theta^k} \right) \left \| \theta^k \z^{k+1} -\theta^k \x^*\right\|^2\notag,
\end{eqnarray}
Step 5: we choose proper step size and obtain Theorem 1 in the paper.\\
\bfseries  Proof of step 1: \mdseries \\
Through Eq.~\eqref{b1}, we have
\begin{eqnarray}
\theta^k \z^k = \y^k-(1-\theta^k)\x^k, \quad k\geq0.
\end{eqnarray}
and through the Step 4 in Algorithm 2 in the paper,
\begin{eqnarray}
\theta^k \z^{k+1} = \x^{k+1}-(1-\theta^k)\x^k, \quad k\geq 0.
\end{eqnarray}
Eliminating $\z^k$, we have
\begin{eqnarray}
\frac{\y^k-(1-\theta^k)\x^k}{\theta^k} = \frac{ \x^{k}-(1-\theta^k)\x^{k-1}}{\theta^{k-1}}, \quad k\geq 1.
\end{eqnarray}
Thus
\begin{eqnarray}\label{extra}
\y^k = \x^k +\frac{\theta^k(1-\theta^k)}{\theta^{k-1}}\left( \x^k-\x^{k-1}\right), \quad k\geq 1.
\end{eqnarray}
Set $a^k =\frac{\theta^k(1-\theta^k)}{\theta^{k-1}}$,  we have $a^k \leq 1$. We have
\begin{eqnarray}\label{22}
\y^k &=& \x^k + a^k (\x^k - \y^{k-1}) +a^k(\y^{k-1}-\x^{k-1})\\
&=& \y^{k-1} +(a^k+1)(\x^k-\y^{k-1}) +a^ka^{k-1}(\x^{k-1}-\x^{k-2}), \quad k\geq2.\notag
\end{eqnarray} 
For $\x^{k-1}-\x^{k-2}$, and  $k\geq j(k)+2\geq 2$, we have
\begin{eqnarray}\label{one11}
&&\x^{k-1}-\x^{k-2}\\
&=&\x^{k-1} -\y^{k-2} +\y^{k-2}-\x^{k-2}\notag\\
&=&\x^{k-1} -\y^{k-2}+a^{k-2}(\x^{k-2}-\x^{k-3})\notag\\
&=&\x^{k-1} -\y^{k-2}+a^{k-2}(\x^{k-2} -\y^{k-3})+a^{k-2}a^{k-3}(\x^{k-3}-\x^{k-4})\notag\\
&=&\x^{k-1} -\y^{k-2}+ \sum_{i=j(k)+1}^{k-2}   \left(\left(\prod_{l=i}^{k-2} a^{l}\right) (\x^{i} -\y^{i-1})\right)+\left(\prod_{l=j(k)}^{k-2} a^{l}\right)(\x^{j(k)}-\x^{j(k)-1}).\notag
\end{eqnarray}
Set $b(l,k) =\prod_{i=l}^{k} a^{i}$, where $l\leq k $. Substituting Eq.~\eqref{one11} into Eq.~\eqref{22}, we have
\begin{eqnarray}\label{222}
&&\y^k \\ 
&=&  \y^{k-1}+(b(k,k)+1)(\x^k-\y^{k-1})+b(k-1,k)(\x^{k-1} -\y^{k-2})\notag\\
&&+\sum_{i=j(k)+1}^{k-2}\left(b(i,k)(\x^i-\y^{i-1})\right) +b(j(k),k)(\x^{j(k)}-\x^{j(k)-1})\notag\\
&=&\y^{k-1}+(\x^k-\y^{k-1})+\sum_{i=j(k)+1}^{k}\left(b(i,k)(\x^i-\y^{i-1})\right) +b(j(k),k)(\x^{j(k)}-\x^{j(k)-1})\notag.
\end{eqnarray}
By checking, when $k=j(k)$ and $k = j(k)+1$, Eq.~\eqref{222} is right. So Eq.~\eqref{222} holds for any $k\geq j(k)\geq0$.
Summing Eq.~\eqref{222} with $k = j(k)+1$ to $k$, we have
\begin{eqnarray}\label{133}
&&\y^k\\
&=& \y^{j(k)}+\sum_{i=j(k)+1}^{k}(\x^{i}-\y^{i-1})+\sum_{l=j(k)+1}^{k}\sum_{i=j(k)+1}^{l} b(i,l)(\x^i-\y^{i-1})\notag\\
&&+\left(\sum_{i=j(k)+1}^{k}b(j(k),i)\right)(\x^{j(k)}-\x^{j(k)-1})\notag\\
&\overset{a}=& \y^{j(k)}+\sum_{i=j(k)+1}^{k}(\x^{i}-\y^{i-1})+\sum_{i=j(k)+1}^{k}\left(\sum_{l=i}^{k} b(i,l)\right)(\x^i-\y^{i-1})\notag\\
&&+\left(\sum_{i=j(k)+1}^{k}b(j(k),i)\right)(\x^{j(k)}-\x^{j(k)-1})\notag\\
&\overset{Eq.~\eqref{extra}}=& \x^{j(k)}+\sum_{i=j(k)+1}^{k}(\x^{i}-\y^{i-1})+\sum_{i=j(k)+1}^{k}\left(\sum_{l=i}^{k} b(i,l)\right)(\x^i-\y^{i-1})\notag\\
&&+\left(\sum_{i=j(k)}^{k}b(j(k),i)\right)(\x^{j(k)}-\x^{j(k)-1})\notag,
\end{eqnarray}
where  $\overset{a}=$ is obtained by rearrange terms.
Then by comparing the results, we obtain Step 1.\\
\bfseries  Proof of step 2: \mdseries \\
Through the optimal solution of  $\z^{k+1}$ in Step 2 of Algorithm 2 in the paper, we have that
\begin{eqnarray}\label{a155}
\theta^k(\z^{k+1}-\z^k)+\gamma\nabla f(\w^{j(k)}) +\gamma\bxi^k =\mathbf{0},
\end{eqnarray}
where $\bxi^k\in \partial h(\z^{k+1})$. And through Eq.~\eqref{b2}, we have
\begin{eqnarray}\label{a166}
(\x^{k+1}-\y^k) +\gamma\nabla f(\w^{j(k)}) +\gamma\bxi^k =\mathbf{0}.
\end{eqnarray}
For $f$ has  Lipschitz continues gradient, we obtain
\begin{eqnarray}\label{44}
f(\x^{k+1}) &\leq& f(\y^k) +\langle \nabla f(\y^k), \x^{k+1}- \y^k\rangle +\frac{L}{2}\left\|\x^{k+1}- \y^k \right\|^2\notag\\
&\overset{a}=& f(\y^k)-\gamma\langle \nabla f(\y^k), \nabla f(\w^{j(k)} )+\bxi^k \rangle +\frac{L}{2}\left\|\x^{k+1}- \y^k \right\|^2\notag\\
&\overset{b}=&f(\y^k)-\gamma\langle \nabla f(\w^{j(k)} )+\bxi^k, \nabla f(\w^{j(k)} )+\bxi^k \rangle +\frac{L}{2}\left\|\x^{k+1}- \y^k\right \|^2\notag\\
&&+\gamma \langle \bxi,\nabla f(\w^{j(k)} )+\bxi^k  \rangle+\gamma\langle    \nabla f(\w^{j(k)} ) - \nabla f(\y^k), \nabla f(\w^{j(k)} )+\bxi^k  \rangle\notag\\
&\overset{c}=&f(\y^k)-\gamma (1-\frac{\gamma L}{2})\left\| \frac{\x^{k+1}-\y^k}{\gamma}\right \|^2-\langle \bxi^k, \x^{k+1}-\y^k\rangle\notag\\
&& -\langle    \nabla f(\w^{j(k)} ) - \nabla f(\y^k),  \x^{k+1}-\y^k   \rangle,
\end{eqnarray}
where $\overset{a}=$, we use Eq.~\eqref{a166}; in $\overset{b}=$, we use $-\nabla f(\y^k) =  -\nabla f(\w^{j(k)} ) - \bxi^k+ \bxi^k-\nabla f(\y^k) +\nabla f(\w^{j(k)})$; in $\overset{c}=$, we reuse Eq.~\eqref{a166}.

For the last term of Eq.~\eqref{44}, applying Cauchy-Schwarzwe inequality, we have
\begin{eqnarray}\label{111}
&&\langle   \nabla f(\w^{j(k)} ) - \nabla f(\y^k),  \x^{k+1}-\y^k     \rangle \notag\\
&\leq& \frac{\gamma}{2C_1}\left\| \nabla f(\w^{j(k)} ) -\nabla f(\y^k) \right\|^2 + \frac{\gamma C_1}{2}\left\|\frac{\x^{k+1}-\y^k }{\gamma}\right\|^2\notag\\
&\leq&\frac{\gamma L^2}{2C_1}\left\| \w^{j(k)} - \y^k \right\|^2 + \frac{\gamma C_1}{2}\left\|\frac{\x^{k+1}-\y^k }{\gamma}\right\|^2.
\end{eqnarray}
Substituting Eq.~\eqref{111} into Eq.~\eqref{44}, we obtain the results of Step 2.\\
\bfseries  Proof of step 3: \mdseries \\
\begin{eqnarray}\label{66}
&&\frac{1}{2\gamma}\left\|\theta^k \z^{k+1} -\theta^k \x^*  \right\|^2\\
&=& \frac{1}{2\gamma}\left\| \theta^k \z^{k} -\theta^k \x^* + \theta^k \z^{k+1}-\theta^k\z^k  \right\|^2\notag\\
&=&\frac{1}{2\gamma} \left\| \theta^k \z^{k} -\theta^k \x^*\right\|^2 +\frac{1}{2\gamma} \left\| \theta^k \z^{k+1}-\theta^k\z^k  \right\|^2 +\frac{1}{\gamma}\langle \theta^k \left(\z^{k+1}-\z^k\right) ,  \theta^k\z^{k}-\theta^k\x^* \rangle\notag\\
&\overset{a}=&\frac{1}{2\gamma}\left\| \theta^k \z^{k} -\theta^k \x^*\right\|^2 +\frac{1}{2\gamma}\left\|\x^{k+1}-\y^k \right\|^2-\langle \nabla f(\w^{j(k)})+\bxi^k,  \theta^k\z^{k}-\theta^k\x^*  \rangle,\notag
\end{eqnarray}
where in $\overset{a}=$, we use Eq.~\eqref{a155}. Then for the last term, we have that
\begin{eqnarray}\label{77}
&&-\langle   \nabla f(\w^{j(k)}),  \theta^k\z^{k}-\theta^k\x^* \rangle\\
&\overset{\eqref{b1}}=&-\langle  \nabla   f(\w^{j(k)}),  \y^k-(1-\theta^k)\x^k-\theta^k\x^* \rangle\notag\\
&\overset{a}=&-\langle  \nabla  f(\w^{j(k)}),  \w^{j(k)}-(1-\theta^k)\x^k-\theta^k\x^* \rangle -\langle\nabla f(\w^{j(k)}), \y^k- \w^{j(k)}\rangle\notag\\
&\overset{b}\leq& (1-\theta^k)f(\x^k)+\theta^k f(\x^*) -f(\w^{j(k)})-\langle \nabla f(\w^{j(k)}), \y^k- \w^{j(k)}\rangle\notag\\
&\overset{c}\leq& (1-\theta^k)f(\x^k)+\theta^k f(\x^*) -f(\y^k)+\langle \nabla f(\y^k)-\nabla f(\w^{j(k)}), \y^k- \w^{j(k)}\rangle\notag,
\end{eqnarray}
where $\overset{a}=$, we insert $\w^{j(k)}$; in $\overset{b}\leq$, we use the convexity of $f$, namely applying
$$  f(\w^{j(k)}) +\langle \nabla f(\w^{j(k)}),   \a -\w^{j(k)}  \rangle \leq f(\a), $$
on $\a = \x^*$, and $\a = \x^k$, respectively; in $\overset{c}\leq$, we use that 
\begin{eqnarray}
-f(\w^{j(k)}) \leq -f(\y^k) +\langle \nabla f(\y^k), \y^k- \w^{j(k)}\rangle.
\end{eqnarray}
Substituting Eq.~\eqref{77} into Eq.~\eqref{66}, after simplifying, we obtain the result of Step 3.\\
\bfseries  Proof of step 4: \mdseries 
Adding Eq.~\eqref{im11}and Eq.~\eqref{im2},  we have that
\begin{eqnarray}\label{888}
f(\x^{k+1}) &\leq&  (1-\theta^k)f(\x^k)+\theta^k f(\x^*)-\gamma (\frac{1}{2}-\frac{\gamma L}{2})\left\| \frac{\x^{k+1}-\y^k}{\gamma} \right\|^2-\langle \bxi^k, \x^{k+1}-\y^k\rangle\notag\\
&&+\frac{\gamma L^2}{2C_1}\left\| \w^{j(k)} - \y^k \right\|^2 + \frac{\gamma C_1}{2}\left\|\frac{\x^{k+1}-\y^k }{\gamma}\right\|^2-\langle  \bxi^k,  \theta^k\z^{k}-\theta^k\x^*  \rangle \notag\\
&&+\langle \nabla f(\y^k)-\nabla f(\w^{j(k)}), \y^k- \w^{j(k)}\rangle+\frac{1}{2\gamma} \left\| \theta^k \z^{k} -\theta^k \x^*\right\|^2-\frac{1}{2\gamma}\left \| \theta^k \z^{k+1} -\theta^k \x^*\right\|^2\notag\\
&\overset{a}\leq& (1-\theta^k)f(\x^k)+\theta^k f(\x^*)-\gamma (\frac{1}{2}-\frac{\gamma L}{2})\left\| \frac{\x^{k+1}-\y^k}{\gamma} \right\|^2-\langle \bxi^k, \x^{k+1}-\y^k\rangle\notag\\
&&+\left(\frac{\gamma L^2}{2C_1}+L\right)\left\| \w^{j(k)} - \y^k\right \|^2 + \frac{\gamma C_1}{2}\left\|\frac{\x^{k+1}-\y^k }{\gamma}\right\|^2-\langle  \bxi^k,  \theta^k\z^{k}-\theta^k\x^*  \rangle\notag\\
&&+\frac{1}{2\gamma}\left \| \theta^k \z^{k} -\theta^k \x^*\right\|^2-\frac{1}{2\gamma} \left\| \theta^k \z^{k+1} -\theta^k \x^*\right\|^2,
\end{eqnarray}
where in $\overset{a}\leq$, we use $\langle \nabla f(\y^k)-\nabla f(\w^{j(k)}), \y^k- \w^{j(k)}\rangle \leq \| \y^k- \w^{j(k)}\|^2 $.
Since $\bxi\in\partial h(\z^{k+1})$, we have that
\begin{eqnarray}\label{88}
&&-\langle \bxi^k, \x^{k+1}-\y^k\rangle  -\langle  \bxi^k,  \theta^k\z^{k}-\theta^k\x^*  \rangle \overset{\eqref{b2}}{=}\theta^k\langle \bxi^k ,\x^*-\z^{k+1}\rangle\\
&\leq& \theta^k h(\x^*) -\theta^k h(\z^{k+1})-\frac{\mu\theta^k}{2}\left\|\z^{k+1}-\x^* \right\|^2.\notag
\end{eqnarray}
For the convexity of $h(\z^{k+1})$, and the step 4 in Algorithm 2 in the paper, we have
\begin{eqnarray}\label{11111}
\theta^k h(\z^{k+1}) +(1-\theta^k)h(\x^k) \geq h(\x^{k+1}). 
\end{eqnarray}
Substituting Eq.~\eqref{88}  into  Eq.~\eqref{888}, and using Eq.~\eqref{11111}, we have
\begin{eqnarray}
F(\x^{k+1})&\leq&(1-\theta^k)F(\x^k)+\theta^k F(\x^*)-\gamma (\frac{1}{2}-\frac{\gamma L}{2})\left\| \frac{\x^{k+1}-\y^k}{\gamma}\right \|^2\notag\\
&&+\left(\frac{\gamma L^2}{2C_1}+L\right)\left\| \w^{j(k)} - \y^k \right\|^2 + \frac{\gamma C_1}{2}\left\|\frac{\x^{k+1}-\y^k }{\gamma}\right\|^2 \notag\\
&&+\frac{1}{2\gamma}\left \| \theta^k \z^{k} -\theta^k \x^*\right\|^2-\left(\frac{1}{2\gamma}+\frac{\mu}{2\theta^k} \right) \left\| \theta^k \z^{k+1} -\theta^k \x^*\right\|^2\notag.
\end{eqnarray}
\bfseries  Proof of step 5: \mdseries \\
We first consider the not-strongly convex case. Through Eq.~\eqref{step1}, we have
\begin{eqnarray}\label{onee}
&&\left\|\w^{j(k)} -\y^k\right\|^2\\
&=&\left\|\sum_{i=j(k)+1}^{k}\left(1+\sum_{l=i}^{k} b(i,l)\right)(\x^i-\y^{i-1})\right\|^2\notag\\
&\overset{a}\leq&   \left( \sum_{i=j(k)+1}^{k}\left(1+ \sum_{l=i}^{k}b(i,l)\right)\right) \sum_{i=j(k)+1}^{k}\left(1+ \sum_{l=i}^{k}b(i,l)\right) \left\|\x^i -\y^{i-1}\right\|^2 \notag\\
&\overset{b}\leq& \left( \sum_{i=j(k)+1}^{k}\left(1+  \sum_{l=1}^{k-i+1} 1 \right)\right) \sum_{i=j(k)+1}^{k}\left(1+ \sum_{l=1}^{k-i+1} 1 \right) \left\|\x^i -\y^{i-1}\right\|^2 \notag\\
&\overset{c}\leq& \left( \sum_{ii=1}^{k-j(k)}\left(1+  \sum_{l=1}^{ii} 1 \right)\right) \sum_{ii=1}^{k-j(k)}\left(1+ \sum_{l=1}^{ii} 1 \right)\left \|\x^{k-ii+1} -\y^{k-ii}\right\|^2 \notag\\
&\overset{d}\leq& \left( \sum_{ii=1}^{\min(\tau,k)}\left(1+  \sum_{l=1}^{ii} 1 \right)\right) \sum_{ii=1}^{\tau}\left(1+ \sum_{l=1}^{ii} 1 \right)\left \|\x^{k-ii+1} -\y^{k-ii}\right\|^2 \notag\\
&\leq& \frac{\tau^2 +3\tau}{2} \sum_{i=1}^{\min(\tau,k)}  (i+1)\left\|\x^{k-i+1} -\y^{k-i}\right\|^2,\notag
\end{eqnarray} 
where in $\overset{a}\leq$, we use the fact that for $ c_i\geq 0, \ 0\leq i\leq n,$,
\begin{eqnarray}\label{eee}
\left \|c_1 \a_1 + c_2 \a_2 +\cdots c_n \a_n\right\|^2 \leq (c_1 +c_2 +\cdots +c_n) (c_1\left\|\a_1\right\|^2 +c_2 \left\|\a_2\right\|^2 +\cdots c_n\left \|\a_n\right\|^2),
\end{eqnarray}
since the function $f(\x) = \|\x\|^2$ is convex, and so 
\begin{eqnarray}
&&\left\| \frac{c_1}{\sum_{i=0}^n c_i} \a_1+ \frac{c_2}{\sum_{i=0}^n c_i} \a_2+\cdots + \frac{c_n}{\sum_{i=0}^n c_i} \a_n \right \|^2\notag\\
&\leq&   \frac{c_1}{\sum_{i=0}^n c_i}\left \|\a_1\right\|^2+ \frac{c_2}{\sum_{i=0}^n c_i} \left\|\a_2\right\|^2+\cdots + \frac{c_n}{\sum_{i=0}^n c_i} \left\|\a_n \right \|^2     ; \notag
\end{eqnarray}
in $\overset{b}\leq$, we use $b(i,l)\leq 1$;  in $\overset{c}\leq$, we change variable $ii = k-i+1$; and in $\overset{d}\leq$, we use $k-j(k)\leq \tau$.

As we are more interested the limited case, namely $k$ is large.
We suppose at the first $\tau$ step, we run our algorithm in serial. Diving $(\theta^k)^2$ on Eq.~\eqref{onee} and summing the results with $k = 0$ to $K$, we have 
\begin{eqnarray}\label{im6}
&&\sum_{k=0}^K\frac{1}{(\theta^k)^2}\left\|\w^{j(k)} -\y^k\right\|^2\\
&=&\sum_{k=\tau}^K\frac{1}{(\theta^k)^2}\left\|\w^{j(k)} -\y^k\right\|^2\notag\\
&\leq&\frac{\tau^2 +3\tau}{2}\sum_{k=\tau}^K \sum_{i=1}^{\min(\tau,k-\tau)}  \frac{(i+1)}{(\theta^k)^2}\left\|\x^{k-i+1} -\y^{k-i}\right\|^2\notag\\
&\overset{a}\leq&\frac{\tau^2 +3\tau}{2}\sum_{k=\tau}^K \sum_{i=1}^{\min(\tau,k-\tau)}  \frac{4(i+1)}{(\theta^{k-i+1})^2}\left\|\x^{k-i+1} -\y^{k-i}\right\|^2\notag\\
&\overset{b}\leq&  (\tau^2 +3\tau)^2\sum_{k=\tau}^{K-1}  \frac{1}{(\theta^k)^2}\left\|\x^{k+1}-  \y^{k}\right\|^2,\notag\\
&\leq&  (\tau^2 +3\tau)^2\sum_{k=0}^{K}  \frac{1}{(\theta^k)^2}\left\|\x^{k+1}-  \y^{k}\right\|^2,\notag
\end{eqnarray}
where in $\overset{a}\leq$, we use that $ (k+j)^2 \leq 4 k^2$, since $k\geq \tau\geq j$, so  $\frac{1}{(\theta^k)^2} \leq  \frac{4}{(\theta^{k-i+1})^2}$ with $k\geq \tau$ and $i\leq \min(\tau,k-\tau)$; $\overset{b}\leq$ is because that  for each $\frac{1}{(\theta^k)^2}\left\|\x^{k}-  \y^{k-1}\right\|^2$ ($1\leq k \leq K$) there are most $\tau$ terms with coefficient from $8$ to $4(\tau+1)$.

Diving $(\theta^k)^2$ on both sides of Eq.~\eqref{im3}, and use $\mu=0$, we have
\begin{eqnarray}\label{im4}
\frac{F(\x^{k+1}) -F(\x^*)}{(\theta^k)^2}&\leq& \frac{(1-\theta^k)(F(\x^{k}) -F(\x^*))}{(\theta^{k})^2} -\frac{\gamma}{(\theta^k)^2} (\frac{1}{2}-\frac{\gamma L}{2})\left\| \frac{\x^{k+1}-\y^k}{\gamma} \right\|^2\\
&&+\frac{\gamma}{(\theta^k)^2}\left(\frac{\gamma^2 L^2}{2C_1}+\gamma L\right)\left\| \frac{\w^{j(k)} - \y^k}{\gamma} \right\|^2 + \frac{\gamma C_1}{2(\theta^k)^2}\left\|\frac{\x^{k+1}-\y^k }{\gamma}\right\|^2 \notag\\
&&+\frac{1}{2\gamma} \left\|  \z^{k} - \x^*\right\|^2-\frac{1}{2\gamma} \left\|  \z^{k+1} - \x^*\right\|^2\notag\\
&\overset{a}\leq& \frac{(F(\x^{k}) -F(\x^*))}{(\theta^{k-1})^2} -\frac{\gamma}{(\theta^k)^2} (\frac{1}{2}-\frac{\gamma L}{2})\left\| \frac{\x^{k+1}-\y^k}{\gamma} \right\|^2\notag\\
&&+\frac{\gamma}{(\theta^k)^2}\left(\frac{\gamma^2 L^2}{2C_1}+L\gamma\right)\left\| \frac{\w^{j(k)} - \y^k}{\gamma}\right \|^2 + \frac{\gamma C_1}{2(\theta^k)^2}\left\|\frac{\x^{k+1}-\y^k }{\gamma}\right\|^2 \notag\\
&&+\frac{1}{2\gamma} \left\|  \z^{k} - \x^*\right\|^2-\frac{1}{2\gamma} \left \|  \z^{k+1} - \x^*\right\|^2\notag,
\end{eqnarray} 
where in $\overset{a}\leq$, we use that $\frac{1-\theta^k}{(\theta^k)^2}\leq \frac{1}{(\theta^{k-1})^2} $ for $k\geq 1$. When $k=0$, we have $1-\theta^0=0$.

Summing Eq.~\eqref{im4} with $k$ from $0$ to $K$,  and applying Eq.~\eqref{im6}, we have that
\begin{eqnarray}
&&\frac{F(\x^{K+1}) -F(\x^*)}{(\theta^k)^2}\\
&\leq&-\sum_{k=0}^K\frac{\gamma}{(\theta^k)^2} (\frac{1}{2}-\frac{\gamma L}{2})\left\| \frac{\x^{k+1}-\y^k}{\gamma} \right\|^2\notag\\
&&+\sum_{k=0}^K\frac{\gamma}{(\theta^k)^2}\left(\frac{\gamma^2 L^2}{2C_1}+\gamma L \right)\left\| \frac{\w^{j(k)} - \y^k}{\gamma} \right\|^2 + \sum_{k=0}^K \frac{\gamma C_1}{2(\theta^k)^2}\left\|\frac{\x^{k+1}-\y^k }{\gamma}\right\|^2 \notag\\
&&+\frac{1}{2\gamma}\left \|  \z^{0} - \x^*\right\|^2-\frac{1}{2\gamma} \left \|  \z^{K+1} - \x^*\right\|^2\notag\\
&\leq&+\frac{1}{2\gamma} \left\|  \z^{0} - \x^*\right\|^2-\frac{1}{2\gamma} \left \|  \z^{K+1} - \x^*\right\|^2\notag\notag\\ 
&&- \left(\frac{1}{2}-\frac{\gamma L}{2}-\frac{C_1}{2}-\left(\frac{\gamma^2 L^2}{2C_1}+\gamma L\right)(\tau^2+3\tau)^2\right)\sum_{k=0}^K \frac{\gamma}{(\theta^k)^2}\left\| \frac{\x^{k+1}-\y^k}{\gamma} \right\|^2\notag.
\end{eqnarray}
Set $C_1=\gamma L$, we have that
$$   2\gamma L  + 3\gamma L(\tau^2+3\tau)^2\leq 1,  $$
So
\begin{eqnarray}
\frac{F(\x^{K+1}) -F(\x^*)}{(\theta^K)^2}+\frac{1}{2\gamma}  \left\|  \z^{K+1} - \x^*\right\|^2\leq \frac{1}{2\gamma}\left \|  \z^{0} - \x^*\right\|^2.
\end{eqnarray}

Now we consider the strongly convex case. In the following, we set $\theta=\theta^k$, and use $\theta^a$ to denote the $a$'s power of $\theta$, instead. Multiply  Eq.~\eqref{onee} with $(1-\theta)^{K-k}$, and summing the results with $k$ from $0$ to $K$, we have
\begin{eqnarray}\label{123}
&&\sum_{k=0}^K (1-\theta)^{K-k}\left\|\w^{j(k)} -\y^k\right\|^2\\
&\leq& \frac{\tau^2 +3\tau}{2}\sum_{k=0}^K \sum_{i=1}^{\min(\tau,k)}  (i+1)(1-\theta)^{K-k}\left\|\x^{k-i+1} -\y^{k-i}\right\|^2,\notag\\
&\leq& \frac{\tau^2 +3\tau}{2} \sum_{k=0}^K\sum_{i=1}^{\min(\tau,k)}  (1-\theta)^{-i}(i+1)(1-\theta)^{K-(k-i)}\left\|\x^{k-i+1} -\y^{k-i}\right\|^2,\notag\notag\\
&\leq& \frac{\tau^2 +3\tau}{2(1-\theta)^{\tau}} \sum_{k=0}^K\sum_{i=1}^{\min(\tau,k)}  (i+1)(1-\theta)^{K-(k-i)}\left\|\x^{k-i+1} -\y^{k-i}\right\|^2,\notag\notag\\
&\overset{a}\leq& \frac{(\tau^2 +3\tau)^2}{4(1-\theta)^{\tau}}\sum_{k=0}^{K-1}(1-\theta)^{K-i}\left\|\x^{i+1} -\y^{i}\right\|^2\notag\\
&\leq& \frac{(\tau^2 +3\tau)^2}{4(1-\theta)^{\tau}}\sum_{k=0}^{K}(1-\theta)^{K-i}\left\|\x^{i+1} -\y^{i}\right\|^2\notag,
\end{eqnarray}
where  $\overset{a}\leq$ is because that for each $(1-\theta)^{K-i}\left\|\x^{i+1} -\y^{i}\right\|^2$ ($1\leq k \leq K$) there are most $\tau$ terms with coefficient from $2$ to $\tau+1$, like Eq.~\eqref{im6}.

By arrange term on Eq.~\eqref{im3}, we have that
\begin{eqnarray}\label{112}
&&F(\x^{k+1}) -F(\x^*) + \left(\frac{\theta^2}{2\gamma}+\frac{\mu\theta}{2} \right)  \left\|  \z^{k+1} - \x^*\right\|^2\notag\\ &\leq&(1-\theta)\left(F(\x^k)- F(\x^*)+ \left(\frac{\theta^2}{2\gamma}+\frac{\mu\theta}{2} \right)  \left\|  \z^{k} - \x^*\right\|^2\right)\\\notag
&&-\gamma (\frac{1}{2}-\frac{\gamma L}{2}-\frac{C_1}{2})\left\| \frac{\x^{k+1}-\y^k}{\gamma} \right\|^2+\gamma\left(\frac{\gamma^2 L^2}{2C_1}+\gamma L\right)\left\| \frac{\w^{j(k)} - \y^k}{\gamma}\right \|^2, 
\end{eqnarray}
since we have set $\theta = \frac{-\gamma\mu +\sqrt{\gamma \mu^2 +4 \gamma \mu}}{2}$, which satisfies that
$$ \left(\frac{\theta^2}{2\gamma}+\frac{\mu\theta}{2} \right) (1-\theta) = \frac{\theta^2}{2\gamma},$$ 
solving it, we will have to solve $g(x)= x^2 + \mu\gamma x - \mu\gamma=0 $, we will have $\sqrt{\gamma\mu}/2\leq\theta\leq \sqrt{\gamma\mu}$, since $\gamma\mu\leq 1$.  For the assumption of $\gamma$, we have
\begin{eqnarray}\label{leqq}
9\gamma L\tau^2 \leq \frac{5}{2}\gamma L  + \gamma L(\tau^2+3\tau)^2\leq 1, 
\end{eqnarray}
We then consider that $\frac{1}{(1-\theta)^{\tau}}$, without loss of generality, we assume that $\tau \geq 2$, we have that 
\begin{eqnarray}\label{1.5}
\frac{1}{(1-\theta)^{\tau}}  \overset{a}\leq \frac{1}{(1-\sqrt{\gamma \mu})^{\tau}}  \overset{b}\leq \frac{1}{(1-\frac{1}{3\tau}\sqrt{ \mu/L})^{\tau}} \overset{c}\leq \frac{1}{(1-\frac{1}{3\tau})^{\tau}}\overset{d}\leq \frac{1}{(1-\frac{1}{3})^{1}} \leq \frac{3}{2},
\end{eqnarray}
where in $\overset{a}\leq$, we use $\theta \leq \sqrt{\gamma \mu}$; in $\overset{b}\leq$, we use Eq.~\eqref{leqq}; $\overset{c}\leq$, we use $\frac{\mu}{L}\leq 1$, and  $\overset{d}\leq$, we use the fact that function $g(x) = (1-\frac{x}{3})^{-x}$ is monotonous increasing when $x\in (0,1]$.

Multiply Eq.~\eqref{112} with $\theta^{K-k}$, and summing the result with $k$ from $0$ to $K$, we have that
\begin{eqnarray}
&&F(\x^{K+1}) -F(\x^*) + \left(\frac{\theta^2}{2\gamma}+\frac{\mu\theta}{2} \right) \left \|  \z^{K+1} - \x^*\right\|^2\\\
&\leq& (1-\theta)^{K+1}\left( F(\x^{0}) -F(\x^*) + \left(\frac{\theta^2}{2\gamma}+\frac{\mu\theta}{2} \right)  \left\|  \z^{0} - \x^*\right\|^2\right)\notag\\
&&-\gamma (\frac{1}{2}-\frac{\gamma L}{2}-\frac{C_1}{2})\sum_{i = 0}^K(1-\theta)^{K-k}\left\|\x^{k+1}-\y^k \right \|^2 +\gamma\left(\frac{\gamma^2 L^2}{2C_1}+\gamma L\right)\sum_{k=0}^K(1-\theta)^{K-k}\left\| \frac{\w^{j(k)} - \y^k}{\gamma} \right\|^2\notag\\
&\overset{a}\leq&(1-\theta)^{K+1}\left( F(\x^{0}) -F(\x^*) + \left(\frac{\theta^2}{2\gamma}+\frac{\mu\theta}{2} \right) \left \|  \z^{0} - \x^*\right\|^2\right)\notag\\
&&-\gamma \left(\frac{1}{2}-\frac{\gamma L}{2}-\frac{C_1}{2}  - \left(\frac{\gamma^2 L^2}{2C_1}+\gamma L\right)\left( \frac{(\tau^2 +3\tau)^2}{4(1-\theta)^{\tau}}\right ) \right)\sum_{i = 0}^K(1-\theta)^{K-k}\left\|\x^{k+1}-\y^k  \right\|^2\notag\\
&\overset{b}\leq&(1-\theta)^{K+1}\left( F(\x^{0}) -F(\x^*) + \left(\frac{\theta^2}{2\gamma}+\frac{\mu\theta}{2} \right) \left \|  \z^{0} - \x^*\right\|^2\right)\notag\\
&&-\gamma \left(\frac{1}{2}-\frac{\gamma L}{2}-\frac{C_1}{2}  - \left(\frac{\gamma^2 L^2}{2C_1}+\gamma L\right)\left( \frac{3(\tau^2 +3\tau)^2}{8}\right ) \right)\sum_{i = 0}^K(1-\theta)^{K-k}\left\|\x^{k+1}-\y^k \right \|^2\notag.
\end{eqnarray}
where in $\overset{a}\leq$, we use Eq.~\eqref{123}, and in $\overset{b}\leq$, we use Eq.~\eqref{1.5}. Setting $C_1$ to be $\frac{3}{2}\gamma L$, we have that
\begin{eqnarray}
1 - \frac{5}{2} \gamma L -(\frac{2}{3}+2)\gamma L \left( \frac{3(\tau^2 +3\tau)^2}{8}\right ) \geq 0,
\end{eqnarray} 
this is the result.

\subsection{AASCD}\label{AASCD-sec}
\begin{lemma}\label{lemm}
	Each $\x^k$ is a convex combination of $\z^0, \cdots,\z^k$, suppose	
	$\x^k \sum_{i=0}^k e_{k,i}\z^i,$, we have $e_{0,0}=1$, $e_{1,0} = 1-n\theta^0$, $e_{1,1} = n\theta$. And for $k>1$, we have
	\begin{equation}\label{ff}
	e_{k+1,i}=
	\left\{
	\begin{array}{rl}
	(1-\theta^k)e_{k,i},&\quad i\leq k-1\\
	n(1-\theta^k)\theta^{k-1}+\theta^k-n\theta^k,&\quad i=k\\
	n\theta^k,& \quad i = k+1.
	\end{array}
	\right.
	\end{equation}
	Supposing $\tth^{k+1} = \sum_{i=0}^k a_{k,i} h(\z^i)$, we have
	\begin{eqnarray}\label{ff2}
	\E_{i_k} (\tth^{k+1}) = (1-\theta^k)\tth^k +\theta^k\sum_{i_k=1}^n h_{i_k}(\z^{k+1}_{i_k}),
	\end{eqnarray}
	where $\E_{i_k}$ denote the random expectation is only taken on $i_k$ under the condition that $\x^k$, $\z^k$ is known. 
\end{lemma}
The proof of Lemma \ref{lemm} is directly taken from \citep{APCG,fercoq2015accelerated}.  For completeness, we provide a proof in the end of the section.

Because all proof only uses the Lipschitz coordinate constant, we use $L$ instead of $L_c$ to represent it for simply.\\
\bfseries Outline of the Proof: \mdseries \\
Step 1:  Set $\y^k = \theta^k\z^k +(1-\theta^k)\x^k$. 
Through the update rule, we have that
\begin{eqnarray}\label{step211}
&&\y^k\\
&=& \y^{j(k)}+\sum_{i=j(k)+1}^{k}(1+c^i)(\x^{i}-\y^{i-1})+\sum_{i=j(k)+1}^{k}\left(\sum_{l=i}^{k-1} c^ib(i+1,l)\right)(\x^i-\y^{i-1})\notag\\
&&+\left(\sum_{i=j(k)+1}^{k}b(j(k),i)\right)(\y^{j(k)}-\x^{j(k)-1}).\notag
\end{eqnarray}
Step 2: By analyzing the function value,  we have
\begin{eqnarray}\label{im21}
f(\x^{k+1}) &\leq& f(\y^k)-\gamma(1-\frac{\gamma L}{2})\left\|\frac{\x^{k+1}_{i_k}-\y^{k}_{i_k}}{\gamma} \right\|^2-\langle \bxi^k_{i_k}, \x^{k+1}_{i_k}-\y^k_{i_k}\rangle\notag\\
&&+\frac{\gamma L^2}{2C_2}\left\| \w_{i_k}^{j(k)} - \y_{i_k}^k \right\|^2 + \frac{\gamma C_2}{2}\left\|\frac{\x_{i_k}^{k+1}-\y_{i_k}^k }{\gamma}\right\|^2,  
\end{eqnarray}
Step 3: By analyzing the $\left\|\z^{k+1}-\x^*\right\|^2$,  we have
\begin{eqnarray}\label{im22}
&&\frac{n^2}{2\gamma}\E_{i_k}\left\|\theta^k \z^{k+1} -\theta^k \x^* \right\|^2\\
&=&\frac{n^2}{2\gamma}\left\|\theta^k \z^{k} -\theta^k \x^* \right\|^2+\frac{1}{2\gamma n}\sum_{i_k=1}^n\left\|\x^{k+1}_{i_k}-\y^k_{i_k} \right\|^2-\sum_{i_k=1}^n\langle \bxi^k_{i_k},  \theta^k\z_{i_k}^{k}-\theta^k\x_{i_k}^*  \rangle.\notag\\
&&+(1-\theta^k)f(\x^k)+\theta^k f(\x^*) -f(\y^k)+\langle \nabla f(\y^k)-\nabla f(\w^{j(k)}), \y^k- \w^{j(k)}\rangle\notag.
\end{eqnarray}
Step 4: Taking expectation on Eq.~\eqref{im21}, and adding and Eq.~\eqref{im22}, and simplifying,  we have
\begin{eqnarray}\label{im23}
&&\E_{i_k} f(\x^{k+1})+\E_{i_k} [\tth^{k+1}] -F(\x^*) +\frac{n^2(\theta^{k})^2+n\theta^k\mu\gamma}{2\gamma}\E_{i_k}\left\| \z^{k+1} -\x^* \right\|^2\notag\\
&\leq&(1-\theta^k) \left(f(\x^k)+ \tth^{k} -F(\x^*) \right)-\gamma(1-\frac{\gamma L}{2}-\frac{C_2}{2})\E_{i_k}\left\|\x^{k+1}-\y^{k}\right \|^2\notag \\
&&+\left(\frac{\gamma L^2}{2nC_2}+L\right)\left\| \w^{j(k)} - \y^k \right\|^2 +\frac{n^2(\theta^{k})^2+(n-1)\theta^k\mu\gamma}{2\gamma}\left\| \z^{k} - \x^* \right\|^2.
\end{eqnarray}
Step 5: we choose proper step size and obtain Theorem 2 in the paper.\\
\bfseries  Proof of step 1: \mdseries \\
Through the  step 4  and  step 5 in Algorithm 4 in the paper, we have
\begin{eqnarray}\label{23}
\theta^k\z^k = \y^k - (1-\theta^k)\x^k,
\end{eqnarray}
and 
\begin{eqnarray}\label{24}
n\theta^k\z^{k+1} = \x^{k+1}  -  (1-\theta^k)\x^k +(n-1)\theta^k \z^k.
\end{eqnarray}
We  have
\begin{eqnarray}\label{11e}
\x^{k+1} =\y^k +n\theta^k(\z^{k+1}-\z^k).
\end{eqnarray}
Multiplying Eq.~\eqref{23} with $(n-1)$, and adding with Eq.~\eqref{24}, we have
\begin{eqnarray}\label{25}
n\theta^k\z^{k+1} =\x^{k+1} - (1-\theta^k)\x^k +(n-1)\y^k -(n-1)(1-\theta^k)\x^k.
\end{eqnarray}
Eliminating $\z^k$ using Eq.~\eqref{25} and Eq.~\eqref{23},  for $k\geq 1$, we have
\begin{eqnarray}\label{yy}
\frac{\y^k-(1-\theta^k)\x^k}{\theta^k} = \frac{\x^{k}-(1-\theta^{k-1})\x^{k-1}+(n-1)\y^{k-1}-(n-1)(1-\theta^{k-1})\x^{k-1}}{n\theta^{k-1}}.
\end{eqnarray}
Computing out  $\y^k$ through Eq.~\eqref{yy}, we have
\begin{eqnarray}
\y^k &=& \x^k -\theta^k \x^k +\frac{\theta^k \x^{k}}{n\theta^{k-1}}-\frac{\theta^k(1-\theta^{k-1}) \x^{k-1}}{n\theta^{k-1}}\\
&&-\frac{(n-1)\theta^k(1-\theta^{k-1}) \x^{k-1}}{n\theta^{k-1}}+\frac{(n-1)\theta^k\y^{k-1}}{n\theta^{k-1}}\notag\\
&=&\x^k +\frac{\theta^k}{\theta^{k-1}}(\frac{1}{n}-\theta^{k-1})(\x^k-\y^{k-1})+\frac{\theta^k(1-\theta^{k-1})}{\theta^{k-1}}(\y^{k-1}-\x^{k-1}).\notag
\end{eqnarray}
Still, we set $a^k=\frac{\theta^k(1-\theta^{k-1})}{\theta^{k-1}}(\y^{k-1}-\x^{k-1})$ and $ b(l,k) =\prod_{i=l}^{k} a^{i}$, where $l\leq k $. Then by setting $c^k = \frac{\theta^k}{\theta^{k-1}}(\frac{1}{n}-\theta^{k-1})$, we have
\begin{eqnarray}\label{299}
\y^k &=&\x^k +c^k(\x^k -\y^{k-1})+a^k(\y^{k-1}-\x^{k-1})\\
&=&\y^{k-1} +(1+c^k)(\x^k -\y^{k-1})+a^k(\y^{k-1}-\x^{k-1})\notag\\
&=&\y^{k-1} +(1+c^k)(\x^k -\y^{k-1})+a^kc^{k-1}(\x^{k-1}-\y^{k-2})+a^ka^{k-1}(\y^{k-2}-\x^{k-2})\notag\\
&=&\y^{k-1} +(1+c^k)(\x^k -\y^{k-1})+\sum_{i = j(k)+1}^{k-1}b(i+1,k)c^i(\x^i-\y^{i-1})\notag\\
&&+b(j(k)+1,k)(\y^{j(k)}-\x^{j(k)-1}), \quad k\geq j(k)+1\geq 1. \notag
\end{eqnarray}
Like Eq.~\eqref{133}, summing Eq.~\eqref{299} with $k = j(k)+1$ to $k$, we have
\begin{eqnarray}
&&\y^k\\
&=& \y^{j(k)}+\sum_{i=j(k)+1}^{k}(1+c^i)(\x^{i}-\y^{i-1})+\sum_{l=j(k)+1}^{k}\sum_{i=j(k)+1}^{l-1} c^i b(i+1,l)(\x^i-\y^{i-1})\notag\\
&&+\left(\sum_{i=j(k)+1}^{k}b(j(k)+1,i)\right)(\y^{j(k)}-\x^{j(k)-1})\notag\\
&=& \y^{j(k)}+\sum_{i=j(k)+1}^{k}(1+c^i)(\x^{i}-\y^{i-1})+\sum_{i=j(k)+1}^{k}\left(\sum_{l=i}^{k-1} c^ib(i+1,l)\right)(\x^i-\y^{i-1})\notag\\
&&+\left(\sum_{i=j(k)+1}^{k}b(j(k),i)\right)(\y^{j(k)}-\x^{j(k)-1})\notag.
\end{eqnarray}
Comparing the result, we obtain step 1.\\
\bfseries  Proof of step 2: \mdseries \\
Through the optimal solution of $\z^{k+1}_{i_k}$ in step $4$, we have that 
\begin{eqnarray}\label{113}
n\theta^k(\z^{k+1}_{i_k} - \z^k_{i_k}) + \gamma \nabla_{i_k} f(\w^{j(k)})+\gamma \bxi^k_{i_k}=\mathbf{0},
\end{eqnarray}
where we denote $\bxi^k_{i_k}$ as a subgradient of $h_{i_k}$, i.e. $ \bxi^k_{i_k} \in h_{i_k}(\z^{k+1}_{i_k})$.
Through Eq.~\eqref{11e},
\begin{eqnarray}\label{212}
\x^{k+1}_{i_k} - \y^k_{i_k} + \gamma \nabla_{i_k} f(\w^{j(k)})+\gamma \bxi_{i_k}^k=\mathbf{0},
\end{eqnarray}
Since $f$ has Lipschitz continue gradient on coordinate $i_k$, we have
\begin{eqnarray}\label{166}
f(\x^{k+1}) &\leq& f(\y^k) +\langle\nabla_{i_k} f(\y^k), \x^{k+1}_{i_k}-\y^{k}_{i_k}\rangle +\frac{L}{2}\left\|\x^{k+1}_{i_k}-\y^{k}_{i_k} \right\|^2 \notag\\
&\overset{a}=&f(\y^k)  -\gamma\langle\nabla_{i_k} f(\y^k),\nabla_{i_k} f(\w^{j(k)})+ \bxi^k_{i_k}\rangle+\frac{L}{2}\left\|\x^{k+1}_{i_k}-\y^{k}_{i_k}\right \|^2 \notag\\
&\overset{b}=&f(\y^k)  -\gamma\langle\nabla_{i_k} f(\w^{j(k)})+ \bxi^k_{i_k},\nabla_{i_k} f(\w^{j(k)})+ \bxi^k_{i_k}\rangle+\frac{L}{2}\left\|\x^{k+1}_{i_k}-\y^{k}_{i_k} \right\|^2 \notag\\
&&+ \gamma \langle \bxi^k_{i_k}, \nabla_{i_k} f(\w^{j(k)})+ \bxi^k_{i_k} \rangle +\gamma \langle  \nabla_{i_k} f(\w^{j(k)}) -\nabla_{i_k} f(\y^k),\nabla_{i_k} f(\w^{j(k)})+ \bxi^k_{i_k} \rangle\notag\\
&=&f(\y^k)-\gamma(1-\frac{\gamma L}{2})\left\|\frac{\x^{k+1}_{i_k}-\y^{k}_{i_k}}{\gamma}\right \|^2-\langle \bxi_{i_k}^k, \x^{k+1}_{i_k}-\y^k_{i_k}\rangle\notag\\
&& -\langle    \nabla_{i_k} f(\w^{j(k)} ) - \nabla_{i_k} f(\y^k),  \x_{i_k}^{k+1}-\y_{i_k}^k   \rangle\notag\\
&\overset{c}\leq&f(\y^k)-\gamma(1-\frac{\gamma L}{2})\left\|\frac{\x^{k+1}_{i_k}-\y^{k}_{i_k}}{\gamma}\right \|^2-\langle \bxi^k_{i_k}, \x^{k+1}_{i_k}-\y^k_{i_k}\rangle\notag\\
&&+\frac{\gamma L^2}{2C_2}\left\| \w_{i_k}^{j(k)} - \y_{i_k}^k \right\|^2 + \frac{\gamma C_2}{2}\left\|\frac{\x_{i_k}^{k+1}-\y_{i_k}^k }{\gamma}\right\|^2, 
\end{eqnarray}
where in $\overset{a}\leq$, we use Eq.~\eqref{212}; in $\overset{b}\leq$, we  insert $\nabla_{i_k} f(\y^k)+\bxi^k_{i_k}$; in $\overset{c}\leq$, we use Cauchy-Schwarzwe inequality.\\
\bfseries  Proof of step 3: \mdseries 
\begin{eqnarray}\label{167}
&&\frac{n^2}{2\gamma}\left\|\theta^k \z^{k+1} -\theta^k \x^*\right \|^2\\
&=& \frac{n^2}{2\gamma}\left\| \theta^k \z^{k} -\theta^k \x^*+ \theta^k \z^{k+1}-\theta^k\z^k \right \|^2\notag\\
&=&\frac{n^2}{2\gamma} \left\| \theta^k \z^{k} -\theta^k \x^*\right\|^2 +\frac{n^2}{2\gamma}\left \| \theta^k \z^{k+1}_{i_k}-\theta^k\z^k_{i_k} \right \|^2 +\frac{n^2}{\gamma}\langle \theta^k \left(\z^{k+1}_{i_k}-\z^k_{i_k}\right) ,  \theta^k\z^{k}_{i_k}-\theta^k\x^*_{i_k} \rangle\notag\\
&\overset{\eqref{113}}=&\frac{n^2}{2\gamma}\left\| \theta^k \z^{k}_{i_k} -\theta^k \x^*_{i_k}\right\|^2 +\frac{1}{2\gamma }\left\|\x^{k+1}_{i_k}-\y^k_{i_k} \right\|^2-n\langle \nabla_{i_k} f(\w^{j(k)})+\bxi^k_{i_k},  \theta^k\z_{i_k}^{k}-\theta^k\x_{i_k}^*  \rangle.\notag
\end{eqnarray}

So taking expectation on Eq.~\eqref{167}, we have
\begin{eqnarray}\label{199}
&&\frac{n^2}{2\gamma}\E_{i_k}\left\|\theta^k \z^{k+1} -\theta^k \x^* \right\|^2\\
&=&\frac{n^2}{2\gamma}\left\|\theta^k \z^{k} -\theta^k \x^* \right\|^2+\frac{1}{2\gamma n}\sum_{i_k=1}^n\left\|\x^{k+1}_{i_k}-\y^k_{i_k}\right \|^2\notag\\
&&-\langle \nabla f(\w^{j(k)}),  \theta^k\z^{k}-\theta^k\x^*  \rangle-\sum_{i_k=1}^n\langle \bxi^k_{i_k},  \theta^k\z_{i_k}^{k}-\theta^k\x_{i_k}^*  \rangle.\notag
\end{eqnarray}
By the same technology of  Eq.~\eqref{77},  for the last second term of  Eq.~\eqref{199},  we have that 
\begin{eqnarray}\label{201}
&&-\langle   \nabla f(\w^{j(k)}),  \theta^k\z^{k}-\theta^k\x^* \rangle\\
&=&-\langle  \nabla   f(\w^{j(k)}),  \y^k-(1-\theta^k)\x^k-\theta^k\x^* \rangle\notag\\
&=&-\langle  \nabla  f(\w^{j(k)}),  \w^{j(k)}-(1-\theta^k)\x^k-\theta^k\x^* \rangle -\langle\nabla f(\w^{j(k)}), \y^k- \w^{j(k)}\rangle\notag\\
&\leq& (1-\theta^k)f(\x^k)+\theta^k f(\x^*) -f(\w^{j(k)})-\langle \nabla f(\w^{j(k)}), \y^k- \w^{j(k)}\rangle\notag\\
&\leq& (1-\theta^k)f(\x^k)+\theta^k f(\x^*) -f(\w^{j(k)})-\langle \nabla f(\w^{j(k)}), \y^k- \w^{j(k)}\rangle\notag\\
&\leq& (1-\theta^k)f(\x^k)+\theta^k f(\x^*) -f(\y^k)+\langle \nabla f(\y^k)-\nabla f(\w^{j(k)}), \y^k- \w^{j(k)}\rangle\notag,
\end{eqnarray}
Substituting Eq.~\eqref{201} into Eq.~\eqref{199}, we have the results of Step 3.\\
\bfseries  Proof of step 4: \mdseries \\
Taking expectation on Eq.~\eqref{im21}, we have
\begin{eqnarray}\label{200}
\E_{i_k} f(\x^{k+1}) &\leq& f(\y^k)-\gamma(1-\frac{\gamma L}{2}-\frac{C_2}{\gamma})\frac{1}{n}\sum_{i=1}^n\left\|\frac{\x^{k+1}_{i_k}-\y^{k}_{i_k}}{\gamma} \right\|^2\\
&&-\frac{1}{n}\sum_{i=1}^n \langle\bxi^k_{i_k}, \x^{k+1}_{i_k}-\y^k_{i_k}\rangle+\left(\frac{\gamma L^2}{2nC_2}\right)\left\| \w^{j(k)} - \y^k\right \|^2.\notag
\end{eqnarray}
Adding Eq.~\eqref{200} and Eq.~\eqref{im22}  we have
\begin{eqnarray}
\E_{i_k} f(\x^{k+1}) &\leq& (1-\theta^k)f(\x^k)+\theta^k f(\x^*)-\frac{\gamma}{n}\sum_{i_k=1}^n(\frac{1}{2}-\frac{\gamma L}{2}-\frac{C_2}{2})\left\|\frac{\x^{k+1}_{i_k}-\y^{k}_{i_k}}{\gamma}\right \|^2 \notag\\
&&+\left(\frac{\gamma L^2}{2nC_2}+L\right)\left\| \w^{j(k)} - \y^k\right \|^2 -\sum_{i_k=1}^n\langle   \bxi^k_{i_k}, \theta^k\z_{i_k}^{k}-\theta^k\x_{i_k}^*+\frac{1}{n} \left(\x^{k+1}_{i_k}-\y^k_{i_k}\right) \rangle\notag\\
&&+\frac{n^2}{2\gamma}\left\|\theta^k \z^{k} -\theta^k \x^*\right \|^2-\frac{n^2}{2\gamma}\E_{i_k}\left\|\theta^k \z^{k+1} -\theta^k \x^*\right \|^2\notag\\
&\overset{a}=&  (1-\theta^k)f(\x^k)+\theta^k f(\x^*)-\frac{\gamma}{n}\sum_{i_k=1}^n(1-\frac{\gamma L}{2}-\frac{C_2}{2})\left\|\frac{\x^{k+1}_{i_k}-\y^{k}_{i_k}}{\gamma} \right\|^2\notag \\
&&+\left(\frac{\gamma L^2}{2nC_2}+L\right)\left\| \w^{j(k)} - \y^k\right \|^2 -\sum_{i_k=1}^n \langle   \bxi^k_{i_k}, \theta^k\z_{i_k}^{k+1}-\theta^k\x_{i_k}^*\rangle\notag\\
&&+\frac{n^2}{2\gamma}\left\|\theta^k \z^{k} -\theta^k \x^* \right\|^2-\frac{n^2}{2\gamma}\E_{i_k}\left\|\theta^k \z^{k+1} -\theta^k \x^*\right \|^2
\end{eqnarray}
where in $\overset{a}=$, we use Eq.~\eqref{11e}.

The same us Eq.~\eqref{88}, since $h_{i_k}$ is convex, we have
\begin{eqnarray}
\theta^k\langle \bxi^k_{i_k} ,\x^*_{i_k}-\z^{k+1}_{i_k}\rangle \leq \theta^k h_{i_k}(\x^*) -\theta^k h_{i_k}(\z^{k+1}_{i_k})-\frac{\mu\theta^k}{2}\left\|\z_{i_k}^{k+1}-\x_{i_k}^* \right\|^2.
\end{eqnarray}
Analyzing the expectation, we have
\begin{eqnarray}\label{ppp}
\E_{i_k} \left\|\z^{k+1}-\x^* \right \|^2 &=& \frac{1}{n}\sum_{i_k=1}^n\left[ \left \|\z^{k+1}_{i_k}-\x^*_{i_k} \right \|^2  + \sum_{j\neq i_k}  \left\|\z^{k}_{j}-\x^*_{i_k} \right \|^2 \right] \notag\\
&=&\frac{1}{n}\sum_{i_k=1}^n\left\|\z^{k+1}_{i_k}-\x^*_{i_k} \right \|^2+\frac{n-1}{n} \left\|\z^{k}-\x^*\right \|^2,
\end{eqnarray}
Since as Eq.~\eqref{ppp}, we can find that
\begin{eqnarray}\label{xk}
\E_{i_k}\left \|\x^{k+1}-\y^k  \right\|^2 &=& \frac{1}{n}\sum_{i_k=1}^n\left[  \left\|\x^{k+1}_{i_k}-\y^*_{i_k}  \right\|^2  + \sum_{j\neq i_k}  \left\|\x^{k}_{j}-\y^k_{i_k}  \right\|^2 \right] \notag\\
&\overset{\eqref{11e}}=&\frac{1}{n}\sum_{i_k=1}^n\left\|\x^{k+1}_{i_k}-\y^k_{i_k}\right\|^2.
\end{eqnarray}
We have that
\begin{eqnarray}
&&\E_{i_k} f(\x^{k+1})+\sum_{i_k=1}^nh(\z^{k+1}_{i_k})\\
&\leq&(1-\theta^k)f(\x^k)+\theta^k F(\x^*)-\frac{\gamma}{n}\sum_{i_k=1}^n(\frac{1}{2}-\frac{\gamma L}{2}-\frac{C_2}{2})\left\|\frac{\x^{k+1}_{i_k}-\y^{k}_{i_k}}{\gamma} \right\|^2\notag \\
&&+\left(\frac{\gamma L^2}{2nC_2}+L\right)\left\| \w^{j(k)} - \y^k\right \|^2 -\sum_{i_k=1}^n\frac{\mu\theta^k}{2}\left\|\z_{i_k}^{k+1}-\x_{i_k}^* \right\|^2\notag\\
&&+\frac{n^2}{2\gamma}\left\|\theta^k \z^{k} -\theta^k \x^* \right\|^2-\frac{n^2}{2\gamma}\E_{i_k}\left\|\theta^k \z^{k+1} -\theta^k \x^* \right\|^2\notag\\
&\overset{a}=&(1-\theta^k)f(\x^k)+\theta^k F(\x^*)-\frac{\gamma}{n}\sum_{i_k=1}^n(\frac{1}{2}-\frac{\gamma L}{2}-\frac{C_2}{2})\left\|\x^{k+1}_{i_k}-\y^{k}_{i_k} \right\|^2\notag \\
&&+\left(\frac{\gamma L^2}{2nC_2}+L\right)\left\| \w^{j(k)} - \y^k \right\|^2 +\frac{n^2(\theta^{k})^2+(n-1)\theta\mu\gamma}{2\gamma}\left\| \z^{k} - \x^* \right\|^2\notag\\
&&-\frac{n^2(\theta^{k})^2+n\theta\mu\gamma}{2\gamma}\E_{i_k}\left\| \z^{k+1} -\x^* \right\|^2\notag,
\end{eqnarray}
where in $\overset{a}=$, we use Eq.~\eqref{ppp}. We obtain the results of Step 4.\\
\bfseries  Proof of step 5: \mdseries \\
Through Eq.~\eqref{step211},  using the same technique of Eq.~\eqref{onee}, we have
\begin{eqnarray}\label{2onee}
&&\left\|\w^{j(k)} -\y^k\right\|^2\\
&=&\left\|\sum_{i=j(k)+1}^{k}\left(1+c^i+c^i\sum_{l=i}^{k-1} b(i,l)\right)(\x^i-\y^{i-1})\right\|^2\notag\\
&\overset{a}\leq&   \left( \sum_{i=j(k)+1}^{k}\left(1+c^i +c^i \sum_{l=i}^{k-1}b(i,l)\right)\right) \sum_{i=j(k)+1}^{k}\left(1+c^i +c^i \sum_{l=i}^{k-1}b(i,l)\right)\left \|\x^i -\y^{i-1}\right\|^2 \notag\\
&\overset{b}\leq& \left( \sum_{i=j(k)+1}^{k}\left(1+  c^i \sum_{l=1}^{k-i+1} 1 \right)\right) \sum_{i=j(k)+1}^{k}\left(1+ c^i\sum_{l=1}^{k-i+1} 1 \right) \left\|\x^i -\y^{i-1}\right\|^2 \notag\\
&\overset{c}\leq& \left( \sum_{ii=1}^{k-j(k)}\left(1+  \frac{1}{n}\sum_{l=1}^{ii} 1 \right)\right) \sum_{ii=1}^{k-j(k)}\left(1+\frac{1}{n} \sum_{l=1}^{ii} 1 \right) \left\|\x^{k-ii+1} -\y^{k-ii}\right\|^2 \notag\\
&\overset{d}\leq& \left( \sum_{ii=1}^{\min(\tau,k)}\left(1+ \frac{1}{n} \sum_{l=1}^{ii} 1 \right)\right) \sum_{ii=1}^{\tau}\left(1+ \frac{1}{n}\sum_{l=1}^{ii} 1 \right) \left\|\x^{k-ii+1} -\y^{k-ii}\right\|^2 \notag\\
&\leq& \left(\frac{\tau^2 +\tau}{2n}+\tau\right) \sum_{i=1}^{\min(\tau,k)}  (\frac{i}{n}+1)\left\|\x^{k-i+1} -\y^{k-i}\right\|^2,\notag
\end{eqnarray}
where in $\overset{a}\leq$, we use Eq.~\eqref{eee}; in $\overset{b}\leq$, we use $b(i,l)\leq 1$;  in $\overset{c}\leq$, we change variable $ii = k-i+1$, and use $\c^k\leq \frac{1}{n}$; and in $\overset{d}\leq$, we use $k-j(k)\leq \tau$. Since $\tau \leq \sqrt{n}\leq 2n$. Diving $(\theta^k)^2$ on Eq.~\eqref{onee} and summing the results with $k = 0$ to $K$, we have 
\begin{eqnarray}\label{2im6}
&&\sum_{k=0}^K\frac{1}{(\theta^k)^2}\left\|\w^{j(k)} -\y^k\right\|^2\\
&\overset{a}\leq&\left(\frac{\tau^2+\tau}{2n} +\tau \right)\sum_{k=0}^K \sum_{i=1}^{\min(\tau,k-\tau)}  \frac{4(\frac{i}{n}+1)}{(\theta^{k-i+1})^2}\left\|\x^{k-i+1} -\y^{k-i}\right\|^2\notag\\
&\overset{b}\leq&  \left(\frac{\tau^2+\tau}{n} +2\tau \right)^2\sum_{k=0}^{K-1}  \frac{1}{(\theta^k)^2}\left\|\x^{k+1}-  \y^{k}\right\|^2,\notag\\
&\leq& \left(\frac{\tau^2+\tau}{n} +2\tau \right)^2\sum_{k=0}^{K}  \frac{1}{(\theta^k)^2}\left\|\x^{k+1}-  \y^{k}\right\|^2,\notag
\end{eqnarray}
where in $\overset{a}\leq$, we use that $ (k+j)^2 \leq 4 k^2$, when $k\geq 2n\geq \tau$ and $j\leq \tau$,     and $\frac{1}{\theta^k} = \frac{(2n+k)}{2}$ for $k\geq 0$;  $\overset{b}\leq$ is because that  for each $\frac{1}{(\theta^k)^2}\left\|\x^{k}-  \y^{k-1}\right\|^2$ ($1\leq k \leq K$) there are most $\tau$ terms.

Diving $(\theta^k)^2$ on both sides of Eq.~\eqref{im23}, and use $\mu=0$, we have
\begin{eqnarray}\label{im24}
&&\frac{\E_{i_k} f(\x^{k+1})+\E_{i_k} [\tth^{k+1}] -F(\x^*)}{(\theta^k)^2} +\frac{n^2}{2\gamma}\E_{i_k}\left\| \z^{k+1} -\x^* \right\|^2\notag\\
&\leq&\frac{1-\theta^k}{(\theta^k)^2} \left(f(\x^k)+ \tth^{k} -F(\x^*) \right)-\frac{\gamma}{(\theta^k)^2}(\frac{1}{2}-\frac{\gamma L}{2}-\frac{C_2}{2})\E_{i_k}\left\|\x^{k+1}-\y^{k} \right\|^2\notag \\
&&+\left(\frac{\gamma L^2}{2nC_2}+L\right)\left\| \w^{j(k)} - \y^k\right \|^2 +\frac{n^2}{2\gamma}\left\| \z^{k} - \x^* \right\|^2\notag\\
&\overset{a}\leq&\frac{1}{(\theta^{k-1})^2} \left(f(\x^k)+ \tth^{k} -F(\x^*) \right)-\frac{\gamma}{(\theta^k)^2}(1-\frac{\gamma L}{2}-\frac{C_2}{2})\E_{i_k}\left\|\x^{k+1}-\y^{k}\right\|^2\notag \\
&&+\left(\frac{\gamma L^2}{2nC_2}+L\right)\left\| \w^{j(k)} - \y^k\right \|^2 +\frac{n^2}{2\gamma}\left\| \z^{k} - \x^*\right \|^2,
\end{eqnarray}
where in $\overset{a}\leq$, we use that $\frac{1-\theta^k}{(\theta^k)^2}\leq \frac{1}{(\theta^{k-1})^2} $.

Taking expectation on the first $k$ iteration for Eq.~\eqref{im24}, and summing it with $k$ from $0$ to $K$,   we have that
\begin{eqnarray}
&&\frac{\E f(\x^{K+1})+\E [\tth^{K+1}] -F(\x^*)}{(\theta^K)^2} +\frac{n^2}{2\gamma}\E\left\| \z^{K+1} -\x^* \right\|^2\\
&\leq&\frac{f(\x^0)+ \tth^{0} -F(\x^*)}{(\theta^{-1})^2} -\frac{\gamma}{(\theta^k)^2}(\frac{1}{2}-\frac{\gamma L}{2}-\frac{C_2}{2})\sum_{k=0}^K\E\left\|\x^{k+1}-\y^{k}\right \|^2\notag \\
&&+\left(\frac{\gamma L^2}{2nC_2}+L\right)\sum_{k=0}^K\E\left\| \w^{j(k)} - \y^k \right\|^2 +\frac{n^2}{2\gamma}\left\| \z^0 - \x^*\right \|^2\notag\\
&\overset{\eqref{2im6}}\leq&\frac{f(\x^0)+ \tth^{0} -F(\x^*)}{(\theta^{-1})^2} +\frac{n^2}{2\gamma}\left\| \z^0 - \x^* \right\|^2\notag \\
&&- \left(\frac{1}{2}-\frac{\gamma L}{2}-\frac{C_2}{2}-\left(\frac{\gamma^2 L^2}{2nC_2}+\gamma L\right)\left(\frac{\tau^2+\tau}{n} +2\tau \right)^2\right)\sum_{k=0}^K\frac{\gamma}{(\theta^k)^2}\E\left\| \frac{\x^{k+1}-\y^k}{\gamma} \right\|^2.\notag
\end{eqnarray}

Set $C_2=\gamma L$, we have that
$$   2\gamma L  + (1+\frac{1}{n})\gamma L \left(\frac{\tau^2+\tau}{n} +2\tau \right)^2\leq 1,  $$
So
\begin{eqnarray}
&&\frac{\E F(\x^{K+1}) -F(\x^*)}{(\theta^K)^2} +\frac{n^2}{2\gamma}\E\left\| \z^{K+1} -\x^* \right\|^2 \notag\\
&\overset{a}\leq&\frac{\E f(\x^{K+1})+\E [\tth^{K+1}] -F(\x^*)}{(\theta^K)^2} +\frac{n^2}{2\gamma}\E\left\| \z^{K+1} -\x^* \right\|^2 \notag\\
&\leq&\frac{f(\x^0)+ \tth^{0} -F(\x^*)}{(\theta^{-1})^2} +\frac{n^2}{\gamma^2}\left\|\z^0 -\x^*\right\|^2\notag\\
&\overset{b}\leq&\frac{F(\x^0)-F(\x^*)}{(\theta^{-1})^2} +\frac{n^2}{2\gamma}\left\|\z^0 -\x^*\right\|^2,
\end{eqnarray}
where in $\overset{a}\leq$, we use $h(\x^{K+1})=h(\sum_{i=0}^{K+1}e_{k+1,i}\z^i)\leq\sum_{i=0}^{K+1}e_{k+1,i}h(\z^i)=\tth^{K+1}$; and in $\overset{b}\leq$, we use $h(\x^0) = \tth^0$.

Now we consider the strongly convex case. In the following, again we set $\theta=\theta^k$, and use $\theta^a$ to denote the $a$'s power of $\theta$. Multiply  Eq.~\eqref{2onee} with $(1-\theta)^{K-k}$, and summing the results with $k$ from $0$ to $K$, we have
\begin{eqnarray}\label{1123}
&&\sum_{k=0}^K (1-\theta)^{K-k}\left\|\w^{j(k)} -\y^k\right\|^2\\
&\leq&\left(\frac{\tau^2 +\tau}{2n}+\tau\right)\sum_{k=0}^K \sum_{i=1}^{\min(\tau,k)}  (\frac{i}{n}+1)(1-\theta)^{K-k}\left\|\x^{k-i+1} -\y^{k-i}\right\|^2,\notag\\
&\leq& \left(\frac{\tau^2 +\tau}{2n}+\tau\right) \sum_{k=0}^K\sum_{i=1}^{\min(\tau,k)}  (1-\theta)^{-i}(\frac{i}{n}+1)(1-\theta)^{K-(k-i)}\left\|\x^{k-i+1} -\y^{k-i}\right\|^2,\notag\notag\\
&\leq& \frac{1}{(1-\theta)^{\tau}}\left( \frac{\tau^2+\tau}{2n}+\tau\right) \sum_{k=0}^K\sum_{i=1}^{\min(\tau,k)}  (\frac{i}{n}+1)(1-\theta)^{K-(k-i)}\left\|\x^{k-i+1} -\y^{k-i}\right\|^2,\notag\notag\\
&\overset{a}\leq& \frac{\left((\tau^2+\tau)/n +2\tau\right)^2}{4(1-\theta)^{\tau}}\sum_{k=0}^{K-1}(1-\theta)^{K-i}\left\|\x^{i+1} -\y^{i}\right\|^2\notag,
\end{eqnarray}
where  $\overset{a}\leq$ is because that for each $(1-\theta)^{K-i}\left\|\x^{i+1} -\y^{i}\right\|^2$ ($1\leq k \leq K$) there are most $\tau$ terms, like Eq.~\eqref{im6}.

By arrange term on Eq.~\eqref{im23}, we have that
\begin{eqnarray}\label{2112}
&&\E_{i_k} f(\x^{k+1})+\E_{i_k} [\tth^{k+1}] -F(\x^*) +\frac{n^2(\theta)^2+n\theta\mu\gamma}{2\gamma}\E_{i_k}\left\| \z^{k+1} -\x^* \right\|^2\\
&\leq&(1-\theta) \left(f(\x^k)+ \tth^{k} -F(\x^*) +\frac{n^2(\theta)^2+n\theta\mu\gamma}{2\gamma}\left\| \z^{k} -\x^* \right\|^2 \right)\notag\\
&&-\gamma(\frac{1}{2}-\frac{\gamma L}{2}-\frac{C_2}{2})\E_{i_k}\left\|\x^{k+1}-\y^{k} \right\|^2
+\left(\frac{\gamma L^2}{2nC_2}+L\right)\left\| \w^{j(k)} - \y^k \right\|^2.\notag
\end{eqnarray}
since we have set $\theta = \frac{-\gamma\mu +\sqrt{\gamma \mu^2 +4 \gamma \mu}}{2n}$, which satisfies that
$$ \left(\frac{\theta^2n^2}{2\gamma}+\frac{n\mu\theta}{2} \right) (1-\theta) = \frac{\theta^2}{2\gamma}+\frac{(n-1)\mu\theta}{2},$$ 
solving it, we will have to solve $g(x)= n^2x^2 + n\mu\gamma x - \mu\gamma=0 $, we assume $\mu /L \leq n^2$, and we will have $\sqrt{\gamma\mu}/2\leq n\theta\leq \sqrt{\gamma\mu}$.  For the assumption of $\gamma$, we have
\begin{eqnarray}\label{2leqq}
3\gamma L\tau^2 \leq 2\gamma L  + (\frac{3}{4}+\frac{3}{8n})\gamma L\left((\tau^2+\tau)/n +2\tau\right)^2\leq 1, 
\end{eqnarray}
We then consider that $\frac{1}{(1-\theta)^{\tau}}$, without loss of generality, we assume that $n \geq 2$, we have that 
\begin{eqnarray}\label{21.5}
\frac{1}{(1-\theta)^{\tau}}  \overset{a}\leq \frac{1}{(1-\sqrt{\gamma \mu}/n)^{\tau}}  \overset{b}\leq \frac{1}{(1-\frac{1}{\tau}\sqrt{ \mu/(3L)}/n)^{\tau}} \overset{c}\leq \frac{1}{(1-\frac{1}{2\sqrt{3}\tau})^{\tau}}\leq \frac{1}{(1-\frac{1}{2\sqrt{3}})} \leq \frac{3}{2},
\end{eqnarray}
where in $\overset{a}\leq$, we use $\theta \leq \sqrt{\gamma \mu}$; in $\overset{b}\leq$, we use Eq.~\eqref{leqq}; $\overset{c}\leq$, we use $\sqrt{\frac{\mu}{L}}/{n}\leq \frac{1}{2}$.

Taking expectation Eq.~\eqref{2112}, and  Multiply Eq.~\eqref{2112} with $\theta^{K-k}$,  then and summing the result with $k$ from $0$ to $K$, we have that
\begin{eqnarray}
&&\E f(\x^{K+1})+\E [\tth^{K+1}] -F(\x^*) +\frac{n^2(\theta)^2+n\theta\mu\gamma}{2\gamma}\E\left\| \z^{K+1} -\x^*\right \|^2\\
&\leq&(1-\theta)^{K+1} \left(f(\x^0)+ \tth^{k} -F(\x^*) +\frac{n^2(\theta)^2+n\theta\mu\gamma}{2\gamma}\left\| \z^{0} -\x^* \right\|^2 \right)\notag\\
&&-\gamma(\frac{1}{2}-\frac{\gamma L}{2}-\frac{C_2}{2})\sum_{i = 0}^K(1-\theta)^{K-k}\E\left\|\x^{k+1}-\y^{k}\right\|^2
+\left(\frac{\gamma L^2}{2nC_2}+L\right)\sum_{i = 0}^K(1-\theta)^{K-k}E\left\| \w^{j(k)} - \y^k \right\|^2\notag\\
&\leq&(1-\theta)^{K+1} \left(f(\x^0)+ \tth^{k} -F(\x^*) +\frac{n^2(\theta)^2+n\theta\mu\gamma}{2\gamma}\left\| \z^{0} -\x^* \right\|^2 \right)\notag\\
&&-\gamma\left(\frac{1}{2}-\frac{\gamma L}{2}-\frac{C_2}{2}-\left(\frac{\gamma^2 L^2}{2nC_2}+\gamma L\right)\left(\frac{\left((\tau^2+\tau)/n +2\tau\right)^2}{4(1-\theta)^{\tau}} \right)\right)\sum_{i = 0}^K(1-\theta)^{K-k}\E\left\|\x^{k+1}-\y^{k} \right\|^2.\notag
\end{eqnarray}
Set $C_2=\gamma L$,  we have 
$$ 2\gamma L+ \left(\frac{\gamma L}{n}+2\gamma L\right)\left(\frac{\left((\tau^2+\tau)/n +2\tau\right)^2}{4(1-\theta)^{\tau}} \right) \leq  2\gamma L+ \left(\frac{\gamma L}{n}+2\gamma L\right)\left(\frac{3\left((\tau^2+\tau)/n +2\tau\right)^2}{8} \right)\leq 1$$.

Then using   $h(\x^{K+1})\leq \tth^{K+1}$ and $h(\x^0) = \tth^0$, we obtain the results.
\begin{eqnarray}
&&\E [F(\x^{K+1})] -F(\x^*) +\frac{n^2(\theta)^2+n\theta\mu\gamma}{2\gamma}\E\left\| \z^{K+1} -\x^*\right \|^2\\
&\leq&(1-\theta)^{K+1} \left(F(\x^0) -F(\x^*) +\frac{n^2(\theta)^2+n\theta\mu\gamma}{2\gamma}\left\| \z^{0} -\x^*\right \|^2 \right).\notag
\end{eqnarray}
\bfseries Proof of Lemma~\ref{lemm}, taken from \citep{APCG,fercoq2015accelerated}. \mdseries \\
We proof $e_{k+1,i}$ first. When $k = 0$ and $1$, it is right. We then proof Eq.~\eqref{ff}.
Since
\begin{eqnarray}
\x^{k+1} &=& (1-\theta^k)\x^k +\theta^k \z^k +n\theta^k (\z^{k+1}-\z^k)\\
&=&(1-\theta^k)\sum_{i=0}^k e_{k,i}\z^i +\theta^k \z^k +n\theta^k(\z^{k+1}-\z^{k})\notag\\
&=&(1-\theta^k)\sum_{i=0}^{k-1} e_{k,i}\z^i +\left((1-\theta^k)e_{k,k}+\theta^k -n\theta^k\right)\z^k +n\theta^k \z^{k+1}.\notag
\end{eqnarray}
Comparing the results, we obtain Eq.~\eqref{ff}. For Eq.~\eqref{ff2}, we have
\begin{eqnarray}
\E_{i_k}[\tth^{k+1}] &\overset{a}=& \sum_{i=0}^k h(\z^i) +\E_{i_k} n\theta^k h(\z^{k+1})\notag\\
&=& \sum_{i=0}^ke_{k+1,i} h(\z^i) +\frac{1}{n}\sum_{i_k} n\theta^k \left(h_{i_k}(\z^{k+1}_{i_k})+\sum_{j\neq i_k}h_j(\z^k) \right)\notag\\
&=& \sum_{i=0}^ke_{k+1,i} h (\z^i) +\theta^k\sum_{i_k} h_{i_k}(\z^{k+1}_{i_k})+(n-1)\theta^k h(\z^k) \notag\\
&\overset{b}=& \sum_{i=0}^{k-1} e_{k+1,i} h(\z^i) +(n(1-\theta^k)\theta^{k-1}+\theta^k-n\theta^k)h(\z^k)+(n-1)\theta^kh(\z^k) +\theta^k\sum_{i_k} h_{i_k}(\z^{k+1}_{i_k})\notag\\
&=& \sum_{i=0}^{k-1}e_{k+1,i}  h(\z^i)+n(1-\theta^k)\theta^{k-1}h(\z^k)+\theta^k\sum_{i_k} h_{i_k}(\z^{k+1}_{i_k})\notag\\
&\overset{c}=& \sum_{i=0}^{k-1}e_{k,i} (1-\theta^k) h(\z^i) +(1-\theta^k)e_{k,k}h(\z^k)+\theta^k\sum_{i_k} h_{i_k}(\z^{k+1}_{i_k})\notag\\
&=& \sum_{i=0}^{k}e_{k,i} (1-\theta^k) h(\z^i) +\theta^k\sum_{i_k} h_{i_k}(\z^{k+1}_{i_k})=\tth^k +\theta^k\sum_{i_k} h_{i_k}(\z^{k+1}_{i_k}),
\end{eqnarray}
where in $\overset{a}=$, we use $e_{k+1,k+1}=n\theta^k$; in $\overset{b}=$, we use $e_{k+1,k}=n(1-\theta^k)\theta^{k-1}+\theta^k-n\theta^k$; and in $\overset{c}=$, we use $e_{k+1,i}=(1-\theta^k)e_{k,i}$ for $i\leq k-1$, and $e_{k,k} =n\theta^{k-1}$.

\subsection{AASVRG}\label{AASVRG-sec}
\begin{lemma}\label{bound}
	Define $f(\x)=\frac{1}{n}\sum_{i=1}^n f_i(\x)$, if $f_i(\x)$'s, $i=1,2,\cdots,n$ have Lipschitz continuous gradients,  for any $\uu$ and $\tx$, defining
	\begin{eqnarray}
	\tna f(\uu) = \nabla f_{k}(\uu) -\nabla f_k(\tx)+ \frac{1}{n}\sum_{i=1}^n \nabla f_i(\tx),
	\end{eqnarray}
	we have
	\begin{eqnarray}\label{lemsec}
	\E\left(\left\|\tna f(\uu) - \nabla  f(\uu)\right\|^2\right) \leq 2L(f(\tx)-f(\uu)+\langle \nabla f(\uu), \uu-\tx \rangle,
	\end{eqnarray}
	where  the expectation is taken on the random number of $k$ under the condition that $\uu$ and $\tx$ are known.
\end{lemma}
The lemma is directly taken from \citep{kat} and \citep{SVRG}. For completeness, we provide a proof of Lemma~\ref{bound} in the end of the section.
We define 
\begin{eqnarray}\label{3e}
\y^s_k = (\theta^s_1)\z^s_k +\theta_2 \tx^s +a^s \x^s_k,
\end{eqnarray}
Through the step 6 in Algorithm 5 in the paper, we have
\begin{eqnarray}\label{4e}
\x^s_{k+1} = \y^s_k + \theta^s_1 (\z^s_{k+1}-\z^s_{k}).
\end{eqnarray}
\bfseries Outline of the Proof: \mdseries \\
Step 1:  and set $b^s(l,k) = (a^{s})^{k-l+1}$. Through the update rule, we have that
\begin{eqnarray}\label{step31}
\y^s_{k}- \w^s_{j(k)} = \sum_{i=j(k)+1}^{k}(\x^s_{i}-\y^s_{i-1})+\sum_{i=j(k)+1}^{k}\left(\sum_{l=i}^{k} b^s(i,l)\right)(\x^s_i-\y^s_{i-1}). 
\end{eqnarray}
Step 2: By analyzing the function value,  we have
\begin{eqnarray}\label{acc1}
&&\E_k f(\x^s_{k+1})\\
&\leq&  f(\y^s_k) - \gamma(1-\frac{\gamma L}{2}- \frac{C_3}{2}-\frac{C_4}{2} )\E\left\|\frac{\x^s_{k+1} - \y^s_{k}}{\gamma}\right\|^2 -\frac{\gamma L^2}{2C_4}\left\|  \w^s_{j(k)} - \y^s_k \right \|^2\notag\\
&&-\E_k\langle \bxi^s_{k+1},  \x^s_{k+1}-\y^s_k\rangle+ \frac{\gamma L }{C_3} \left( f(\tx^s)-f(\w^s_{j(k)})+\langle\nabla f(\w^s_{j(k)}), \w^s_{j(k)} -\tx^s  \rangle  \right).\notag
\end{eqnarray} 
Step 3: By analyzing the $\left\|\z^s_{k+1}-\x^*\right\|^2$,  we have
\begin{eqnarray}\label{3im2}
&&\frac{1}{2\gamma}\E_k\left\|\theta^s_1\z^{s}_{k+1}-\theta^s_1\x^*\right\|^2\\
&\leq&\frac{1}{2\gamma}\left\|\theta^s_1\z^s_k -\theta^s_1\x^*  \right\|^2 + \frac{\gamma}{2}\E_k\left\| \frac{\y^s_{k}-\x^s_{k+1}}{\gamma}\right\|^2\notag\\
&&-\E_k\left\langle  \bxi^s_{k+1},     \y^s_{k}-a^s\x^s_{k}-\theta_{2}\tx^s-\theta^s_1\x^* \right\rangle\notag\\
&&+(1-\theta_{2}-\theta^s_1)f(\x^s_k)+\theta^s_1f(\x^*) -f(\y^s_k) + \theta_2f(\w_{j(k)}^s)\notag\\
&&+L\left\|\w^s_{j(k)}-\y^s_k\right\|^2+\langle  \nabla f(\w^s_{j(k)}), \theta_2(\tx^s -\w^s_{j(k)}) \rangle\notag,
\end{eqnarray}
Step 4: By adding Eq.~\eqref{acc1} and Eq.~\eqref{3im2}, and simplifying, we have
\begin{eqnarray}\label{3im7}
&&\E_k F(\x^s_{k+1})+\frac{1+\frac{\mu \gamma}{\theta^s_1}}{2\gamma}\left\|\theta^s_1\z^{s}_{k+1}-\theta^s_1\x^*\right\|^2\\
&\leq& a^kF(\x^s_k)+\theta^s_1F(\x^*)+ \theta_2 F(\tx^s)+\frac{1}{2\gamma}\left\|\theta^s_1\z^{s}_{k}-\theta^s_1\x^*\right\|^2\notag\\
&&- \gamma(\frac{1}{2}-\frac{3\gamma L}{2}-\frac{C_3}{2} )\E\left\|\frac{\x^s_{k+1} - \y^s_{k}}{\gamma}\right\|^2 -\left(\frac{\gamma L^2}{2C_3}+L\right)\left\|  \w^s_{j(k)} - \y^s_k \right \|^2.\notag
\end{eqnarray}
Step 5: we choose proper step size and obtain Theorem 3 in the paper.\\
\bfseries  Proof of step 1: \mdseries \\
Proof: Through Eq.~\eqref{3e}, we have that
\begin{eqnarray}\label{311}
\theta_1^s \z^s_k + \theta_2 \tx^s = \y^s_k -(1-a^s)\x^s_k, \quad k\geq 0 
\end{eqnarray}
and
\begin{eqnarray}\label{312}
\theta_1^s \z^s_{k+1} + \theta_2 \tx^s = \x^s_{k+1} -(1-a^s)\x^s_{k}, \quad k\geq 0 
\end{eqnarray}
Thus we obtain
\begin{eqnarray}\label{1111}
\y^s_k = \x^s_k +a^s (\x^s_k -\x^s_{k-1}), \quad k\geq 1.
\end{eqnarray}
Eq.~\eqref{1111} is the same with Eq.~\eqref{extra},  thus by the same proof, we can obtain that:
\begin{eqnarray}\label{1232}
&&\y^s_k\\
&=& \x^s_{j(k)}+\sum_{i=j(k)+1}^{k}(\x^s_{i}-\y^s_{i-1})+\sum_{i=j(k)+1}^{k}\left(\sum_{l=i}^{k} b^s(i,l)\right)(\x^s_i-\y^s_{i-1})\notag\\
&&+\left(\sum_{i=j(k)}^{k}b^s(j(k),i)\right)(\x^s_{j(k)}-\x^s_{j(k)-1})\notag,
\end{eqnarray}
Comparing Eq.~\eqref{1232} with the definition of $\w^s_{j(k)}$, we obtain the results.\\
\bfseries  Proof of step 2: \mdseries \\
Through the optimal solution of $\z^s_{k+1}$ in Step 4 of Algorithm 5 in the paper, there exists $\bxi^{s}_{k+1}\in \partial h(\z^s_{k+1}) $, satisfying
\begin{eqnarray}
\theta^s_1(\z^s_{k+1}-\z^s_{k}) +\gamma \tna^s_k  +\gamma \bxi^s_{k+1} = \mathbf{0}.
\end{eqnarray}
Through Eq.~\eqref{311} and Eq.~\eqref{312}, eliminating $\x^s_k$, we have
$$ \theta^s_1(\z^s_{k+1}-\z^s_{k}) + \x^s_{k+1}-\y^s_{k}=\mathbf{0}.$$
So we have
\begin{eqnarray}\label{optimal}
\x^s_{k+1}-\y^s_{k} +\gamma \tna^s_k  +\gamma \bxi^s_{k+1} = \mathbf{0}.
\end{eqnarray}

For $f(\cdot)$ has Lipschitz continuous gradients, we have
\begin{eqnarray}\label{all1}
f(\x^s_{k+1}) &\leq& f(\y^s_k) + \langle \nabla f(\y^s_k), \x^s_{k+1}-\y^s_k\rangle + \frac{L}{2}\left\|\x^s_{k+1}-\y^s_k\right\|^2 \\\notag
&=& f(\y^s_k) -  \gamma \langle \nabla f(\y^s_k),  \tna^s_k +\bxi^{s}_{k+1}  \rangle + \frac{L}{2}\left\|\x^s_{k+1}-\y^s_k \right\|^2\notag\\
&\overset{a}=& f(\y^s_k) -  \gamma \langle \tna^s_k + \bxi^{s}_{k+1}, \tna^s_k +\bxi^{s}_{k+1}  \rangle\notag\\
&& + \frac{L}{2}\left\|\x^s_{k+1}-\y^s_k \right\|^2 + \gamma \langle \tna^s_k + \bxi^{s}_{k+1}-\nabla f(\y^s_k), \tna^s_k+\bxi^{s}_{k+1}   \rangle \notag\\
&\overset{b}=& f(\y^s_k)-\gamma( 1-\frac{\gamma L}{2})\E\left\|\frac{\x^s_{k+1}-\y^s_k}{\gamma} \right\|^2\notag\\
&& - \langle \tna^s_k -\nabla f(\y^s_k), \x^s_{k+1}-\y^s_k   \rangle -\langle \mathbf{\bxi}^s_{k+1}, \x^s_{k+1}-\y^s_k\rangle\notag\\
&\overset{c}=& f(\y^s_k)-\gamma( 1-\frac{\gamma L}{2})\left\|\frac{\x^s_{k+1}-\y^s_k}{\gamma} \right\|^2- \langle \tna^s_k-\nabla f(\w^s_{j(k)}) ,\x^s_{k+1}-\y^s_k   \rangle\notag\\
&& -\langle \nabla f(\w^s_{j(k)}) - \nabla f(\y^s_k),\x^s_{k+1}-\y^s_k   \rangle -\langle \bxi^s_{k+1},  \x^s_{k+1}-\y^s_k\rangle,\notag
\end{eqnarray}
where in equality $\overset{a}=$, we add and subtract the term $\langle \tna^s_k+\bxi^s_k, \tna^s_k+\bxi^s_{k+1} \rangle$;  equality $\overset{b}=$ uses the equality Eq.~\eqref{optimal}; equality $\overset{c}=$, we  add and subtract $\langle \nabla f (\w^s_{j(k)}),\x^s_{k+1}-\y^s_{k} \rangle$.

For the third last term of Eq.~\eqref{all1}, we have
\begin{eqnarray}\label{all2}
&&\E_k  \langle \tna^s_k-\nabla f(\w^s_{j(k)}) ,\y^s_{k}-\x^s_{k+1}   \rangle\\
&\overset{a}\leq&  \frac{\gamma}{2C_3}\E_k\left\|  \tna^s_k-\nabla f(\w^s_{j(k)}) \right\|^2 +\frac{\gamma C_3}{2}\E_k\left\|\frac{\x^s_{k+1}-\y^s_{k}}{\gamma}\right\|^2\notag\\
&\overset{b}\leq& \frac{\gamma L }{C_3} \left( f(\tx^s)-f(\w^s_{j(k)})+\langle\nabla f(\w^s_{j(k)}), \w^s_{j(k)} -\tx^s  \rangle  \right) +\frac{\gamma C_3}{2}\E_k\left\|\frac{\x^s_{k+1}-\y^s_{k}}{\gamma}\right\|^2,\notag
\end{eqnarray}
where we use $\E_k$ to denote that  expectation is taken  on the random number of $i^s_{k}$ (step $k$ and epoch  $s$) under the condition that $\y^s_k$ and $\w^s_{j(k)}$ are known; in  $\overset{a}\leq$, we use the Cauchy-Schwarz inequality;  $\overset{b}\leq$ uses Eq.~\eqref{lemsec}.

For the second last term of Eq.~\eqref{all1}, we have
\begin{eqnarray}\label{all3}
&&\langle  \nabla f(\w^s_{j(k)}) - \nabla f(\y^s_k) ,\y^s_k -\x^s_{k+1}  \rangle\notag\\
&\overset{a}\leq& \frac{\gamma}{2C_4}\E_k\left\|  \nabla f(\w^s_{j(k)}) - \nabla f(\y^s_k)  \right\|^2 + \frac{\gamma C_4}{2}\left\|\frac{\x^s_{k+1} - \y^s_k}{\gamma} \right\|^2\notag\\
&\overset{b}\leq& \frac{\gamma L^2}{2C_4}\left\|  \w^s_{j(k)} - \y^s_k \right \|^2 + \frac{\gamma C_4}{2}\left\|\frac{\x^s_{k+1} - \y^s_k}{\gamma}\right \|^2,
\end{eqnarray}
where in inequality $\overset{a}\leq$, we use Cauchy-Schwarz inequality; in  inequality $\overset{b}\leq$, we use the fact that $f(\cdot)$ has Lipschitz continuous gradients.

Taking expectation for Eq.~\eqref{all1}  and  Eq.~\eqref{all3} on the random number $i^s_{k}$, and adding Eq.~\eqref{all2}, we obtain the results.\\
\bfseries  Proof of step 3: \mdseries \\
We have that 
\begin{eqnarray}\label{all5}
&&\left\|\theta_1^s\z^s_{k+1}-\theta_1^s\x^*\right\|^2\notag\\
&=&\left\|\x^s_{k+1}-a^s\x^s_{k}-\theta_{2}\tx^s-\theta^s_1\x^*\right\|^2\\
&=&\left\|\y^s_{k}-a^s\x^s_{k}-\theta_{2}\tx^s-\theta^s_1\x^* - (\y^s_{k}-\x^s_{k+1})\right\|^2\notag\\
&=&\left\| \y^s_{k}-a^s\x^s_{k}-\theta_{2}\tx^s-\theta^s_1\x^*  \right\|^2 + \left\| \y^s_{k}-\x^s_{k+1}\right\|^2-2\gamma\langle  \bxi^s_{k+1} + \tna_k^s,     \y^s_{k}-a^s\x^s_{k}-\theta_{2}\tx^s-\theta^s_1\x^* \rangle\notag\\
&=&\left\| \theta^s_1\z^s_{k}-\theta^s_1\x^* \right \|^2 +\left \| \y^s_{k}-\x^s_{k+1}\right\|^2-2\gamma\left\langle  \bxi^s_{k+1} + \tna_k^s,     \y^s_{k}-a^s\x^s_{k}-\theta_{2}\tx^s-\theta^s_1\x^* \right\rangle\notag.
\end{eqnarray}
For the last term of Eq.~\eqref{all5}, we have
\begin{eqnarray}\label{all10}
&&\E_k\left\langle  \tna^s_k ,  a^s\x^s_{k}+\theta_{2}\tx^s+\theta^s_1\x^*-\y^s_k \right\rangle\\
&=&\E_k\left\langle  \tna^s_k,  a^s\x^s_{k}+(\theta_{2}-\frac{\gamma L}{C_4})\tx^s+\theta^s_1\x^* -(1-\frac{\gamma L}{C_4})\w^s_{j(k)}\right\rangle\notag\\
&& +\E_k\left\langle \tna^s_k, \w^s_{j(k)}-\y^s_k + \frac{\gamma L}{C_4}(\tx^s -\w^s_{j(k)}) \right\rangle\notag\\
&=&\left\langle  \nabla f(\w_{j(k)}^s),  a^s\x^s_{k}+(\theta_{2}-\frac{\gamma L}{C_4})\tx^s+\theta^s_1\x^* -(1-\frac{\gamma L}{C_4})\w^s_{j(k)}\right\rangle\notag\\
&&+\left\langle \nabla f(\w_{j(k)}^s), \w^s_{j(k)}-\y^s_k + \frac{\gamma L}{C_4}(\tx^s -\w^s_{j(k)}) \right\rangle\notag\\
&\overset{a}\leq& a^sf(\x^s_k)+\theta^s_1f(\x^*) - (1-\theta_2) f(\w^s_{j(k)})\notag\\
&&+\langle\nabla f(\w_{j(k)}^s), \w^s_{j(k)}-\y^s_k\rangle +\left\langle  \nabla f(\w_{j(k)}^s),\theta_2(\tx^s -\w^s_{j(k)}) \right\rangle\notag.
\end{eqnarray}
where in inequality $\overset{a}\leq$, we use the convexity of $f(\cdot)$ and so for any vector $\uu$,
$$\langle\nabla f(\w_{j(k)}^s),  \uu -\w^s_{j(k)}\rangle \leq f(\uu)-f(\w^s_{j(k)}),  $$
and set that $C_4= \frac{\gamma L}{\theta_2}$. For $f(\w^s_{j(k)})$, through the convexity of $f(\cdot)$, we have
\begin{eqnarray}\label{all123}
-f(\w^s_{j(k)})\leq   -f(\y^s_k)+\langle \nabla f(\y^s_{k}), \y^s_{k}- \w^s_{j(k)} \rangle,
\end{eqnarray}
Adding Eq.~\eqref{all123} with Eq.~\eqref{all10}, we have
\begin{eqnarray}\label{all222}
&&\E_k\langle \tna^s_k,  (a^s\x^s_{k}+\theta_{2} \tx^s+\theta^s_1\x^*-\y^s_k \rangle\\
&\leq& (1-\theta_{2}-\theta^s_1)f(\x^s_k)+\theta^s_1f(\x^*) -f(\y^s_k) + \theta_2 f(\w_{j(k)}^s) \notag\\
&&+\langle\nabla f(\w_{j(k)}^s) - \nabla f(\y^s_k), \w^s_{j(k)}-\y^s_k\rangle +\langle  \nabla f(\w_{j(k)}^s), \theta_2 (\tx^s -\w^s_{j(k)}) \rangle\notag\\
&\overset{a}\leq&(1-\theta_{2}-\theta^s_1)f(\x^s_k)+\theta^s_1f(\x^*) -f(\y^s_k) +\theta_2 f(\w_{j(k)}^s)\notag\\
&&+L\left\|\w^s_{j(k)}-\y^s_k\right\|^2+\langle  \nabla f(\w_{j(k)}^s), \theta_2 (\tx^s -\w^s_{j(k)}) \rangle\notag,
\end{eqnarray}
where in $\overset{a}\leq$, we use the $\langle \tna^s_k-\nabla f(\y^s_k), \w^s_{j(k)}-\y^s_k \rangle\leq L\left\|\w^s_{j(k)}-\y^s_k\right\|^2$;  Thus dividing Eq.~\eqref{all5} by $2\gamma$,  and taking expectation on the random number of $i^s_{k}$, we have
\begin{eqnarray}\label{all22}
&&\frac{1}{2\gamma}\E_k\left\|\theta^s_1\z^{s}_{k+1}-\theta^s_1\x^*\right\|^2\\
&\leq&\frac{1}{2\gamma}\left\| \theta^s_1\z^{s}_{k}-\theta^s_1\x^*\right \|^2 + \frac{\gamma}{2}\E_k\left\| \frac{\y^s_{k}-\x^s_{k+1}}{\gamma}\right\|^2\notag\\
&&-\E_k\left\langle  \bxi^s_{k+1} + \tna^k_s,     \y^s_{k}-a^s\x^s_{k}-\theta_{2}\tx^s-\theta^s_1\x^* \right\rangle\notag\\
&\overset{a}\leq&\frac{1}{2\gamma}\left\|\theta^s_1\z^s_k -\theta^s_1\x^*  \right\|^2 + \frac{\gamma}{2}\E_k\left\| \frac{\y^s_{k}-\x^s_{k+1}}{\gamma}\right\|^2\notag\\
&&-\E_k\left\langle  \bxi^s_{k+1},     \y^s_{k}-a^s\x^s_{k}-\theta_{2}\tx^s-\theta^s_1\x^* \right\rangle\notag\\
&&+(1-\theta_{2}-\theta^s_1)f(\x^s_k)+\theta^s_1f(\x^*) -f(\y^s_k) + \theta_2 f(\w_{j(k)}^s)\notag\\
&&+L\left\|\w^s_{j(k)}-\y^s_k\right\|^2+\langle  \nabla f(\w^s_{j(k)}), \theta_2 (\tx^s -\w^s_{j(k)}) \rangle\notag,
\end{eqnarray}
where $\overset{a}\leq$ uses Eq.~\eqref{all222}. This is the result.\\
\bfseries  Proof of step 4: \mdseries \\
Adding Eq.~\eqref{all22} and Eq.~\eqref{acc1},  we obtain the that:
\begin{eqnarray}\label{11231}
&&\E_k f(\x^s_{k+1})+\frac{1}{2\gamma}\left\|\theta^s_1\z^{s}_{k+1}-\theta^s_1\x^*\right\|^2\\
&\leq& a^sf(\x^s_k)+\theta^s_1f(\x^*)  + \theta_2f(\w_{j(k)}^s)+\langle  \nabla f(\w^s_{j(k)}),\theta_2(\tx^s -\w^s_{j(k)}) \rangle\notag\\
&&- \gamma(\frac{1}{2}-\frac{\gamma L}{2}- \frac{C_3}{2}-\frac{C_4}{2} )\E\left\|\frac{\x^s_{k+1} - \y^s_{k}}{\gamma}\right\|^2 -\left(\frac{\gamma L^2}{2C_3}+L\right)\left\|  \w^s_{j(k)} - \y^s_k \right \|^2\notag\\
&&-\E_k\langle \bxi^s_{k+1},  \x^s_{k+1}-\y^s_k\rangle+ \theta_2 \left( f(\tx^s)-f(\w^s_{j(k)})+\langle\nabla f(\w^s_{j(k)}), \w^s_{j(k)} -\tx^s  \rangle  \right)\notag\\
&&-\E_k\left\langle  \bxi^s_{k+1} ,     \y^s_{k}-a^s\x^s_{k}-\theta_{2}\tx^s-\theta^s_1\x^* \right\rangle+\frac{1}{2\gamma}\left\|\theta^s_1\z^{s}_{k}-\theta^s_1\x^*\right\|^2\notag\\
&\overset{a}\leq&a^sf(\x^s_k)+\theta^s_1f(\x^*) + \theta_2 f(\tx^s)+\frac{1}{2\gamma}\left\|\theta^s_1\z^{s}_{k}-\theta^s_1\x^*\right\|^2\notag\\
&&- \gamma(\frac{1}{2}-\frac{3\gamma L}{2}- \frac{C_3}{2})\E\left\|\frac{\x^s_{k+1} - \y^s_{k}}{\gamma}\right\|^2 -\left(\frac{\gamma L^2}{2C_3}+L\right)\left\|  \w^s_{j(k)} - \y^s_k  \right\|^2\notag\\
&&-\E_k\left\langle  \bxi^s_{k+1} ,     \x^s_{k+1}-a^s\x^s_{k}-\theta_{2}\tx^s-\theta^s_1\x^* \right\rangle\notag,
\end{eqnarray}
where in $\overset{a}\leq$, we use $C_4=\frac{\gamma L}{\theta_2}$. For the last term of Eq.~\eqref{11231}, we have
\begin{eqnarray}\label{qq}
&&-\left\langle  \bxi^s_{k+1} ,     \x^s_{k+1}-a^s\x^s_{k}-\theta_{2}\tx^s-\theta^s_1\x^* \right\rangle\\
&= & -\left\langle  \bxi^s_{k+1} ,     \theta^s_1 \z^s_{k+1}-\theta^s_1\x^* \right\rangle\notag\\
&\leq& \theta^s_{1} h(\x^*)-\theta^s_{1} h(\z^s_{k+1})-\frac{\mu}{2}\left\|\z^{k+1}-\x^* \right\|^2\notag\\
&\overset{a}\leq& \theta^s_{1} h(\x^*)- h(\x^s_{k+1}) + \theta_2 h(\tx^s)+ a^s h(\x^s_k) -\frac{\mu}{2}\left\|\z^{k+1}-\x^* \right\|^2,\notag
\end{eqnarray}
where in $\overset{a}\leq$, we use $\x^s_{k+1}=a^s\x^s_{k}+\theta_{2}\tx^s+\theta^s_{1}\x^s_{k+1}$, and the convexity of $h(\cdot)$. Substituting Eq.~\eqref{qq} into Eq.~\eqref{11231}, we obtain the result.\\
\bfseries  Proof of step 5: \mdseries \\
\begin{eqnarray}\label{3onee}
&&\left\|\w^{j(k)} -\y^k\right\|^2\\
&\leq&   \left( \sum_{i=j(k)+1}^{k}\left(1+ \sum_{l=i}^{k}b(i,l)\right)\right) \sum_{i=j(k)+1}^{k}\left(1+ \sum_{l=i}^{k}b(i,l)\right)\left \|\x^i -\y^{i-1}\right\|^2 \notag\\
&\overset{a}\leq& \left( \sum_{i=j(k)+1}^{k}\left(1+  \sum_{l=1}^{k-i+1} \frac{1}{2} \right)\right) \sum_{i=j(k)+1}^{k}\left(1+ \sum_{l=1}^{k-i+1} \frac{1}{2} \right) \left\|\x^i -\y^{i-1}\right\|^2 \notag\\
&\leq& \left( \sum_{ii=1}^{\min(\tau,k)}\left(1+  \sum_{l=1}^{ii}  \frac{1}{2} \right)\right) \sum_{ii=1}^{\tau}\left(1+ \sum_{l=1}^{ii} \frac{1}{2} \right)\left \|\x^{k-ii+1} -\y^{k-ii}\right\|^2 \notag\\
&\leq& 4\tau \sum_{i=1}^{\min(\tau,k)}\left\|\x^{k-i+1} -\y^{k-i}\right\|^2,\notag
\end{eqnarray}
where in $\overset{a}\leq$, we use $b(i,l)\leq \frac{1}{2}$ when $l\geq i$.
We first consider the not-strongly convex case. Using the same technique of Eq.~\eqref{onee}, summing Eq.~\eqref{3onee} with $k = 0$ to $m-1$, we have 
\begin{eqnarray}\label{3im6}
&&\sum_{k=0}^m\left\|\w^{j(k)} -\y^k\right\|^2\\
&\leq&4\tau\sum_{k=0}^{m-1} \sum_{i=1}^{\min(\tau,k-\tau)} \left \|\x^{k-i+1} -\y^{k-i}\right\|^2\notag\\
&\leq& 4\tau^2 \sum_{k=0}^{m-1} \left \|\x^{k+1}-  \y^{k}\right\|^2.
\end{eqnarray}
Taking expectation on Eq.~\eqref{3im7} for the first $k-1$ iteration (all the random numbers coming from epoch s), and summing it with $k=0$ to $m-1$, we have
\begin{eqnarray}\label{ac2}\small
&&\sum_{k = 0}^{m-1}  \E\left( F(\x^s_{k+1})-F(\x^*)\right)+ \frac{1}{2\gamma}\E\left\|\theta^s_1\z^{s}_{m}-\theta^s_1\x^*\right\|^2  \\
&\leq&  \sum_{k = 0}^{m-1}  a^s\E\left(F(\x^s_{k}) -F(\x^*)\right) +m\theta_2 \left( F(\tx^s) - F(\x^*) \right)+\frac{1}{2\gamma}\E\left\|\theta^s_1\z^s_0-\theta^s_1\x^*\right\|^2 \notag\\ &&- \gamma(\frac{1}{2}-\frac{3\gamma L}{2}-\frac{C_3}{2} )\sum_{k=0}^m\E\left\|\frac{\x^s_{k+1} - \y^s_{k}}{\gamma}\right\|^2   +\left(\frac{\gamma L^2}{2C_3}+L\right)\sum_{k=0}^m\left\|  \w^s_{j(k)} - \y^s_k \right \|^2\notag\\
&\overset{a}\leq&  \sum_{k = 0}^{m-1}  a^s\E\left(F(\x^s_{k}) -F(\x^*)\right) +m\theta_2 \left( F(\tx^s) - F(\x^*) \right) +\frac{1}{2\gamma}\E\left\|\theta^s_1\z^s_0-\theta^s_1\x^*\right\|^2\notag\\ &&- \gamma(\frac{1}{2}-\frac{3\gamma L}{2}-\frac{C_3}{2}-4\tau^2\left( \gamma L+\frac{\gamma^2 L^2}{C_3}\right) )\sum_{k=0}^m\E\left\|\frac{\x^s_{k+1} - \y^s_{k}}{\gamma}\right\|^2,
\end{eqnarray}
where $\overset{a}\leq$ uses Eq.~\eqref{3im6}.  By setting $C_3 = 2\gamma L$, we obtain the that
$$1-\frac{3\gamma L}{2}-\frac{C_3}{2}-4\tau^2\left( \gamma L+\frac{\gamma^2 L^2}{2C_3}\right)\leq 0$$

The rest proof is similar to \citep{kat}.
Diving $(\theta^s_1)^2$ on both side of Eq.~\eqref{ac2} and arranging terms, we have
\begin{eqnarray}\label{ac5}
&&\frac{1}{(\theta^s_1)^2} \E\left( F(\x^s_{m})-F(\x^*)\right)+\frac{\theta_{2}+\theta^s_1}{(\theta^s_1)^2}\sum_{k=1}^{m-1}\E\left(F(\x^s_{k})-F(\x^*) \right)\\
&\leq&\frac{1-\theta^s_1-\theta_2}{(\theta^s_1)^2}  (F(\x^s_0)-F(\x^*))+ \frac{m\theta_{2}}{(\theta^s_1)^2}\left(F(\tx^s)-F(\x^*)\right)\notag\\
&&+\frac{1}{2\gamma}\left\|\z^{s}_0-\x^*\right\|^2-\frac{1}{2\gamma}\E\left\|\z^{s}_m-\x^*\right\|^2 \notag.
\end{eqnarray}
When $s>0$, through the definition of $\tx^s$, we have
\begin{eqnarray}\label{ac4}
F(\tx^s) = F(\frac{1}{m}\sum_{k=0}^{m-1} \x^{s-1}_{k} )\leq \frac{1}{m}\sum_{k=0}^{m-1} F(\x^{s-1}_k) = \frac{1}{m}F(\x^s_0) +\frac{1}{m}\sum_{k=1}^{m-1} F(\x^{s-1}_k).
\end{eqnarray}
Through the definition of $\z^s_0$, we have
\begin{eqnarray}\label{ac3}
\z^{s}_0 = \z^{s-1}_m.
\end{eqnarray}
Substituting Eq.~\eqref{ac4} and Eq.~\eqref{ac3} into Eq.~\eqref{ac5}, we have
\begin{eqnarray}
&&\frac{1}{(\theta^s_1)^2} \E\left( F(\x^s_{m})-F(\x^*)\right)+\frac{\theta_{2}+\theta^s_1}{(\theta^s_1)^2}\sum_{k=1}^{m-1}\E\left(F(\x^s_{k})-F(\x^*) \right)\\
&\leq&\frac{1-\theta^s_1}{(\theta^s_1)^2}  (F(\x^{s-1}_m)-F(\x^*))+\frac{\theta_{2}}{(\theta^s_1)^2}\sum_{k=1}^{m-1}\left(F(\x^{s-1}_{k})-F(\x^*) \right) \notag\\
&&+\frac{1}{2\gamma}\left\|\z^{s}_0-\x^*\right\|^2-\frac{1}{2\gamma}\E\left\|\z^{s+1}_0-\x^*\right\|^2, \quad s>0. \notag
\end{eqnarray}
Since $\theta^s_1 = \frac{2}{s+4}\leq \frac{1}{2}$, we have
\begin{eqnarray}
\frac{1}{(\theta^s_1)^2} \geq \frac{1-\theta^{s+1}_1}{(\theta^{s+1}_1)^2},\quad s\geq 0,
\end{eqnarray}
and
\begin{eqnarray}
\frac{\theta_{2}+\theta^s_1}{(\theta^s_1)^2}\geq \frac{\theta_{2}}{(\theta^{s+1}_1)^2},\quad s\geq 0.
\end{eqnarray}
So
\begin{eqnarray}\label{ac6}
&&\frac{1}{(\theta^s_1)^2} \E\left( F(\x^s_{m})-F(\x^*)\right)+\frac{\theta_{2}+\theta^s_1}{(\theta^s_1)^2}\sum_{k=1}^{m-1}\E\left(F(\x^s_{k})-F(\x^*) \right)\\
&\leq&\frac{1}{(\theta^{s-1}_1)^2}  (F(\x^{s-1}_m)-F(\x^*))+\frac{\theta_{2}+\theta^{s-1}_1}{(\theta^{s-1}_1)^2}\sum_{k=1}^{m-1}\left(F(\x^{s-1}_{k})-F(\x^*) \right) \notag\\
&&+\frac{1}{2\gamma}\left\|\z^{s}_0-\x^*\right\|^2-\frac{1}{2\gamma}\E\left\|\z^{s+1}_0-\x^*\right\|^2, \quad s>0. \notag
\end{eqnarray}
When $s=0$, through Eq.~\eqref{ac5}, use $\tx^0 = \x^0_0$, we have
\begin{eqnarray}\label{ac7}
&&\frac{1}{(\theta^0_1)^2} \E\left( F(\x^0_{m})-F(\x^*)\right)+\frac{\theta_{2}+\theta^0_1}{(\theta^0_1)^2}\sum_{k=1}^{m-1}\E\left(F(\x^s_{0})-F(\x^*) \right)\\
&\leq&\frac{1-\theta_{1,0}+(m-1)\theta_{2}}{(\theta^0_1)^2}  (F(\x^0_0)-F(\x^*))\notag\\
&&+\frac{1}{2\gamma}\left\|\z^{s}_0-\x^*\right\|^2-\frac{1}{2\gamma}\E\left\|\z^{s+1}_0-\x^*\right\|^2, \notag.
\end{eqnarray}
Taking expectation for Eq.~\eqref{ac6} with $s$ from $1$ to $S$ (random numbers coming from the $0$ to $s-1$ epochs) and summing the result with $S$ from $1$ to $S-1$, and adding Eq.~\eqref{ac7}, we obtain
\begin{eqnarray}\label{ac8}
&&\frac{1}{(\theta^S_1)^2} \E\left( F(\x^S_{m})-F(\x^*)\right)+\frac{\theta_{2}+\theta^S_1}{(\theta^S_1)^2}\sum_{k=1}^{m-1}\E\left(F(\x^S_{k})-F(\x^*) \right)\\
&\leq&\frac{1-\theta^0_1+(m-1)\theta_{2}}{(\theta^0_1)^2}  (F(\x^0_0)-F(\x^*))+\frac{1}{2\gamma}\E\left\|\z^0_0-\x^*\right\|^2-\frac{1}{2\gamma}\E\left\|\z^{S+1}_0-\x^*\right\|^2 \notag\\
&\overset{a}\leq&2m  (F(\x^0_0)-F(\x^*))+\frac{1}{2\gamma}\left\|\z^0_0-\x^*\right\|^2,\notag
\end{eqnarray}
where in $\overset{a}\leq$, we use $\theta^0_1= \theta_2 = \frac{1}{2}$.

Now we consider the strongly convex case. Through the definition  of $\gamma$, we have
\begin{eqnarray}
8\gamma L\tau^2 \leq 5\gamma L + \frac{95}{8}\gamma L\tau^2 \leq 1, 
\end{eqnarray}
and $\theta^s_1= \frac{1}{\tau}\sqrt{\frac{n\mu}{L}}$. Set $\theta_3 = \frac{\mu\gamma}{\theta^s_1}+1 \leq \frac{1}{8\tau}\sqrt{\frac{\mu}{Ln}}+1$.
Multiply  Eq.~\eqref{3onee} with $\theta_3^k$, and summing the results with $k$ from $0$ to $m-1$, we have
\begin{eqnarray}
&&\sum_{k=0}^K \theta_3^k\left\|\w^{j(k)} -\y^k\right\|^2\\
&\leq&\sum_{k=0}^K 4\tau\theta_3^k \sum_{i=1}^{\min(\tau,k)}\left \|\x^{k-i+1} -\y^{k-i}\right\|^2\notag\\
&\leq&\sum_{k=0}^K 4\tau\sum_{i=1}^{\min(\tau,k)}\theta_3^{i} \theta_3^{k-i}\left\|\x^{k-i+1} -\y^{k-i}\right\|^2\notag\\
&\leq&\sum_{k=0}^K 4\tau\theta_3^{\tau} \sum_{i=1}^{\min(\tau,k)} \theta_3^{k-i}\left\|\x^{k-i+1} -\y^{k-i}\right\|^2\notag\\
&\overset{a}\leq&\sum_{k=0}^K 4\tau(1+\frac{3}{16}) \sum_{i=1}^{\min(\tau,k)} \theta_3^{k-i}\left\|\x^{k-i+1} -\y^{k-i}\right\|^2\notag\\
&\overset{b}\leq&  \frac{19\tau^2}{4}\sum_{k=0}^{K-1}  \theta_3^{k}\left\|\x^{k+1} -\y^{k}\right\|^2\leq  \frac{19\tau^2}{4}\sum_{k=0}^{K}  \theta_3^{k}\left\|\x^{k+1} -\y^{k}\right\|^2\notag,
\end{eqnarray}
where   $\overset{b}\leq$ is because that for each $\theta_3^i\left\|\x^{i+1} -\y^{i}\right\|^2$ ($1\leq k \leq K$) there are most $\tau$ terms, like Eq.~\eqref{im6}; in $\overset{a}\leq$, we use the fact that for the function $g(x)=(1+x)^a \leq 1+\frac{3}{2}ax$, when $a\leq 1$, and $x\leq\frac{1}{a}$.
To proof it, we can use Taylor expansion at point $x=0$ to obtain
\begin{eqnarray}\label{end2}
(1+x)^a  = 1 + ax +\frac{a(a-1)}{2} \xi^2 \leq   1 + ax +\frac{a(a-1)}{2}  \frac{1}{a} x \leq 1+\frac{3}{2}ax,
\end{eqnarray}
where $\xi\in [0, x]$, and
$$\theta_3^\tau \leq (1+\frac{1}{8\tau}\sqrt{\frac{\mu}{nL}})^\tau \overset{a}{\leq}  (1+\frac{1}{8\tau}\frac{\tau}{n})^\tau \leq 1+ \frac{3\tau}{16n}\overset{b}{\leq}  \frac{19}{16},$$
where in $\overset{a}\leq$, we use the assumption that $n\mu \leq  \tau^2L$; and in  $\overset{b}\leq$, we use $\tau\leq n$.
Taking expectation on Eq.~\eqref{3im7} for the first $k-1$ iterations, and then multiply it with $\theta_3^k$, and summing the results with $k$ from $0$ to $m$, we have
\begin{eqnarray}\label{end1}
&&\sum_{k=1}^{m} \theta_3^{k-1} \left(F(\x^s_{k})-F(\x^*) \right) - (a^s)\sum_{k=0}^{m-1} \theta_3^k \left(F(\x^s_{k})-F(\x^*) \right)  -  \theta_2\sum_{k=0}^{m-1} \theta_3^k \left(F(\tx^s)-F(\x^*) \right) \notag\\
&& + \frac{1}{2\gamma}\left\| \theta_1^s \z^s_0 - \theta_1^s\x^*\right\|^2 - \frac{\theta_3^{m}}{2\gamma}\left\| \theta_1^s\z^s_m- \theta_1^s\x^*\right\|^2 \z^s_0 \notag\\
&\leq& -\gamma(1-\frac{3\gamma L}{2}-\frac{C_3}{2} )\sum_{k=0}^{m-1}\theta_3^k\E\left\|\frac{\x^s_{k+1} - \y^s_{k}}{\gamma}\right\|^2 -\sum_{k=0}^{m-1}\left(\frac{\gamma L^2}{2C_3}+L\right)\left\|  \w^s_{j(k)} - \y^s_k  \right\|^2\notag\\
&\leq&- \gamma\left(\frac{1}{2}-\frac{3\gamma L}{2}-\frac{C_3}{2} - \left(\frac{\gamma^2 L^2}{2C_3}+L\gamma\right)\frac{19\tau^2}{4}\right) \sum_{k=0}^{m-1}\theta_3^k\E\left\|\frac{\x^s_{k+1} - \y^s_{k}}{\gamma}\right\|^2.
\end{eqnarray}
Set $C_3 = 2\gamma L$, we have
we have 
$$\frac{1}{2}-\frac{5\gamma L}{2} - \frac{95\tau^2\gamma^2 L^2}{16}\geq 0.$$
The rest proof is similar to \citep{kat}.  By arranging the terms of Eq.~\eqref{end1}, we have
\begin{eqnarray}
&& (\theta_1^s+\theta_{2}  -(1-1/\theta_3))\sum_{i=1}^{m} \theta_3^k \left(F(\x^s_{k})-F(\x^*) \right)  + \theta_3 ^ma^s \left(F(\x^s_{m})-F(\x^*)\right)\notag\\
&\leq&  \theta_2\sum_{k=0}^{m-1} \theta_3^k \left(F(\tx^s)-F(\x^*) \right) +a^s \left(F(\x^s_{0})-F(\x^*)\right)\notag\\
&&+ \frac{1}{2\gamma}\left\| \theta_1^s \z^s_0 - \theta_1^s\x^*\right\|^2 - \frac{\theta_3^{m}}{2\gamma}\left\| \theta_1^s\z^s_m- \theta_1^s\x^*\right\|^2.
\end{eqnarray}
Through the definition of $\tx^{s+1} = (\sum_{j=0}^{m-1}\theta_3^j)^{-1}\sum_{j=0}^{m-1}\x^s_j\theta_3^j$,  we have
\begin{eqnarray}
&& (\theta_1^s+\theta_{2}  -(1-1/\theta_3))\theta_3 \sum_{k=0}^{m-1} \theta_3^k \left(F(\tx^{s+1})-F(\x^*) \right)  + \theta_3 ^ma^s \left(F(\x^s_{m})-F(\x^*)\right)\notag\\
&\leq&  \theta_2\sum_{k=0}^{m-1} \theta_3^k \left(F(\tx^s)-F(\x^*) \right) +a^s \left(F(\x^s_{0})-F(\x^*)\right)\notag\\
&&+ \frac{1}{2\gamma}\left\| \theta_1^s \z^s_0 - \theta_1^s\x^*\right\|^2 - \frac{\theta_3^{m}}{2\gamma}\left\| \theta_1^s\z^s_m- \theta_1^s\x^*\right\|^2.
\end{eqnarray}
Since
\begin{eqnarray}\label{end3}
&&\theta_2 (\theta_3^{m-1}-1)+ (1-1/\theta_3) \notag\\
&\leq&\frac{1}{2}\left( \left(1+ \frac{1}{8\tau}\sqrt\frac{\mu}{nL}\right)^{m-1}-1 \right) + \frac{\frac{1}{8\tau}\sqrt\frac{L}{n\mu}}{\theta_3}\notag\\
&\overset{a}\leq& \frac{3}{4} \frac{1}{8\tau}\sqrt\frac{n\mu}{L} +\frac{\frac{1}{8\tau}\sqrt\frac{mu}{nL}}{\theta_3}\notag\\
&\leq& \frac{1}{\tau}\sqrt{\frac{n\mu}{L}} = \theta^s_1,
\end{eqnarray}
where in $\overset{a}\leq$, we  use $\mu n \leq L\tau^2$, and Eq.~\eqref{end2}. Eq.~\eqref{end3} indicates that $(\theta_1^s+\theta_{2}  -(1-1/\theta_3))\leq \theta_2\theta_3^{m-1}$, so
\begin{eqnarray}
&&\theta_2\theta_3^m \sum_{k=0}^{m-1} \theta_3^k \left(F(\tx^{s+1})-F(\x^*) \right)  + \theta_3 ^ma^s \left(F(\x^s_{m})-F(\x^*)\right)\notag\\
&\leq&  \theta_2\sum_{k=0}^{m-1} \theta_3^k \left(F(\tx^s)-F(\x^*) \right) +a^s \left(F(\x^s_{0})-F(\x^*)\right)\notag\\
&&+ \frac{\gamma}{2}\left\| \theta_1^s \z^s_0 - \theta_1^s\x^*\right\|^2 - \frac{\theta_3^{m}\gamma}{2}\left\| \theta_1^s\z^s_m- \theta_1^s\x^*\right\|^2.
\end{eqnarray}
By telescope the above inequality from $s=0,\cdots, S$, we have that
\begin{eqnarray}
&&\theta_2 \sum_{k=0}^{m-1} \theta_3^k \left(F(\tx^{s+1})-F(\x^*) \right)  +  (1-\theta_1^s-\theta_2) \left(F(\x^s_{m})-F(\x^*)\right)\notag\\
&\leq& (\theta_3)^{-Sm}\left(  \left(\theta_2\sum_{k=0}^{m-1} \theta_3^k +  (1-\theta_1^s-\theta_2)\right) \left(F(\x^0_0)-F(\x^*) \right)  +\frac{(\theta^s_{1})^2\gamma}{2}\left\|\z^0_0-\x^* \right\|^2  \right).\notag\\
\end{eqnarray}
Since $\theta_3^k \geq 1$, and $\theta_2 = \frac{1}{2}$, and so $\sum_{k=0}^{m-1}\theta^k_3 \geq n$, $\theta_1^s\leq\frac{1}{2}$, we have
\begin{eqnarray}
\left(  F(\tx^{S+1})-F(\x^*)   \right) \leq   (\theta_3)^{-Sn} \left(    (1+\frac{1}{n}) \left(F(\x^0_0)-F(\x^*) \right)  +\frac{\gamma}{4n}\left\|\z^0_0-\x^* \right\|^2 \right).
\end{eqnarray}
This ends proof.

\begin{lemma}\label{lips}
	Suppose $f(\x)$ has Lipschitz continuous gradients, for any $\x$ and $\y$, we have
	\begin{eqnarray}
	\left\| \nabla f(\x) -\nabla f(\y) \right \|^2 \leq 2 L \left( f(\x) -f(\y)+\langle\nabla f(\y), \y-\x  \rangle  \right).
	\end{eqnarray}
\end{lemma}
Lemma~\ref{lips} is Theorem 2.1.5 of the textbook of Nesterov~\citep{nesterov2013introductory}.

\bfseries  Proof of Lemma~\ref{bound}: \mdseries
\begin{eqnarray}
&& \E\left(\left\|\tna f(\uu) - \nabla f(\uu)\right\|^2\right)\notag\\
&=& \E \left(\left\|\left(\nabla f_k(\uu) - \nabla f_k(\tx) + \nabla f(\tx) - \nabla f(\uu)\right)\right\|^2 \right)\notag\\
&=&   \E \left(\left\|\left(\nabla f_k(\uu) -\nabla f_k(\tx) - \left( \nabla f(\uu) - \nabla f(\tx)\right)\right\|^2 \right)\right)\notag\notag\\
&\overset{a}\leq& \E \left(\left\|\nabla f_k(\uu) -\nabla f_k(\tx)\right\|^2\right)\label{one111}.
\end{eqnarray}
where in  inequality $\overset{a}\leq$, we use that $$\E\left(\nabla f_k(\uu) -\nabla f_k(\tx) \right) = \nabla f(\uu) - \nabla f(\tx), $$ and $\E\left(\left\|\x-\E(\x)\right\|^2\right)=\E \left\| \x\right\|^2-\left\|\E(\x)\right\|^2\leq \E\left \| \x\right \|^2$;  Then by directly applying Lemma~\ref{lips} to Eq.~\eqref{one111}, we obtain Eq.~\eqref{lemsec}.

\subsection{ASVRG}\label{ASVRG-sec}
The algorithm of ASVRG is shown in Algorithm~\ref{ASVRG1}. We analyze ASVRG in the wild scheme. For the wild scheme, there is no locks in implementation.  So different coordinates of $\x$ read by any child node may at different iteration steps. So $\x_{j(k)}$ is not a real state of $\x^k$.  However, the update
on a single coordinate can be considered to be atomic~\citep{HOGWILD,ASGD1,asvrg1}.  Through~\citep{asvrg1}, we represented $\x^s_{j(k)}$ as follows:
\begin{equation}\label{im}
\x^s_{j(k)} = \x^s_{k} -\sum_{l=1}^{k-1} \gamma I^s_{k(l)}\left( \vv^s_{j(l)}\right),
\end{equation}
where $I^s_{k(l)}$ is an  $\mathcal{R}^d \to\mathcal{R}^d$ function, indicating whether the elements of $\vv^s_{j(l)}$ have been returned from the local memory and written into $\x$ at the ``read'' step in the $k$-th iteration  and $d$ is the dimension of the variable $\x$. Suppose $\vv^s_{j(l)}(p)$ is the $p$-th element of $\vv^s_{j(l)}$ with $p$ ranging from $1$ to $d$. We have
\begin{eqnarray}\label{im1}
I^s_{k(l)}\!\!\left( \vv^s_{j(l)}\right)\!\!(p) \!=\!\!\left\{ \begin{array}{ll}
0, &\!\!\!\!\! \text{if}\ \vv^s_{j(l)}(p) \ \text{has been returned}, \\
\vv^s_{j(l)}(p),&\!\!\!\text{otherwise}.
\end{array}
\right.
\end{eqnarray}
We can find that atom scheme also satisfies  Eq.~\eqref{im}. Now we begin our proof. The proof can be also consider as an extension of \citep{asvrg1}. Lemma~\ref{aabound} is similar to Lemma 1 in \citep{asvrg1}. The variant that we adopt is to consider the term $\E(\left\|\nabla f(\x^s_k)\right\|^2)+4L\E\left(f(\x^s_k)+f(\x^s)-2f(\x^*)\right)$.  For simply, we assume $h(\x) \equiv \mathbf{0}$.  We first prove Lemma~\ref{aabound}.

\begin{algorithm}[tb]
	\caption{ASVRG}
	\label{ASVRG1}
	\begin{algorithmic}
		\STATE $\mathbf{Input}$ $\x^0_0$. Set epoch length $m=2n$, $\tx^0 = \x^0_0$, step size $\gamma$.\\
		\STATE1 $\mathbf{for}$ $s =0$ $\mathbf{to}$ $S-1$ $\mathbf{do}$\\
		\STATE2  $\quad \mathbf{for}$ $k=0$ $\mathbf{to}$ $m-1$ $\mathbf{do}$\\
		\STATE3  $\quad\ \ $ Randomly sample $i_{k,s}$ from ${1,2,\cdots,n}$,
		\STATE5 $\quad\quad \tna f(\x^s_{j(k)}) = \nabla f_{i_{k,s}} (\x^s_{j(k)}) - \nabla f_{i_{k,s}} (\tx^s)+\frac{1}{n}\sum_{i=1}^n \nabla f_{i}(\tx^s)   ,$\\
		\STATE6 $\quad\quad \x_{k+1}^{s} = \x_{k}^{s}-\gamma \tna f(\x^s_{j(k)}),$\\
		\STATE7 $\quad$  $\mathbf{end}$ $\mathbf{for}$ $k$.\\
		\STATE8  $\quad \x^{s+1}_0 = \frac{1}{m}\sum_{k=1}^{m} \x^s_k,$\\
		\STATE9  $\mathbf{end}$ $\mathbf{for}$ $s$.\\
	\end{algorithmic}
\end{algorithm}
\begin{lemma}\label{aabound}
	Suppose $f(\x)$ has Lipschitz continuous gradients, for ASVRG, if the step size $\gamma$ satisfies
	\begin{equation}\label{condition}
	\gamma \leq \min  \frac{1}{L}\left\{\frac{\rho_1-1}{10 \rho_1\sqrt{\rho_2}}, \frac{\rho_2-1}{10 \rho_1^{\frac{1}{2}}\rho_2^{\frac{3}{2}}\frac{\rho_1^{\frac{\tau}{2}}-1}{\sqrt{\rho_1}-1}} \right\},
	\end{equation}
	for some $\rho_1 > 1$ and $\rho_2>1$, then for any $s\geq 0$ and $k \geq 0 $, we have
	\begin{eqnarray}\label{condition1}
	&& \E\left(\left\|\nabla f(\x^s_{k})\right\|^2 \right) + 4L \E\left(f(\x^s_k) +f(\tx^s) -2f(\x^*)\right)\notag\\
	& \leq& \rho_1 \left[ \E\left(\left\|\nabla f(\x^s_{k+1})\right\|^2 \right)\! + 4L \E\left(f(\x^s_{k+1}) +f(\tx^s) -2f(\x^*)\right)\right],
	\end{eqnarray}
	and
	\begin{eqnarray}
	&&  \E\left(\left\|\nabla f(\x^s_{j(k)})\right\|^2 \right) + 4L \E\left(f(\x^s_{j(k)}) +f(\tx^s) -2f(\x^*)\right)\notag\\
	&\leq& \rho_2  \E\left(\left\|\nabla f(\x^s_{k})\right\|^2\right ) + 4L \E\left(f(\x^s_k) +f(\tx^s) -2f(\x^*)\right).
	\end{eqnarray}
\end{lemma}
\bfseries  Proof of Lemma~\ref{aabound}: \mdseries \\
We analyze $\left\|\nabla f(\x^s_{k}) \right\|^2$ and $f(\x^s_{k})-f(\x^*)+f(\tx^s)-f(\x^*)$, respectively. For  $\left\|\nabla f(\x^s_{k}) \right\|^2$, we have
\begin{eqnarray}\label{one1}
&& \E\left(\left\|\nabla f(\x^s_k)\right\|^2-(\left\|\nabla f(\x^s_{k+1})\right\|^2\right)\notag\\
&\overset{a}\leq& 2\E\left( \left\|\nabla f(\x^s_k)\right\|\left\|\nabla f(\x^s_k)- \nabla f(\x^s_{k+1})\right\|   \right)\notag\\
&\overset{b}\leq& 2L\E\left(\left \|\nabla f(\x^s_k)\right\|\left\|\x^s_k- \x^s_{k+1}\right\|   \right)\notag\\
&\overset{c}\leq& L\gamma \E \left(\frac{1}{C_3}  \left\|\nabla f(\x^s_k)\right\|^2 +C_3\left\| \tna f(\x^s_{j(k)})\right \|^2  \right), \quad (C_3\geq 0).
\end{eqnarray}
where in inequality $\overset{a}\leq$, we use the equality that $\|\mathbf{a} \|^2-\|\mathbf{b}\|^2\leq 2\|\mathbf{a} \|\| \mathbf{a}-\mathbf{b} \|$;  inequality $\overset{b}\leq$ uses the fact that $f(\cdot)$ is L-smooth;  inequality $\overset{c}\leq$ uses the Cauchy-Schwarz inequality. For $f(\x^s_{k})-f(\x^*)+f(\tx^s)-f(\x^*)$, we have
\begin{eqnarray}\label{one2}
&&\E\left(f(\x^s_{k})-f(\x^*)+f(\tx^s)-f(\x^*)\right)-\E\left(f(\x^s_{k+1})-f(\x^*)+f(\tx^s)-f(\x^*)\right)\notag\\
&=& \E \left( (f(\x^s_{k} )-  f(\x^s_{k+1})\right)\notag\\
&\overset{a}\leq& \E \left(  \langle  \nabla f(\x^s_k), \x^s_k -\x^s_{k+1}  \rangle \right)\notag\\
&\overset{b}\leq&\gamma\frac{1}{C_3}  \E \left \|\nabla f(\x^s_k)\right\|^2 +  \gamma C_3\E\left\| \tna f(\x^s_{j(k)})\right \| ^2,
\end{eqnarray}
where in inequality $\overset{a}\leq$ we use the convexity of $f(\cdot)$; inequality $\overset{b}\leq$ uses the Cauchy-Schwarz inequality.

Then similarly, we have
\begin{eqnarray}\label{dis3}
&& \E \left(\left\|\nabla f(\x^s_{j(k+1)})\right\|^2 - \left\|\nabla f(\x^s_{k+1}) \right\|^2 \right) \\
&\leq& 2 \E \left( \left \|\nabla f(\x^s_{j(k+1)})\right \| \left\|\nabla f(\x^s_{k+1}) - \nabla f(\x^s_{j(k+1)})\right\| \right) \notag \\
&\leq& \frac{L\gamma}{C_4} \E \left( \left\|\nabla f(\x^s_{j(k+1)})\right\|^2\right) +\frac{LC_4}{\gamma}\E \left( \left\| \x^s_{k+1} - \x^s_{j(k+1)}\right\|^2    \right)\quad \left(C_4 > 0\right)\notag,
\end{eqnarray}
and
\begin{eqnarray}\label{dis4}
&&\E\left(f(\x^s_{j(k)})-f(\x^*)+f(\tx^s)-f(\x^*)\right)-\E\left(f(\x^s_{k+1})-f(\x^*)+f(\tx^s)-f(\x^*)\right)\notag\\
&=& \E \left( (f(\x^s_{j(k)} )-  f(\x^s_{k+1})\right)\notag\\
&\overset{a}\leq& \E \left(  \langle  \nabla f(\x^s_{j(k)}), \x^s_{j(k)} -\x^s_{k+1}  \rangle \right)\notag\\
&\overset{b}\leq&\gamma\frac{1}{C_4}  \E \left \|\nabla f(\x^s_{j(k+1)})\right\|^2 +  \gamma C_4 \left\| \x^s_{k+1} - \x^s_{j(k+1)}\right\|^2.
\end{eqnarray}

For convenience, we set
\begin{eqnarray}\label{bk}
B_k = \E \left\|\nabla f(\x^s_{k-1})\right\|^2 + 4L\E \left(f(\x^s_{k-1})-f(\x^*)+f(\tx^s)-f(\x^*)\right),
\end{eqnarray}
which has omitted the superscript $s$. For the fact that $\E\left\|\x \right\|^2 = \left(\E \left\| \x\right\|\right)^2+D(\x) $, we have
\begin{eqnarray}
\E\left\| \tna f(\x^s_{j(k)}) \right\| ^2  = \left\| \nabla  f(\x^s_{j(k)})\right\|^2 +  \E\left\|\tna  f(\x^s_{j(k)})- \nabla  f(\x^s_{j(k)})\right\| ^2.
\end{eqnarray}
Then from  Lemma~\ref{bound}, we have
\begin{eqnarray}
\E \left( \left\|\tna f(\x^s_k)\right\|^2   \right)\leq B_k.
\end{eqnarray}
Multiplying Eq.~\eqref{one2} by $4L$ and adding Eq.~\eqref{one1}, we obtain
\begin{eqnarray}\label{temp1}
&&B_{k}-B_{k+1}\notag\\
&\leq& 5L\gamma \E \left(\frac{1}{C_3} \left \|\nabla f(\x^s_k)\right\|^2 +C_3\left\| \tna f(\x^s_{j(k)})\right \|^2  \right)\notag\\
&\leq& 5L\gamma \frac{1}{C_3} B_k + 5L\gamma C_3 B_{j(k)}.
\end{eqnarray}
Now we use induction to prove $B_{k-1} \leq \rho_1 B_{k}$, and $B_{j(k)}\leq \rho_2 B_k$.
Suppose $k=1$,  $B_{j(k)}=B_0$, we have
\begin{eqnarray}\label{one33}
B_0 - B_1 \leq 10 L\gamma B_0,
\end{eqnarray}
where we set $C_3=1$. Simplifying Eq.~\eqref{one33}, we have
\begin{eqnarray}
B_0 \leq \frac{1}{1-10L\gamma}B_1.
\end{eqnarray}
Recalling the $\gamma$, we have
\begin{eqnarray}
L\gamma \leq \frac{\rho_1-1}{10\rho_1\sqrt{\rho_2}}\leq \frac{\rho_1-1}{10\rho_1} = \frac{1}{10}\left( 1-\frac{1}{\rho_1}\right).
\end{eqnarray}
so
\begin{eqnarray}
B_0  \leq \frac{1}{1 - 10 L \gamma} B_1 \leq \rho_1 B_1.
\end{eqnarray}
On the other hand, multiplying  Eq.~\eqref{dis4} by $4L$ and then adding Eq.~\eqref{dis3}, we have
\begin{eqnarray}\label{big}
B_{j(k+1)} - B_{k+1} &\leq& \frac{5L\gamma}{C_4} B_{j(k+1)} + \frac{5LC_4}{\gamma}\E \left(\left \| \x^s_{k+1} - \x^s_{j(k+1)}\right\|^2    \right).
\end{eqnarray}
When $k=1$,
\begin{eqnarray}\label{smal1}
&&\E \left(\left \| \x^s_{1} - \x^s_{j(1)}\right\|^2    \right) \notag\\
&\overset{a}\leq&  \gamma ^2 \E \left( \left\|I^s_{0(0)} (\tna f(\x^s_0))\right\|^2  \right)  \notag\\
&\leq& \gamma ^2 \E \left( \left\|\tna f(\x^s_0)\right\|^2  \right)\notag\\
&\leq& \gamma ^2 \rho_1 B_1,
\end{eqnarray}
where in inequality $\overset{a}\leq$ we use the definition of $\x^s_{j(k)}$ in Eq.~\eqref{im} and Eq.~\eqref{im1}. Substituting Eq.\eqref{smal1} into Eq.\eqref{big}, we have
\begin{eqnarray}
B_{j(1)} - B_{1} \leq \frac{5L\gamma}{C_4} B_{j(1)} + 5\gamma L C_4 \rho_1 B_1 .
\end{eqnarray}
Setting $C_4$ to be $\frac{1}{\sqrt{\rho_1}}$,
\begin{eqnarray}
B_{j(1)} - B_{1}  \leq 5\sqrt{\rho_1} L \gamma B_{j(1)}+5\sqrt{\rho_1} L \gamma B_{1}.
\end{eqnarray}
Then
\begin{eqnarray}
B_{j(1)}  \leq \frac{1+5\sqrt{\rho_1} L\gamma }{1 - 5\sqrt{\rho_1} L \gamma} B_1.
\end{eqnarray}
Recalling the assumption on $\gamma$, we have
\begin{eqnarray}\label{rho21}
L\gamma \leq \frac{\rho_2-1}{10\sqrt { \rho_1} \rho_2^{\frac{3}{2}}\frac{\rho_1^{\frac{\tau}{2}}-1}{\sqrt{\rho_1}-1}} \leq  \frac{\rho_2-1}{10\sqrt { \rho_1}\rho_2}.
\end{eqnarray}
So we have
\begin{eqnarray}
10\sqrt{\rho_1} L \gamma \leq 1-\frac{1}{\rho_2} < 1.
\end{eqnarray}
Then
\begin{eqnarray}\label{rho22}
&& B_{j(1)}\notag\\\
& \leq&  \frac{1+5\sqrt{\rho_1} L\gamma }{1 - 5\sqrt{\rho_1} L \gamma} B_1 \notag\\\
&\leq&   \frac{1 }{1 - 10\sqrt{\rho_1} L \gamma} B_1 \notag\\\
&\leq& \rho_2 B_1,
\end{eqnarray}
where we use the fact that $\frac{1+x}{1-x}\leq \frac{1}{1-2x}$ when $2x<1$ in the second inequality.

When $B_{k}$ satisfies $B_{k-1} \leq \rho_1 B_{k}$, and $B_{j(k)}\leq \rho_2 B_k$, we consider $B_{k+1}$.  From Eq.~\eqref{temp1},
\begin{eqnarray}
B_k -B_{k+1} &\leq&    \frac{5L\gamma}{C_3} B_k + 5LC_3\gamma \rho_2 B_{k}.
\end{eqnarray}
Setting $C_3 = \frac{1}{\sqrt{\rho_2}}$, we have
\begin{eqnarray}
B_k - B_{k+1} \leq  10\sqrt{\rho_2} L \gamma  B_k.
\end{eqnarray}
Then
\begin{eqnarray}
B_k  \leq \frac{1}{1- 10\sqrt{\rho_2}L \gamma } B_{k+1}.
\end{eqnarray}
From the assumption on $\gamma$, we have $B_k \leq \rho_1 B_{k+1}$.  The same as  Eq.~($28$) and Eq.~($29$) in \citep{asvrg1}, we have
\begin{eqnarray}
&&\E \left(\left\|\x^s_{k+1} -\x^s_{j(k+1)}\right\|^2\right) \notag\\
&=&  \gamma ^2 \E \left( \left\|\sum_{l=k-\tau+1}^{k} I_{k(l)}\left( \tna f(\x^s_{j(l)})\right)\right\|^2  \right)  \notag\\
&\leq&  \gamma ^2 \E \left( \sum_{p=1}^{d}\left(\sum_{l=k-\tau+1}^{k}  |\tna f(\x^s_{j(l)})(p)|\right)^2  \right). \notag\\
&\leq& \gamma^2\E \left( \sum_{p=1}^{d}\left( \sum_{i=0}^{\tau-1}\sum_{z=0}^{\tau-1}|\tna f(\x^s_{j(k-i)})(p)|\times |\tna f(\x^s_{j(k-z)})(p)|\right) \right),
\end{eqnarray}
where $\tna f(\x^s_{k})(p)$ is the $p$-th coordinate of vector $\tna f(\x^s_{k})$, the first inequality uses the  inequality that $(a_1+a_2+\cdots+a_\tau)^2 \leq (|a_1|+|a_2|+\cdots+|a_\tau|)^2$ on each dimension.
For any $i=0,1,\dots, \tau-1$ and $z=0,1,\dots, \tau-1$, we have
\begin{eqnarray}
&&\E \left( \sum_{p=1}^{d}\left(2 |\tna f(\x^s_{j(k-i)})(p)|\times|\tna f(\x^s_{j(k-z)})(p)|\right) \right)\notag\\
&\overset{a}\leq& \E \left(\sum_{p=1}^{d} \left( \rho_1^{(z-i)/2} |\tna f(\x^s_{j(k-i)})(p)|^2 +  \rho_1^{(i-z)/2} |\tna f(\x^s_{j(k-z)})(p)|^2 \right)     \right) \\\notag
&\leq& \E \left( \rho_1^{(z-i)/2} \left\|\tna f(\x^s_{j(k-i)})\right\|^2 +  \rho_1^{(i-z)/2} \left\|\tna f(\x^s_{j(k-z)})\right\|^2      \right)\\\notag
&\leq&  \rho_1^{(z-i)/2} B_{j(k-i)} +  \rho_1^{(i-z)/2} B_{j(k-z)}      \\\notag
&\leq& \rho_2 \rho_1^{(z-i)/2} \rho_1^i B_k  + \rho_2\rho_1^{(i-z)/2} \rho_1^z B_k \\\notag
&\leq& 2 \rho_2\rho_1^{(i+z)/2}B_k,
\end{eqnarray}
where in inequality $\overset{a}\leq$, we use Cauchy-Schwarz.
So
\begin{eqnarray}\label{span}
&&\E \left(\left\|\x^s_{k+1} -\x^s_{j(k+1)}\right\|^2\right) \notag\\
&\leq& \gamma^2 \rho_2  \sum_{i=0}^{\tau-1}\sum_{z=0}^{\tau-1} \rho_1^{(i+z)/2}B_k \notag\\
&\leq& \gamma ^2 \rho_2 \left(\sum_{i=0}^{\tau-1} \rho_1^{i/2}  \right)^2B_k\notag\\
&\leq& \gamma ^2 \rho_2  \frac{(\rho_1^{\tau/2}-1)^2}{(\sqrt{\rho_1} -1)^2}B_k.
\end{eqnarray}
Substituting Eq.~\eqref{span} into Eq.~\eqref{big}, we have
\begin{eqnarray}
&& B_{j(k+1)} - B_{k+1}  \notag\\
&\leq& \frac{5L\gamma}{C_4} B_{j(k+1)} + 5LC_4\gamma\rho_2 \frac{(\rho_1^{\tau/2}-1)^2}{(\sqrt{\rho_1} -1)^2} B_{k}  \notag\\
&\leq& \frac{5L\gamma}{C_4} B_{j(k+1)} + 5LC_4\gamma\rho_2\rho_1 \frac{(\rho_1^{\tau/2}-1)^2}{(\sqrt{\rho_1} -1)^2} B_{k+1}.
\end{eqnarray}
Setting $C_4 = \frac{1}{\sqrt{\rho_1\rho_2} \frac{\rho_1^{\tau/2}-1}{\sqrt{\rho_1} -1} }$, we have
\begin{equation}
B_{j(k+1)} - B_{k+1}\leq 5L\gamma \sqrt{\rho_1\rho_2} \frac{\rho_1^{\tau/2}-1}{\sqrt{\rho_1} -1} (B_{j(k+1)} + B_{k+1} ).
\end{equation}
Considering the assumption on $\gamma$, like Eq.~\eqref{rho21}, we have
\begin{eqnarray}
10L\gamma \sqrt{\rho_1\rho_2} \frac{\rho_1^{\tau/2}-1}{\sqrt{\rho_1} -1} \leq 1-\frac{1}{\rho_2}< 1,
\end{eqnarray}
then like Eq.~\eqref{rho22}, we have
\begin{eqnarray}
&&B_{j(k+1)} \notag\\
&\leq& \frac{1+ 5L\gamma \sqrt{\rho_1\rho_2} \frac{\rho_1^{\tau/2}-1}{\sqrt{\rho_1} -1}  }{1-5L\gamma \sqrt{\rho_1\rho_2} \frac{\rho_1^{\tau/2}-1}{\sqrt{\rho_1} -1} }B_{k+1} \notag\\
&\leq& \frac{1  }{1-10L\gamma \sqrt{\rho_1\rho_2} \frac{\rho_1^{\tau/2}-1}{\sqrt{\rho_1} -1} }B_{k+1}\notag\\
&\leq& \rho_2 B_{k+1}.
\end{eqnarray}
So Lemma~\ref{aabound} is proved.\\\\
\bfseries  Proof of the convergence results: \mdseries \\
\begin{theorem}\label{theorem1}
	Suppose the step size $\gamma$ in ASVRG satisfies $\gamma \leq \{\frac{(\sqrt{5}-\sqrt{2})\sqrt{2}}{20\cdot 5^\frac{3}{4}\sqrt{e}(\sqrt{e}-1)\tau L}, \frac{1}{12\sqrt{5}e(e-1)\tau^2 L}  \}$, we have
	\begin{eqnarray}
	\E \left(F(\tx^{s+1}) -F(\x^*)\right) \leq  \left(F(\x^0_0)-F(\x^*)\right)+\frac{9}{16\gamma m}\left\|\x_0^0-\x^* \right\|^2.
	\end{eqnarray}
\end{theorem}
Recalling Eq.~\eqref{span}, we have
\begin{eqnarray}\label{boundx}
&&\E\left \|\x^s_{k}-\x^s_{j(k)} \right \|^2\notag\\
&\leq&  \gamma ^2 \rho_2  \frac{(\rho_1^{\tau/2}-1)^2}{(\sqrt{\rho_1} -1)^2}B_{k-1}\notag\\
&\leq&  \gamma ^2 \rho_2 \rho_1 \frac{(\rho_1^{\tau/2}-1)^2}{(\sqrt{\rho_1} -1)^2}B_{k}.
\end{eqnarray}

We first consider $f(\x)$. For $f(\x)$ has Lipschitz continuous gradients, we have
\begin{eqnarray}\label{svrg1}
\E_k f(\x^s_{k+1}) &\leq& f(\x^s_k) + \E_k\langle \nabla f(\x^s_k), \x^s_{k+1}-\x^s_k\rangle + \frac{L}{2}\E_k\left\|\x^s_{k+1}-\x^s_k\right\|^2 \\\notag
&\overset{a}=& f(\x^s_k) -  \gamma \E_k\langle \nabla f(\x^s_k),  \tna f(\x^s_{j(k)})  \rangle + \E_k\frac{\gamma^2L}{2}\left\|\tna f(\x_{j(k)}) \right\|^2\notag\\
&\overset{b}=&  f(\x^s_k)-  \gamma \langle  \nabla f(\x^s_k), \nabla f(\x^s_k) \rangle +\E_k\frac{\gamma^2 L}{2}\left\| \tna f(\x^s_{j(k)})\right \|^2+ \gamma \E_k\langle   \nabla f(\x^s_k), \nabla f(\x^s_k) - \tna f(\x^s_{j(k)})  \rangle  \notag\\
&\overset{c}=&  f(\x^s_k)-  \gamma \langle  \nabla f(\x^s_k), \nabla f(\x^s_k) \rangle +\frac{\gamma^2 L}{2}\E_{i_k}\left\| \tna f(\x^s_{j(k)})\right \|^2+ \gamma \E_k\langle   \nabla f(\x^s_k), \nabla f(\x^s_k) - \tna f(\x^s_{j(k)})   \rangle, \notag
\end{eqnarray}
where the expectation $\E_k$ is taken over the random numbers of  $i_{k,s}$ under the condition that $\x^s_k$ is known; in equality $\overset{a}=$, we use  $\x^s_{k+1} =\x^s_k -\gamma \tna f(\x^s_{j(k)}) $; in equality $\overset{b}=$, we replace $\langle  \nabla f(\x^s_k), \tna f(\x^s_k) \rangle$ with $\langle  \nabla f(\x^s_k), \nabla f(\x^s_k) -\nabla f(\x^s_k) +\tna f(\x^s_k) \rangle$; in $\overset{c}=$, we use $\E_k\langle \nabla f(\x^s_k), \tna f(\x^s_k)   \rangle = \langle \nabla f(\x^s_k), \nabla f(\x^s_k)   \rangle$. Taking expectation on all the random numbers on Eq.~\eqref{svrg1}, we have
\begin{eqnarray}\label{f1}
&&\E f(\x^s_{k+1})\\
&\leq&  \E f(\x^s_k)-\gamma \E\langle  \nabla f(\x^s_k), \nabla f(\x^s_k) \rangle +\frac{\gamma^2 L}{2}\E\left\| \tna f(\x^s_{j(k)}) \right\|^2+\gamma\E\langle   \nabla f(\x^s_k), \nabla f(\x^s_k) - \tna f(\x^s_{j(k)})   \rangle \notag\\
&\overset{a}\leq&\E f(\x^s_k)-\gamma \E\left\| \nabla f(\x^s_k) \right \|^2 +\frac{\gamma^2 L}{2}\E\left\| \tna f(\x^s_{j(k)}) \right\|^2 +\gamma \frac{ C_5}{2}\E\left\|\nabla f(\x^s_k)\right \|^2 + \frac{\gamma L^2}{2C_5} \E\left\|\x^s_{j(k)} - \x^s_k\right\|^2\notag \\
&\overset{b}\leq&\E f(\x^s_k)-\gamma\E \left\| \nabla f(\x^s_k) \right \|^2 +\frac{\gamma^2 L\rho_2}{2}B_k +\gamma \frac{ C_5}{2}\E\left\|\nabla f(\x^s_k)\right \|^2 + \frac{\gamma ^3L^2 \rho_2 \rho_1}{2C_5} \frac{(\rho_1^{\tau/2}-1)^2}{(\sqrt{\rho_1} -1)^2}B_k\notag \\
&\overset{c}\leq&\E f(\x^s_k) - \gamma\left(1-\frac{\rho_2 \gamma L}{2}-\frac{C_5}{2}- \frac{\gamma ^2L^2 \rho_2 \rho_1}{2C_5} \frac{(\rho_1^{\tau/2}-1)^2}{(\sqrt{\rho_1} -1)^2}  \right) \E\left\|\nabla f(\x^s_k) \right \|^2\notag\\
&&+4L^2 \gamma^2\left( \frac{\rho_2}{2} + \frac{\gamma L\rho_2 \rho_1}{2C_5} \frac{(\rho_1^{\tau/2}-1)^2}{(\sqrt{\rho_1} -1)^2} \right) \E\left(f(\x^s_k)-f(\x^*)+f(\tx^s)-f(\x^*) \right),\notag\\
&\overset{d}\leq&\E f(\x^s_k) - \gamma\left(1-\frac{\rho_2 \gamma L}{2}-\frac{\gamma L}{2}- \frac{\gamma L \rho_2 \rho_1}{2} \frac{(\rho_1^{\tau/2}-1)^2}{(\sqrt{\rho_1} -1)^2}  \right) \E\left\|\nabla f(\x^s_k)  \right\|^2\notag\\
&&+4L^2 \gamma^2\left( \frac{\rho_2}{2} + \frac{\rho_2 \rho_1}{2} \frac{(\rho_1^{\tau/2}-1)^2}{(\sqrt{\rho_1} -1)^2} \right) \E\left(f(\x^s_k)-f(\x^*)+f(\tx^s)-f(\x^*) \right),\notag
\end{eqnarray}
where in $\overset{a}\leq$, we use  Cauchy-Schwarz inequality and the smoothness of $f(\cdot)$, i.e. $\left\|\nabla f(\x^s_k)-\tna f(\x^s_{j(k)}\right  \|^2\leq L^2\left \|\x^s_{j(k)} - \x^s_k \right \|^2$, in $\overset{b}\leq$, we substitute $\E \left\|\tna f(\x^s_{j(k)})\right\|^2 \leq \rho_2 B_k$ and Eq.~\eqref{boundx}; in  $\overset{c}\leq$, we use the definition of $B_k$ in Eq.~\eqref{bk}; in  $\overset{d}\leq$, we set $C_5 = \gamma L$.

On the other hand, for $\left\|\x^s_{k+1}-\x^*\right\|^2$, we have
\begin{eqnarray}\label{one44}
&&\left\|\x^s_{k+1}-\x^*\right\|^2\\
&=& \left\| \x^s_{k} -\x^* -\gamma \tna f(\x^s_{j(k)}) \right\|^2 \notag\\
&=& \left\|\x^s_k -\x^*\right\|^2 -2\gamma \langle \tna f(\x^s_{j(k)}), \x^s_k-\x^* \rangle +\gamma^2 \left\| \tna f(\x^s_{j(k)})\right \|^2\notag\\
&\overset{a}=&\left\|\x^s_k -\x^*\right\|^2 -2\gamma \langle \tna f(\x^s_{j(k)}), \x^s_{j(k)}-\x^*  \rangle +\gamma^2\left\|\tna f(\x^s_{j(k)}) \right \|^2 -2\gamma \langle \tna f(\x^s_{j(k)}),  \x^s_k - \x^s_{j(k)}  \rangle\notag\\
&\overset{b}\leq& \left\|\x^s_k -\x^*\right\|^2 +2\gamma\left( f(\x^*)-f(\x^s_{j(k)})\right)+\gamma^2\left\|\tna f(\x^s_{j(k)})\right \|^2 - 2\gamma \langle \tna f(\x^s_{j(k)}),  \x^s_k - \x^s_{j(k)}  \rangle\notag\\
&\overset{c}\leq& \left\|\x^s_k -\x^*\right\|^2+2\gamma\left( f(\x^*)-f(\x^s_{j(k)})\right)-2\gamma \langle \nabla f(\x^s_k), \x^s_k - \x^s_{j(k)} \rangle\notag\\
&&+\gamma^2\left\|\tna f(\x^s_{j(k)}) \right \|^2-2\gamma \langle \tna f(\x^s_{j(k)})- \nabla f(\x^s_k),  \x^s_k - \x^s_{j(k)}  \rangle\notag\\
&\overset{d}\leq& \left\|\x^s_k -\x^*\right\|^2+2\gamma\left( f(\x^*)-f(\x^s_{k})\right)\notag\\
&&+\gamma^2\left\|\tna f(\x^s_{j(k)})\right\|^2-2\gamma \langle \tna f(\x^s_{j(k)})- \nabla f(\x^s_k),  \x^s_k - \x^s_{j(k)}  \rangle\notag
\end{eqnarray}
where in equality $\overset{a}=$, we replace $\langle \tna f(\x^s_{j(k)}), \x^s_{k}  \rangle$ with $\langle \tna f(\x^s_{j(k)}), \x^s_{j(k)} - \x^s_{j(k)}+ \x^s_{k}  \rangle$;  in inequality $\overset{b}=$, we use the convexity of $f(\cdot)$:
$$f(\x^*)- f(\x^s_{j(k)})\geq  \langle \nabla f(\x^s_{j(k)}), \x^* - \x^s_{j(k)} \rangle;$$
in inequality $\overset{c}\leq$, we add and subtract the term $\langle \nabla f(\x^s_k), \x^s_k - \x^s_{j(k)} \rangle$; in inequality $\overset{d}\leq$, we use the fact that
$$f(\x^s_{j(k)}) - f(\x^s_k)\geq  -\langle \nabla f(\x^s_{k}), \x^s_{k}- \x^s_{j(k)} \rangle.$$

Taking expectation only on the random number $i_{k,s}$ on Eq.~\eqref{one44}, and use the fact that
\begin{eqnarray}
\E_k\langle \nabla f(\x_{j(k)}) -\tna f(\x_{j(k)}),  \x_k - \x_{j(k)}   \rangle=0,
\end{eqnarray}
we have
\begin{eqnarray}\label{one55}
&&\E_k\left\|\x^s_{k+1}-\x^*\right\|^2\\
&\leq& \E_k\left\|\x^s_k -\x^*\right\|^2+2\gamma\E_k\left( f(\x^*)-f(\x^s_{k})\right)\notag\\
&&+\gamma^2\E_k\left\|\tna f(\x^s_{j(k)}) \right \|^2-2\gamma \E_k\langle \nabla f(\x^s_{j(k)})- \nabla f(\x^s_k),  \x^s_k - \x^s_{j(k)}  \rangle,\notag
\end{eqnarray}
Taking expectation on all the random numbers on Eq.~\eqref{one55}, we have
\begin{eqnarray}\label{one66}
&&\E\left\|\x^s_{k+1}-\x^*\right\|^2\\
&\leq&  \E\left\|\x^s_k -\x^*\right\|^2+2\gamma\E\left( f(\x^*)-f(\x^s_{k})\right)\notag\\
&&+\gamma^2\E\left\|\tna f(\x^s_{j(k)})  \right\|^2-2\gamma \E\langle \nabla f(\x^s_{j(k)})- \nabla f(\x^s_k),  \x^s_k - \x^s_{j(k)}  \rangle,\notag\\
&\overset{a}\leq&\E\left\|\x^s_k -\x^*\right\|^2+2\gamma\E\left( f(\x^*)-f(\x^s_{k})\right)\notag\\
&&+\gamma^2\E\left\|\tna f(\x^s_{j(k)})  \right\|^2+ 2\gamma ^3L \rho_2 \rho_1 \frac{(\rho_1^{\tau/2}-1)^2}{(\sqrt{\rho_1} -1)^2}B_{k},\notag
\end{eqnarray}
where in $\overset{a}\leq$, we use Eq.~\eqref{boundx} and the smoothness of $f(\cdot)$.
Diving Eq.~\eqref{one66} by $2\gamma$ on both size,  and using the definition of $B_k$ in Eq.~\eqref{bk}, we have
\begin{eqnarray}\label{one4}
&&\E f(\x^s_k)-f(\x^*) +\frac{1}{2\gamma}\E\left\|\x^s_{k+1}-\x^*\right\|^2\\
&\leq& \frac{1}{2\gamma}\E\left\|\x^s_{k}-\x^*\right\|^2+\frac{\gamma}{2}\E\left\|\tna f(\x^s_{j(k)})\right \|^2+\gamma ^2L \rho_2 \rho_1 \frac{(\rho_1^{\tau/2}-1)^2}{(\sqrt{\rho_1} -1)^2}B_k\notag\\
&\leq&  \frac{1}{2\gamma}\E\left\|\x^s_{k}-\x^*\right\|^2+\gamma\left(\frac{ \rho_2}{2}  +\gamma L \rho_2 \rho_1 \frac{(\rho_1^{\tau/2}-1)^2}{(\sqrt{\rho_1} -1)^2} \right)\E\left\| \nabla f(\x^s_k)\right \|^2\notag\\
&&+4\gamma L\left(\frac{ \rho_2}{2}  + \gamma L \rho_2 \rho_1 \frac{(\rho_1^{\tau/2}-1)^2}{(\sqrt{\rho_1} -1)^2} \right) \E\left( f(\x^s_k)-f(\x^*)+f(\tx^s) -f(\x^*) \right).\notag
\end{eqnarray}
Multiply Eq.~\eqref{one4} by $\frac{9}{8}$ and add it to Eq.~\eqref{f1}, we have
\begin{eqnarray}
&&\E f(\x^s_{k+1})-f(\x^*) + \left(\frac{1}{8}-C_6 \right)\E\left(f(\x^s_{k})-f(\x^*)\right) + \frac{9}{16\gamma}\E\left\|\x^s_{k+1}-\x^*\right\|^2\notag\\
&\leq& C_6\E\left(f(\tx^s)-f(\x^*)\right)+\frac{9}{16\gamma}\E\left\|\x_k^s-\x^* \right \|^2 -\gamma C_7 \E\left\|\nabla f(\x^s_k)\right\|^2.
\end{eqnarray}
where 
\begin{equation}
C_6 =  (\frac{9}{4}\gamma L+ 2\gamma^2L^2)\rho_2  + \frac{13}{2}\gamma^2 L^2 \rho_2 \rho_1 \frac{(\rho_1^{\tau/2}-1)^2}{(\sqrt{\rho_1} -1)^2}, 
\end{equation}
and
\begin{eqnarray}
C_7 &= &1-\frac{\rho_2 \gamma L}{2}-\frac{\gamma L}{2} - \frac{9\rho_2}{16} - \frac{13}{8} \gamma L \rho_2 \rho_1 \frac{(\rho_1^{\tau/2}-1)^2}{(\sqrt{\rho_1} -1)^2}\notag.
\end{eqnarray}
We first verify that$\gamma \leq \{\frac{(\sqrt{5}-\sqrt{2})\sqrt{2}}{20\cdot 5^\frac{3}{4}\sqrt{e}(\sqrt{e}-1)\tau L}, \frac{1}{12\sqrt{5}e(e-1)\tau^2 L}  \}$,  $\rho^\tau_1 = e$ and $\rho_2 = \frac{\sqrt{5}}{2}$ satisfies the condition of Lemma~\ref{aabound}.
\begin{eqnarray}
\frac{\rho_1 -1}{10\rho_1 \rho_2^{\frac{1}{2}}}\geq \frac{e^{\frac{1}{\tau}}-1}{10 e^\frac{1}{\tau}\left(\frac{\sqrt{5}}{2}\right)^\frac{1}{2}}\overset{a}\geq\frac{\sqrt{2}}{10\cdot 5^\frac{1}{4}e \tau}\geq L\gamma,
\end{eqnarray}
where we use the fact that $e^x-1\geq x$ for $x\geq 0$ and $\frac{1}{e^\frac{1}{\tau}}\geq\frac{1}{e}$ in inequality $\overset{a}\geq$. In addition, we have
\begin{eqnarray}
\frac{\rho_2 -1}{10\rho_1^{\frac{1}{2}}\rho_2^{\frac{3}{2}}\frac{\rho_1^{\frac{\tau}{2}}-1}{\sqrt{\rho_1}-1}}\geq \frac{\rho_2-1}{10\cdot \rho_2^\frac{3}{2}\sqrt{e}}\frac{\rho_1^{\frac{1}{2}}-1}{\rho_1^{\frac{\tau}{2}}-1}\overset{a}\geq \frac{\rho_2-1}{20\cdot \rho_2^\frac{3}{2}\sqrt{e}(\sqrt{e}-1)\tau} = \frac{(\sqrt{5}-2)\sqrt{2}}{20\cdot 5^\frac{3}{4}\sqrt{e}(\sqrt{e}-1)\tau} \geq L\gamma,\notag
\end{eqnarray}
where in $\overset{a}\geq$, we use $\rho_1^{\frac{1}{2}}-1 \geq \frac{1}{2\tau}$. 
Since 
\begin{eqnarray}
&&\rho_2 \rho_1 \frac{(\rho_1^{\tau/2}-1)^2}{(\sqrt{\rho_1} -1)^2}\notag\\
&\leq& \frac{\sqrt{5}}{2} e (e-1) \frac{1}{(e^\frac{1}{2\tau}-1)^2}\notag\\
&\leq& 2\sqrt{5}e(e-1)\tau^2,
\end{eqnarray}
from the assumption of $\gamma\leq \frac{1}{12\sqrt{5}e(e-1)\tau^2 L}$, we can also get
\begin{eqnarray}
&&L\gamma\rho_2 \rho_1 \frac{(\rho_1^{\tau/2}-1)^2}{(\sqrt{\rho_1} -1)^2}\notag\\
&\leq& \frac{1}{12\sqrt{5}e(e-1)\tau^2}\cdot 2\sqrt{5} e (e-1)\tau^2\notag\\
&=&\frac{1}{6}.
\end{eqnarray}
From the assumption of $\gamma\leq \frac{(\sqrt{5}-\sqrt{2})\sqrt{2}}{20\cdot 5^\frac{3}{4}\sqrt{e}(\sqrt{e}-1)\tau L}$, we have $L\gamma\leq \frac{1}{100}$, so
\begin{eqnarray}
C_6\leq \frac{\sqrt{5}}{2}\left(\frac{1}{100}\cdot \frac{9}{4} + 2\cdot \frac{1}{100} \cdot \frac{1}{100} \right) + \frac{13}{2} \cdot \frac{1}{100} \cdot \frac{1}{6} \leq\frac{1}{8},
\end{eqnarray}
and
\begin{eqnarray}
C_7\geq 1-\frac{\sqrt{5}}{400}-\frac{1}{200} -\frac{9\sqrt{5}}{32}-\frac{13}{8} \cdot \frac{1}{6} \geq 0.
\end{eqnarray}
We obtain
\begin{eqnarray}
&&\E f(\x^s_{k+1})-f(\x^*)+ \frac{9}{16\gamma}\E\left\|\x^s_{k+1}-\x^*\right\|^2\notag\\
&\leq& C_6\E \left(f(\tx^s)-f(\x^*)\right) + \left(C_6 - \frac{1}{8} \right)\E\left(f(\x^s_{k})-f(\x^*)\right) + \frac{9}{16\gamma}\E\left\|\x^s_{k}-\x^*\right\|^2\notag\\
&\leq& C_6\E\left(f(\tx^s)-f(\x^*)\right)+\frac{9}{16\gamma}\E\left\|\x_k^s-\x^* \right \|^2 \notag\\
&\leq& \E\left(f(\tx^s)-f(\x^*)\right)+\frac{9}{16\gamma}\E\left\|\x_k^s-\x^* \right \|^2. \notag\\
\end{eqnarray}
Summing $k$ from $0$ to $m-1$ and using the fact that
\begin{eqnarray}
f(\tx^{s+1}) \leq \sum_{k=1}^{m} f(\x^{s}_k),
\end{eqnarray}
we have
\begin{eqnarray}\label{final}
&& m\E \left(f(\tx^{s+1}) -f(\x^*)\right)+ \frac{9}{16\gamma}\E\left\|\x^s_{m}-\x^*\right\|^2\notag\\
&\leq&  m\E\left(f(\tx^s)-f(\x^*)\right)+\frac{9}{16\gamma}\E\left\|\x_0^s-\x^*  \right\|^2.
\end{eqnarray}
Summing Eq.~\eqref{final} with $s$ from $0$ to $s$, and using $\x^s_0=\x^{s-1}_m$, we have the results that
\begin{eqnarray}
&&m\E \left(f(\tx^{s+1}) -f(\x^*)\right)+ \frac{9}{16\gamma}\E\left\|\x^s_{m}-\x^*\right\|^2\notag\\
&\leq&  m\left(f(\x^0_0)-f(\x^*)\right)+\frac{9}{16\gamma}\left\|\x_0^0-\x^* \right \|^2.
\end{eqnarray}
So
\begin{eqnarray}
&&\E \left(f(\tx^{s+1}) -f(\x^*)\right)\notag\\
&\leq&  \left(f(\x^0_0)-f(\x^*)\right)+\frac{9}{16\gamma m}\left\|\x_0^0-\x^*  \right\|^2.
\end{eqnarray}

\subsection{Other Material}\label{Other Material-sec}
\subsubsection{Sparse Update}
Proof of the algorithm 1 and Algorithm 2 are equivalent.
We use $\w^{j(k)}_1$ and $\bdelta^k_1$ to denote  $\w^{j(k)}$ and $\bdelta^k$ generated by Algorithm 1,  and use $\w^{j(k)}_2$ and $\bdelta^k_2$ to denote  $\w^{j(k)}$ and $\bdelta^k$ generated by Algorithm 2. To prove the results, we use induction to show that $\z^k=\uu^k$, $\x^k = \uu^k + a^k\vv^k $.

When $k=0$, we have $\z^k=\uu^k=\mathbf{0}$, $\x^k = \uu^k + a^k\vv^k=\mathbf{0} $. For $j(0)=0$, then $\w_1^{j(0)}=\y^0=\w_2^{j(0)} = \mathbf{0} = \uu^0 +\d^1 \vv^0$. So we have $\bdelta^0_1$    = $\bdelta^0_2$. Then we have that $\z^1 = \uu^1$.  So
\begin{eqnarray}\label{one1}
\x^{1} = \y^0 + \theta^0 \bdelta^0_1 = \uu^0 + d^1 \vv^0 +\theta^0 \bdelta_2 = \uu^1 +d^1 \vv^1 -\bdelta_2^0 + \frac{d^1\bdelta^0_2}{d^0} +\theta^0\bdelta^0_2=\uu^1 +d^1 \vv^1,
\end{eqnarray}
where in the third equality, we use $d^{k+1}= d^{k}(1-\theta^k)$.

When $k>0$, suppose we have  $\z^k=\uu^k$, and $\x^k = \uu^k + d^k\vv^k$, then
\begin{eqnarray}\label{two1}
\y^k = (1-\theta^k)\x^k +\theta^k\z^k = (1-\theta^k)d^k \vv^k +\uu^k=d^{k+1} \vv^k +\uu^k.
\end{eqnarray}
If we obtain $\w^{j(k)}_1 =\w^{j(k)}_2$, then $\bdelta^k_1 = \bdelta^k_2$ and $\z^{k+1} = \uu^{k+1}$. For $\x^{k+1}$, we have
\begin{eqnarray}\label{thrid1}
\x^{k+1} &=& \y^k + \theta^k \bdelta^k_1 = \uu^k + d^{k+1} \vv^k +\theta^k \bdelta_2 \\
&& = \uu^{k+1} +d^{k+1} \vv^{k+1} -\bdelta_2^k + \frac{d^{k+1}\bdelta^k_2}{d^k} +\theta^k\bdelta^k_2=\uu^{k+1} +d^{k+1} \vv^{k+1}.\notag
\end{eqnarray}
Now we are to prove $\w^{j(k)}_1 =\w^{j(k)}_2$. We introduce an auxiliary algorithm, shown in Algorithm 3.
\begin{algorithm}[!h]
	\caption{AAGD-auxiliary}
	\label{AAGD-auxiliary}
	\begin{algorithmic}
		\STATE  $\mathbf{Input}$ $\x_1^{j(k)}= \x^{j(k)} $ and $\z_1^{j(k)}= \z^{j(k)}$.
		\STATE  $  \mathbf{for}$ $k=j(k)$ $\mathbf{to}$ $K-1$ $\mathbf{do}$\\
		\STATE1  $\quad\y_1^{k} = (1-\theta^k)\x^k +\theta^k \z^k$
		\STATE2  $\quad \z_1^{k+1} =  \z_1^k  $. 
		\STATE3  $\quad  \x_1^{k+1} = \y_1^k + \theta^k(\z_1^{k+1}-\z_1^k)  $. 
		\STATE $\mathbf{end \ for}$
	\end{algorithmic}
\end{algorithm}
The algorithm is the  serial AGD by setting $\bdelta_1^k=0$.  The result of Eq.~\eqref{step1} can be directly used by setting $\x_1^i -\y_1^{i-1}= \mathbf{0}$ when $i>j(k)$. So we obtain that $\y_1^k =\w^{j(k)}_1$. Now we are to prove that $\y_1^k = \w^{j(k)}_2$, that is to prove that
\begin{eqnarray}
\y_1^k  = \uu^{j(k)} +d^{k+1}\vv^{j(k)}.
\end{eqnarray}
To proof this, we show that Algorithm 3 is equivalent to Algorithm 4.
\begin{algorithm}[!h]
	\caption{AAGD-auxiliary2}
	\label{AAGD-auxiliary2}
	\begin{algorithmic}
		\STATE  $\mathbf{Input}$ $\theta^k$, $\uu_1^{j(k)}= \uu^{j(k)} $ and $\vv_1^{j(k)}= \vv_1^{j(k)}$, $d_1^{j(k)} =d^{j(k)}$.
		\STATE  $  \mathbf{for}$ $k=0$ $\mathbf{to}$ $K-1$ $\mathbf{do}$\\
		\STATE1  $\quad d^{k+1}= d^{k}(1-\theta^k)$, 
		\STATE2  $\quad\y_2^k =  \uu_1^{k} + d^{k+1} \vv_1^{k} .$
		\STATE3  $\quad \uu_1^{k+1} =  \uu_1^{k}  $. 
		\STATE4  $\quad \vv_1^{k+1} =  \vv_1^{k} $. 
		\STATE $\mathbf{end \ for}$
	\end{algorithmic}
\end{algorithm}
By the induction same as Eq.~\eqref{one1}, \eqref{two1}, \eqref{thrid1}, we can obtain that $\x_1^k  =  \uu_1^k +d^k \vv_1^k$, $\z_1^k  = \uu_1^k$, and $\y_1^k = \uu_1^k +d^{k+1} \vv_1^k$. As $\uu^k  = \uu^{j(k)}$ and $\vv^k  = \vv^{j(k)}$, we obtain $\y_1^k  = \uu^{j(k)} +d^{k+1}\vv^{j(k)}$. This ends proof.

\subsubsection{Pre-define Update Order}
Our technique need to predefine the update order to obtain $\w^{j(k)}$.  Once such an order has been set, each thread may update the gradient estimator accordingly. If one thread returns the gradient  early, the gradient can be stored and it will go on for the next iteration. The master thread will use the gradient to update parameters after receiving all the required gradient.

However, though the threads will never be hanged up, the large inconsistency of real order will amplify the delay effect. We found that for dense datasets, the computation costs are roughly the same for each child node, so one may directly set $k = j(k)+\zeta-1$, where $\zeta$ is the number of cores. This works well in practice. While for sparse datasets, simply setting the predefined order is not advised.  We introduce way to avoid predefining the order. 

Through our algorithm, we can find that when smooth-part of the objective function is quadratic, such as $f(\x) = \|\A\x\|^2$, then
\begin{eqnarray}
\nabla f(\w^{j(k)}) =  \nabla f(\uu^{j(k)}) + d^{k+1} \nabla f(\vv^{j(k)}).
\end{eqnarray}
So we can first compute the  the $ \nabla f(\uu^{j(k)})$ and  $\nabla f(\vv^{j(k)})$, and then add them together. This can avoid predefining the order. When the  smooth-part of the objective function are not quadratic, we can  uses Hessian-Vector~\citep{pearlmutter1994fast} product to approximate the gradient. Set AASVRG as an example, applying Hessian-Vector product~\citep{pearlmutter1994fast} to approximate $\nabla f_{i^k_s}(\w^s_{j(k)})$,  we have
\begin{eqnarray}\label{1333}
\nabla f_{i^k_s}(\w^s_{j(k)}) 
\approx \nabla f_{i^k_s} (\p^s_{j(k)})  - \alpha^s_k \HH_{i^k_s}(\p^s_{j(k)})(\x^s_{j(k)} -\x^s_{j(k)-1}), 
\end{eqnarray}
where $\p^s_{j(k)} =  \x^s_{j(k)}+ \frac{a^s\left(1-(a^s)^{\tau+1}\right)}{1-a^s} (\x^s_{j(k)} -\x^s_{j(k)-1})$, and $$\alpha^s_k= \frac{(a^s)^{k-j(k)+2}\left(1-(a^s)^{\tau-(k-j(k))}\right)}{1-a^s},$$  $\HH_{i_{k,s}}(\ww^s_{j(k)})$ denotes the Hessian Matrix of $f_{i_{k,s}}$ at point $\ww^s_{j(k)}$.  In this way, $k-j(k)$ does not need to be known before computing the gradient estimator.

We can find that Eq.~\eqref{1333} has the following property: 1) $\alpha^s_k$ decreases exponentially with respect to the growth of delay $k-j(k)$. For severely delayed system~(lots of cores are running), we can assume that $k-j(k)$ is large, so $\alpha$ is small; 2) when $f$ is quadratic, Eq.~\eqref{1333} holds strictly, so E-ASVRG also achieves the accelerated convergence rate. For lots of machine learning problems, the Hessian-Vector product can be efficiently computed through  Hessian Free techniques~\citep{pearlmutter1994fast,martens2010deep}, which is in $O(d)$ time, the same as computing the gradient, where $d$ is the dimension of the parameter.

\subsubsection{Implementation Details}
\bfseries Deadlock Avoidance \mdseries  To avoid deadlock,  we  associate an ordering for all the locks such that each thread follows the same ordering to acquire the locks.\\
\bfseries Sparse Update \mdseries   We can find that by changing variable,   it is able to spare date on the sparse dataset. For ASCDA, like \citep{APCG}, we can introduce $x_1 = \A \uu$, and $\x_2=\A_vv$ to fast obtain the gradient. When $\theta$ is fixed, such as for SC and AASVRG, the update of $\vv^k$ will cause numerical problems because $d^k\to 0$, we can store $d^k \vv^k$  as $\vv_1^k$ and $\vv_2^k$, with the first one store the value, and the second store the power.\\
\bfseries Spin locks \mdseries   Also observed by \citep{ADSCD}, when there are no locks,  due to the memory conflict, $\x_1 \neq \A \uu$ and $\x_2 \neq \A \vv$, this is harmful and will lead the algorithm solving a deflected problem. To tackle it,  we  create $n$'s spin lock, and add lock when the corresponding coordinate of $\uu_1$ and $\uu_2$ are updated.

\subsubsection{ Sparsity $ \triangle$-assumption\citep{asvrg}}
The sparsity $ \triangle$-assumption in \citep{asvrg} is as follows: for problem of composite finite-sum problem, i.e. Eq. (16) in the paper, suppose $f_i$ only depends on $\x_{e_i}$, where $e_i \subseteq [d]$, i.e.,  $f_i$ acts only on  the components of $\x$ indexed by the set $e_i$. Let $\left\|\x\right\|_i^2$ denote $\sum_{j=\in e_j}\left\|\x_j\right \|^2$; then the convergence depends on $\triangle$, the smallest constant such that $\E_i [ \left \|\x\right\|_i^2] = \triangle \left\|\x\right\|^2$, and $\triangle\ll1$.

One can find by the assumption that changes of  each update are small, so through the proof,  step 1, e.g. Eq.~\eqref{step31}, $\left\|\y^k -\w^{j(k)}\right\|^2$ will $\triangle$ times smaller. So Proposition 2 in the paper is obtained.

\subsection{More Experimental Results}
\begin{table}[h]
	\centering
	\caption{Details of the dense datasets. (Dim., is short for dimensionality)}
	\begin{tabular}{rrrrr}
		\toprule
		Datasets   & \#training  & Dim. & Class & \#mini-batch \\
		\midrule
		USPS      & 7291   & 256 & 10 & 50  \\
		MNIST     & 60000  & 784 & 10 & 50  \\
		SENSIT  & 78823  & 100 & 3  & 50  \\
		$\!\!$EPSILON   & 400000 & 2000&  2 & 200  \\			
		\bottomrule
	\end{tabular}
	\label{tab:dataset}
\end{table}

\begin{figure*}[!h]
	\center
	
	\subfigure[news20~(time)]{\includegraphics[width=2.5in]{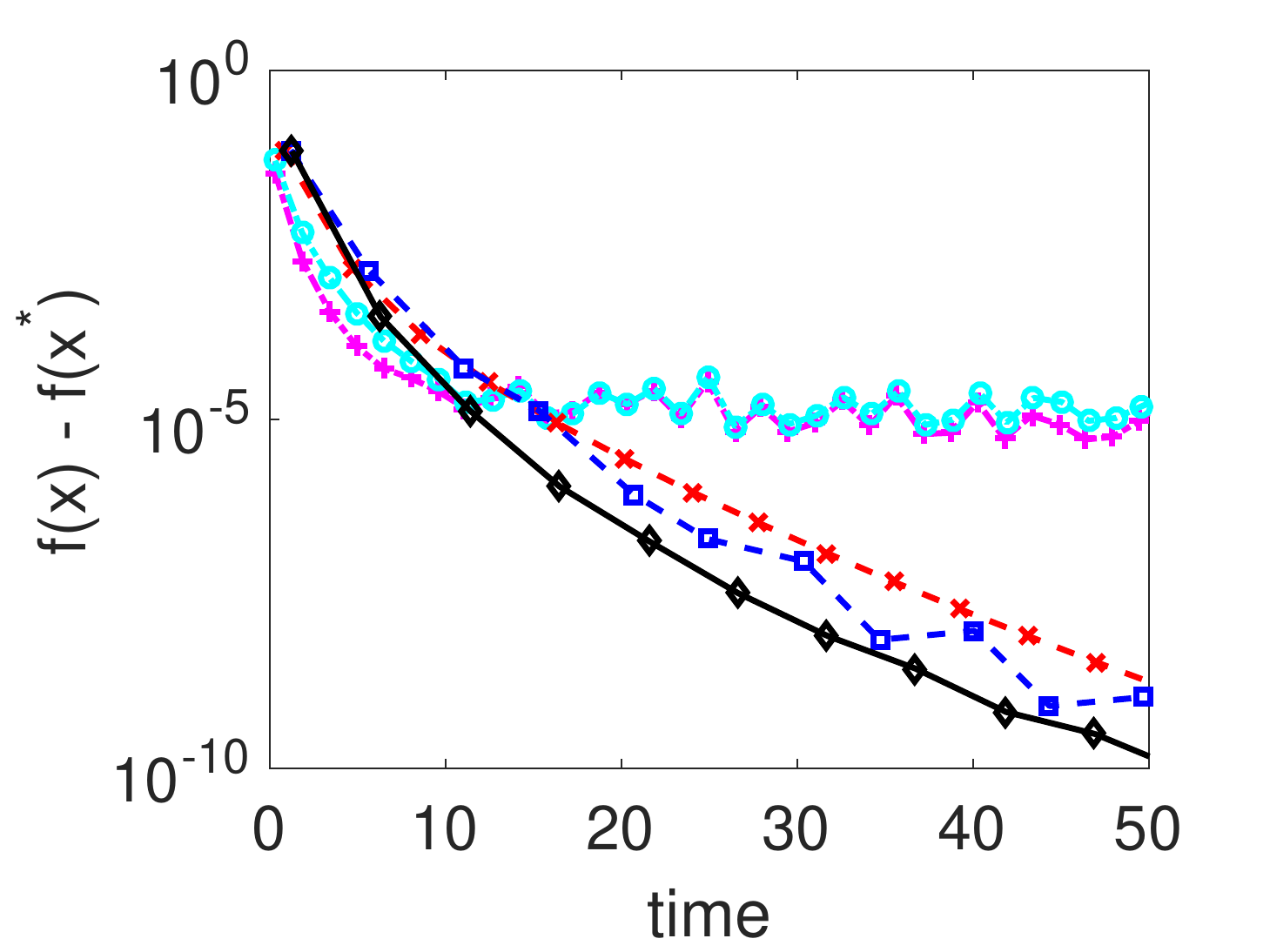}}
		\hspace{-0.15in}
	\subfigure[news20~(iters)]{\includegraphics[width=2.5in]{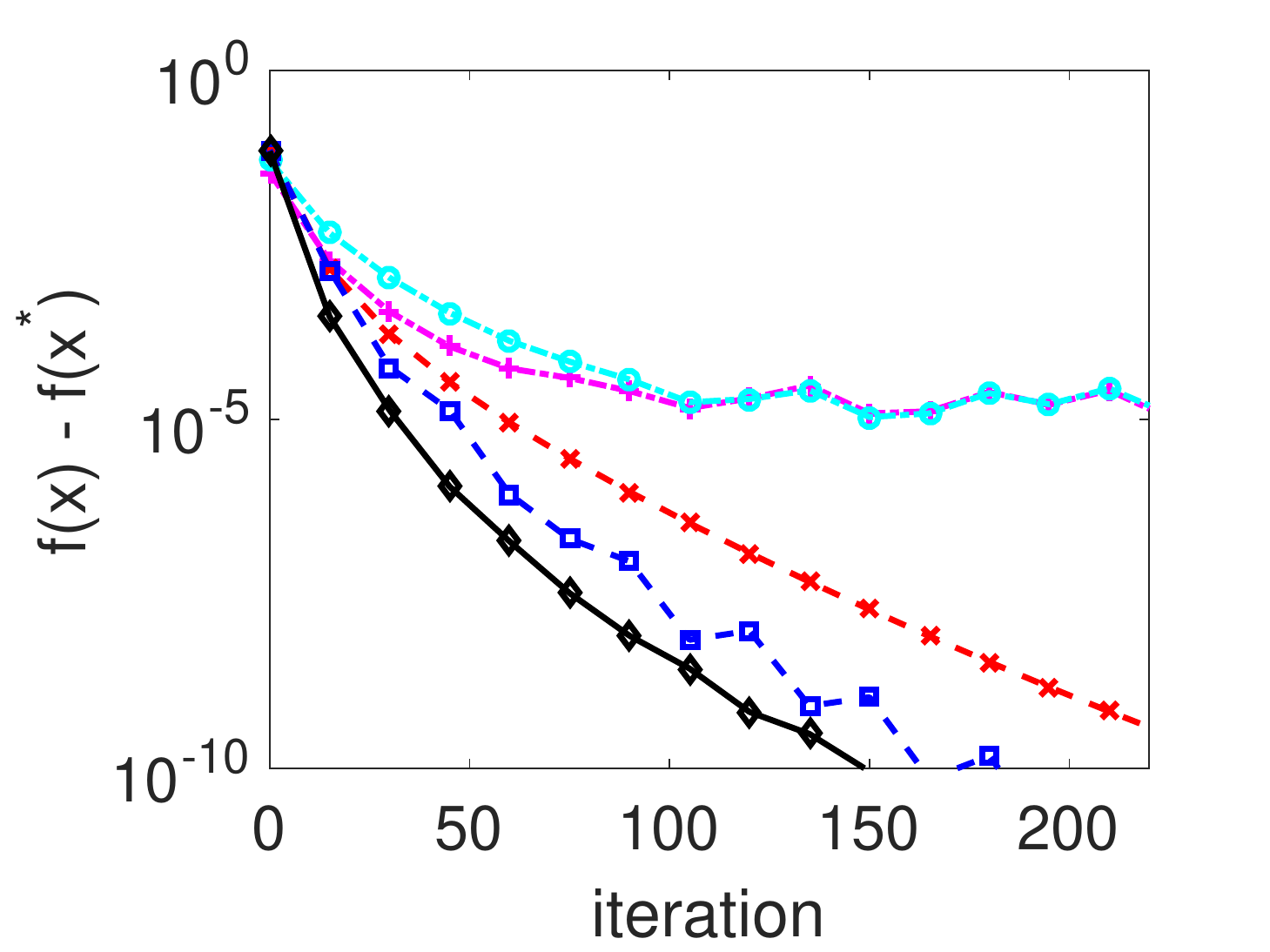}}
		\hspace{-0.15in}
	\subfigure[usps~(time)]{\includegraphics[width=2.5in]{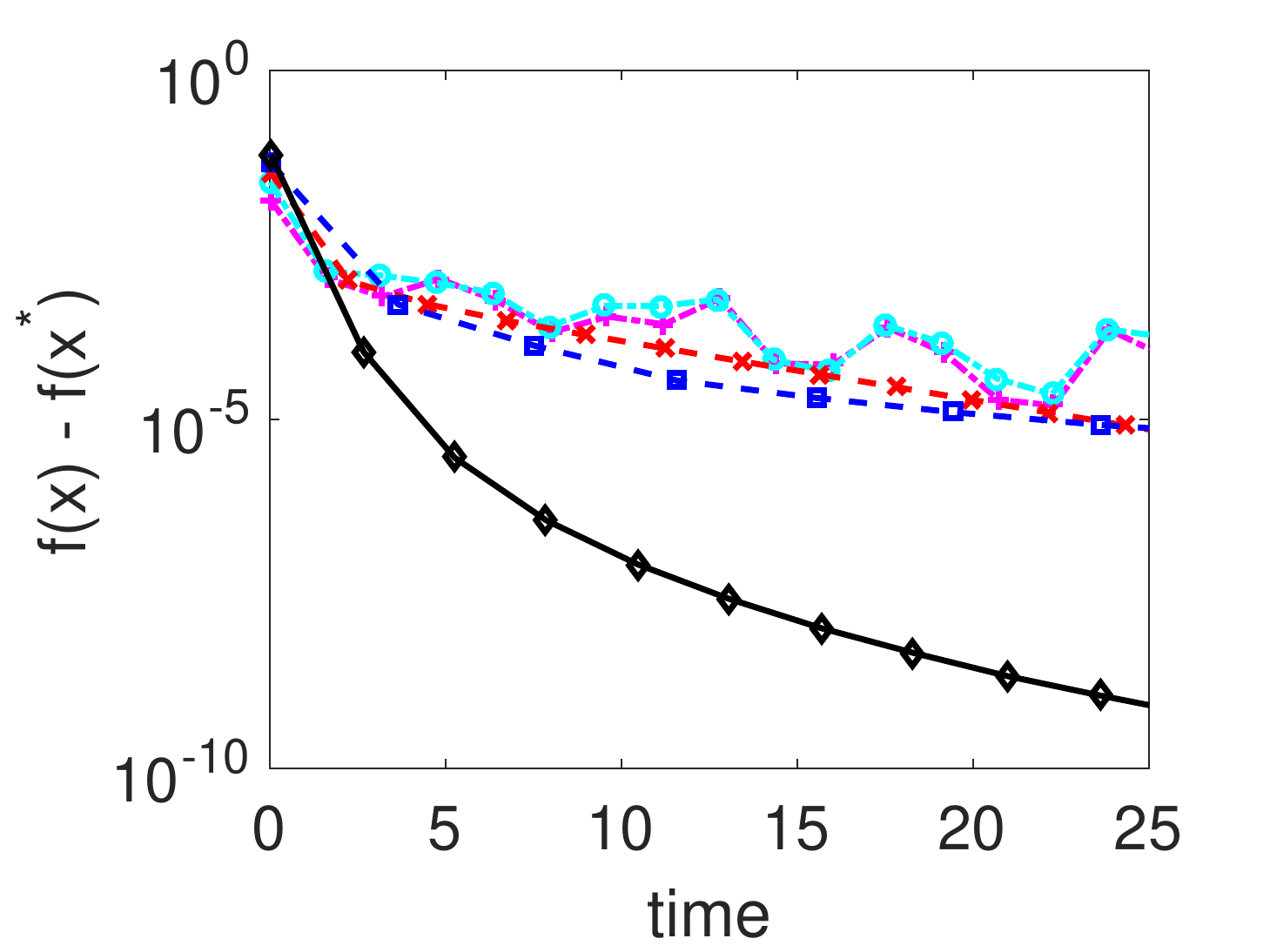}}
		\hspace{-0.15in}
	\subfigure[usps~(iters)]{\includegraphics[width=2.5in]{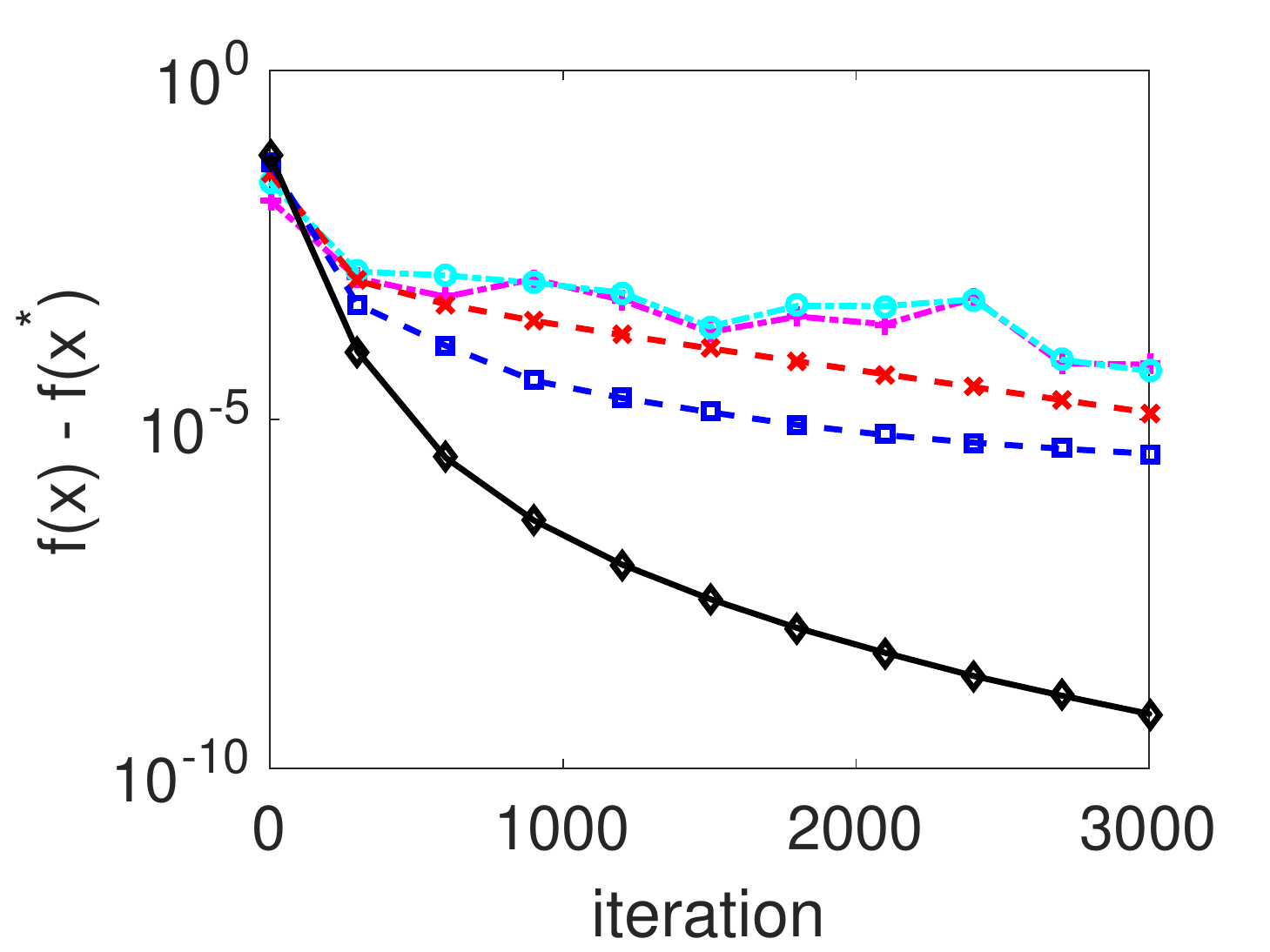}}
	\subfigure[combined~(time)]{\includegraphics[width=2.5in]{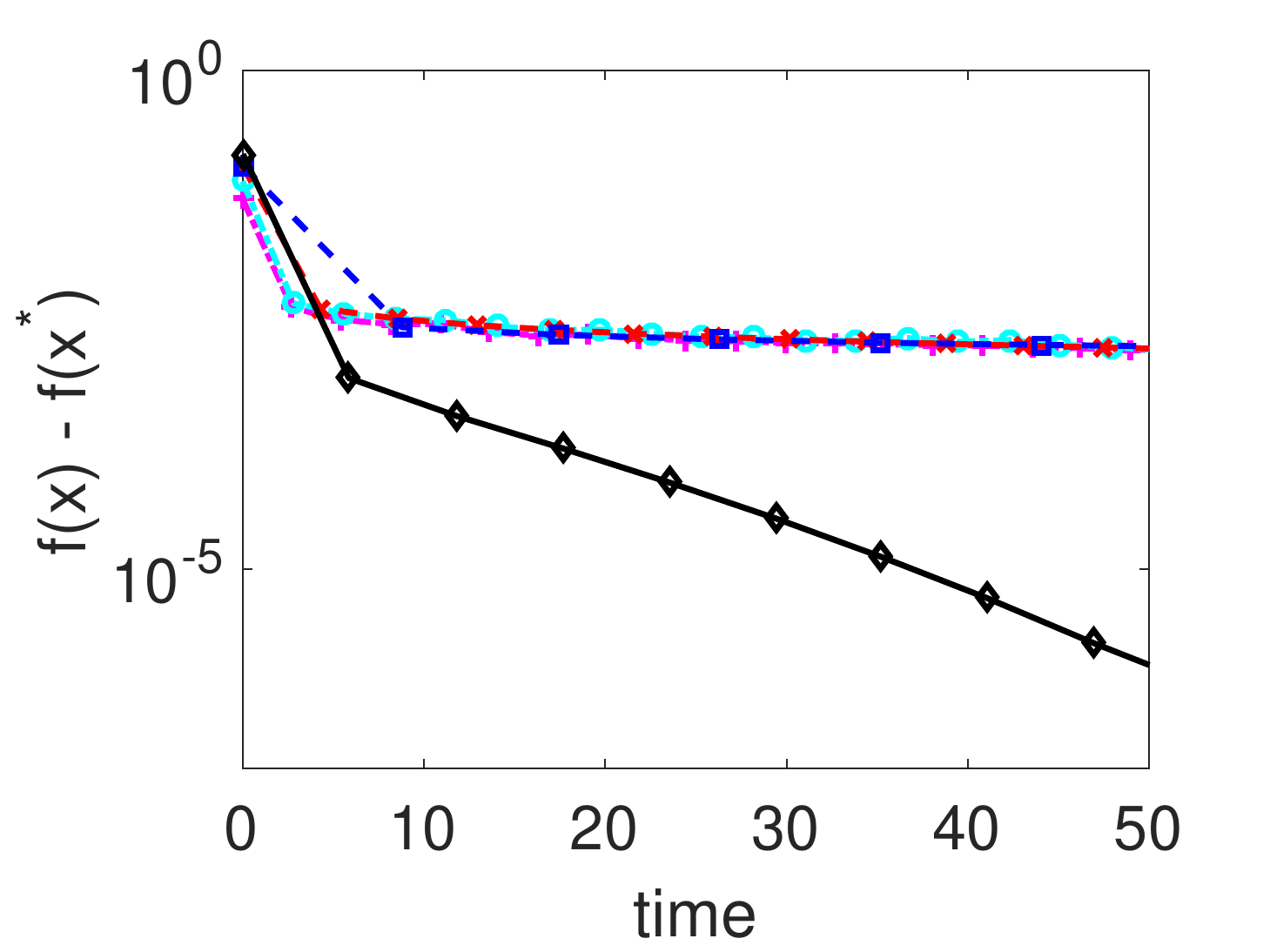}}
		\hspace{-0.15in}
	\subfigure[combined~(iters)]{\includegraphics[width=2.5in]{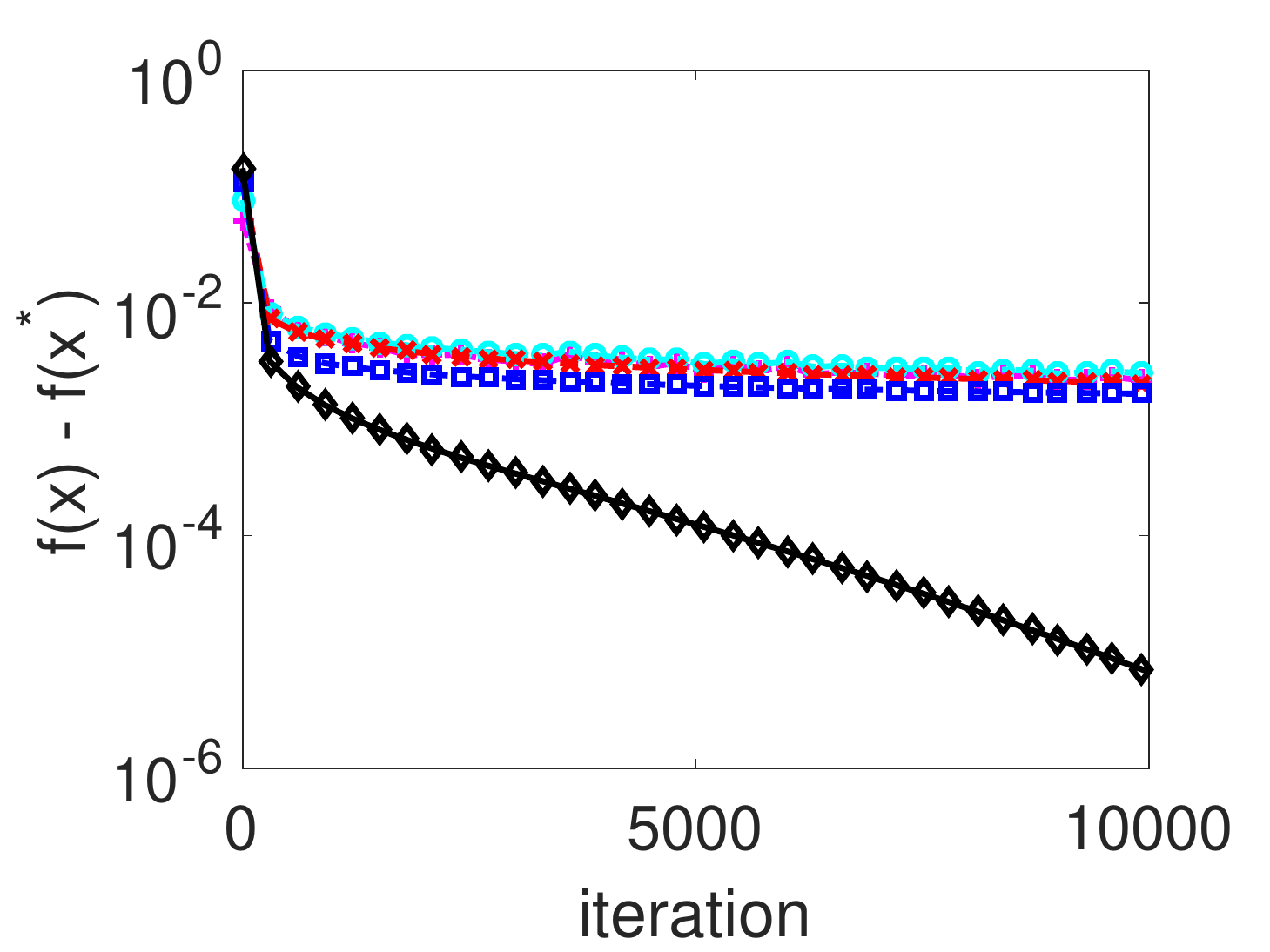}}
		\hspace{-0.15in}
	\includegraphics[height=0.4in]{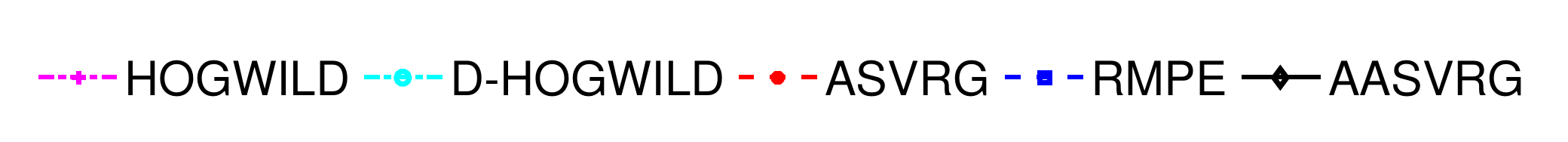}
	\vspace{-0.3cm}
	\caption{Residuals vs CPU training time~(s) and iterations for solving Ridge Regression problem. ``news20'' is a sparse dataset, while ``usps'' and ``combined'' are dense.}
	\label{fig:primal}
	\vspace{-0.3cm}
\end{figure*}

We have also verified the convergence speed for AASVRG on another three datasets, namely the sparse dataset \textit{new20} and dense datasets \textit{usps}, \textit{combined}. The results are shown in Fig.~\ref{fig:primal}. It turns out that our algorithm has competitive results on all of these datasets.

\begin{figure*}[!h]
	\center
	\hspace{-0.15in}
	\subfigure[$\lambda = \frac{1}{n}$]    {\includegraphics[width=2.5in]{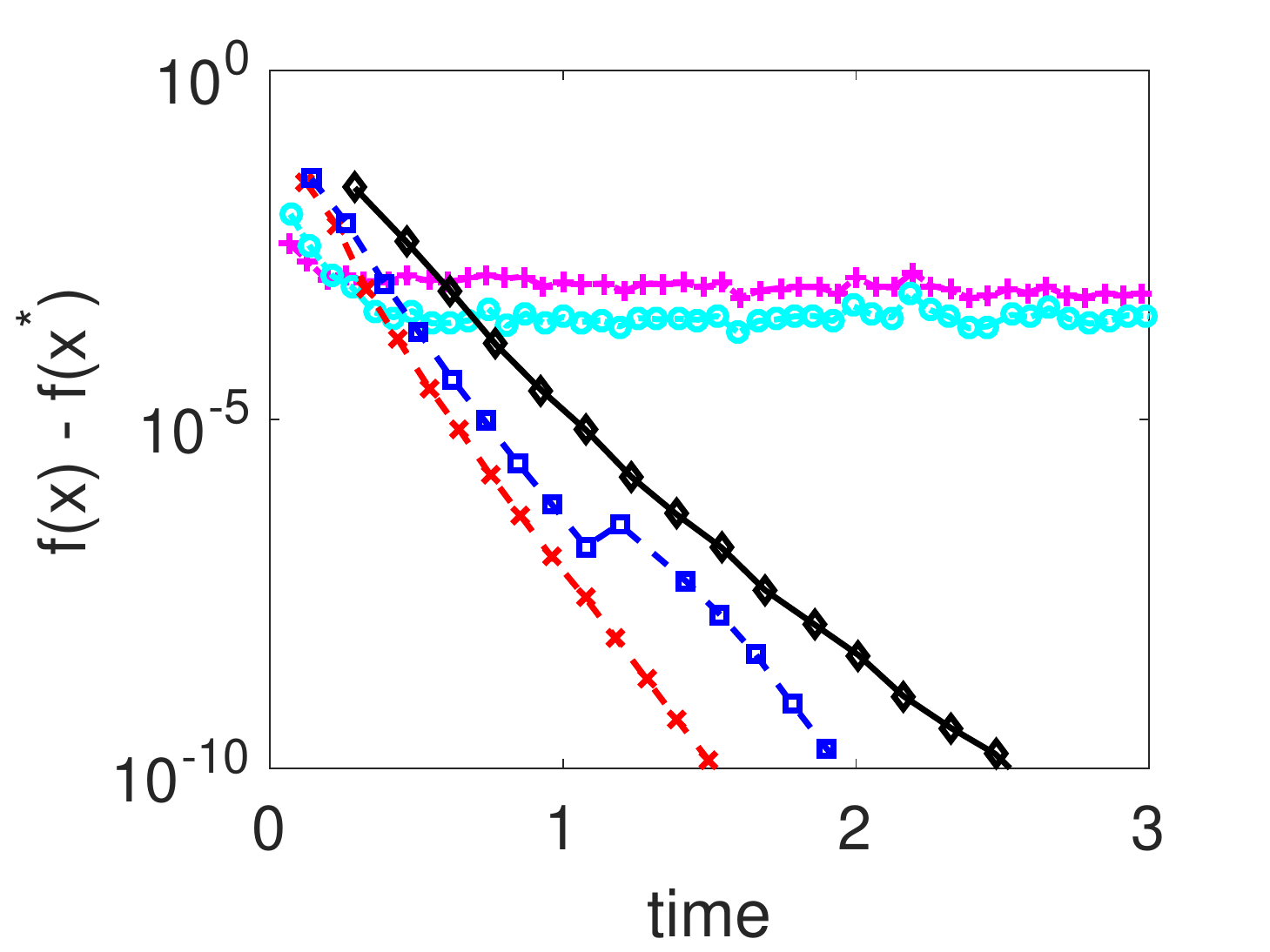}}
	\hspace{-0.15in}
	\subfigure[$\lambda = \frac{1}{n}$]    {\includegraphics[width=2.5in]{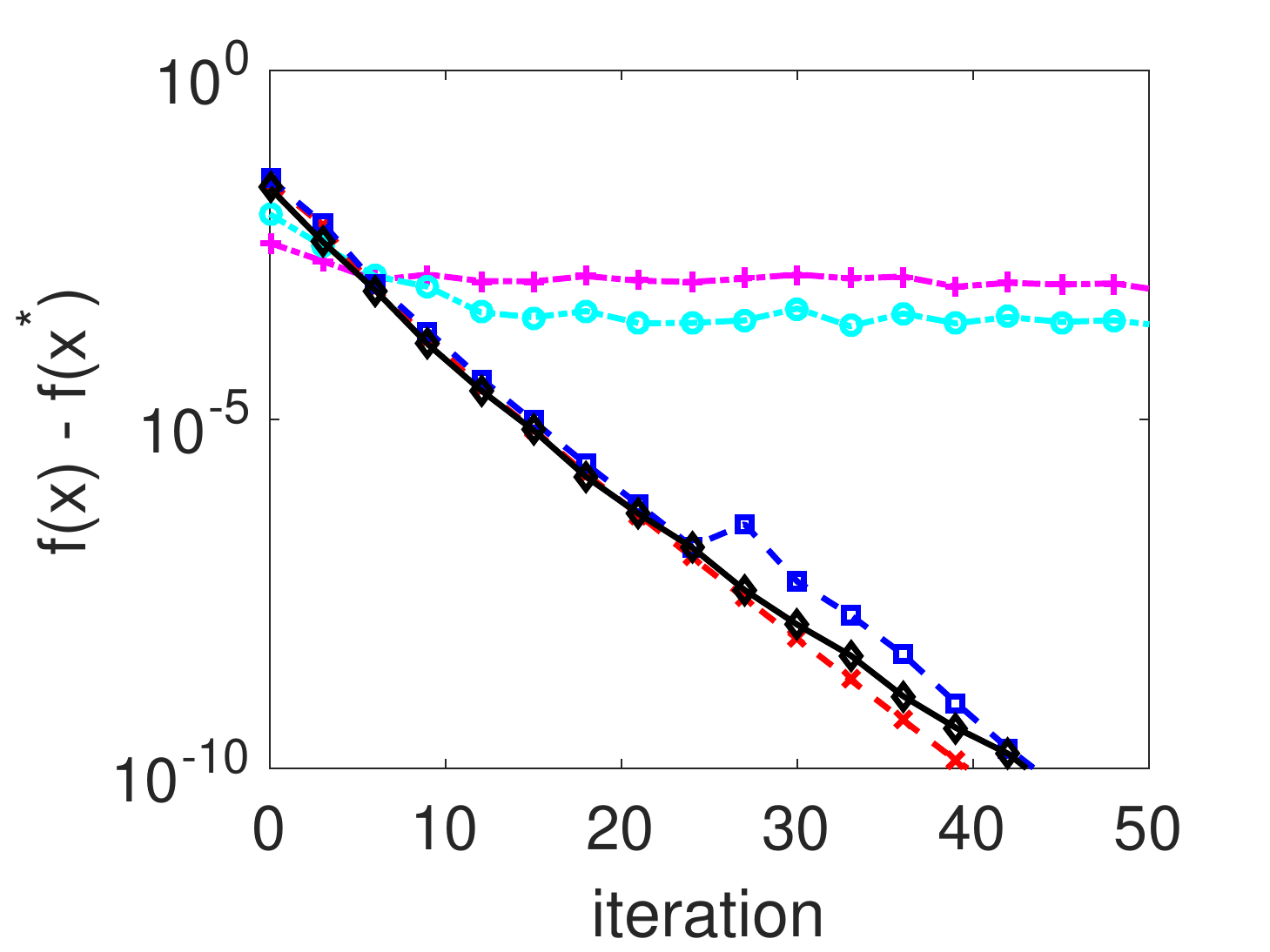}}
	\hspace{-0.15in}
	\subfigure[$\lambda = \frac{1}{10n}$]  {\includegraphics[width=2.5in]{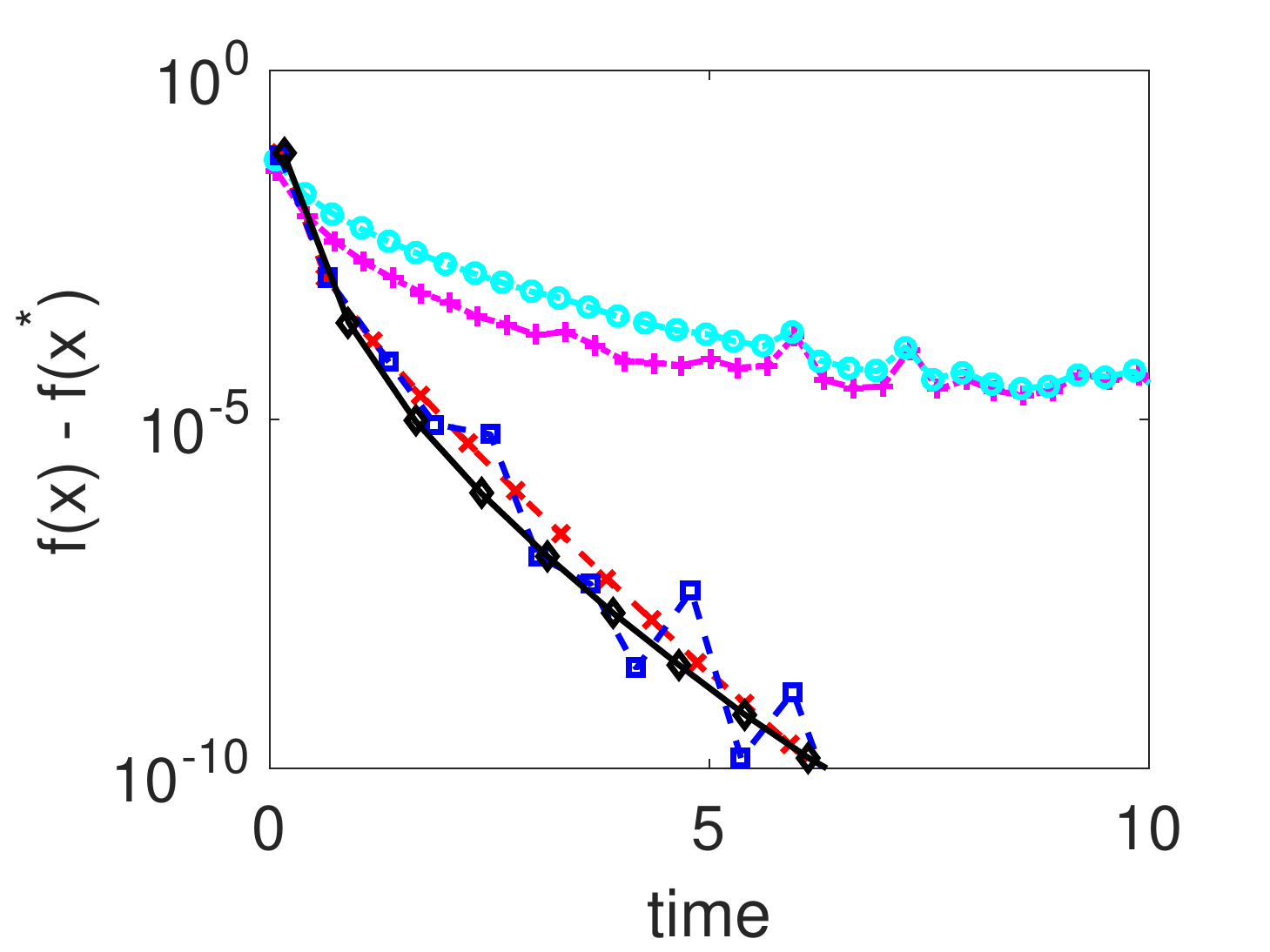}}
	\hspace{-0.15in}
	\subfigure[$\lambda = \frac{1}{10n}$]  {\includegraphics[width=2.5in]{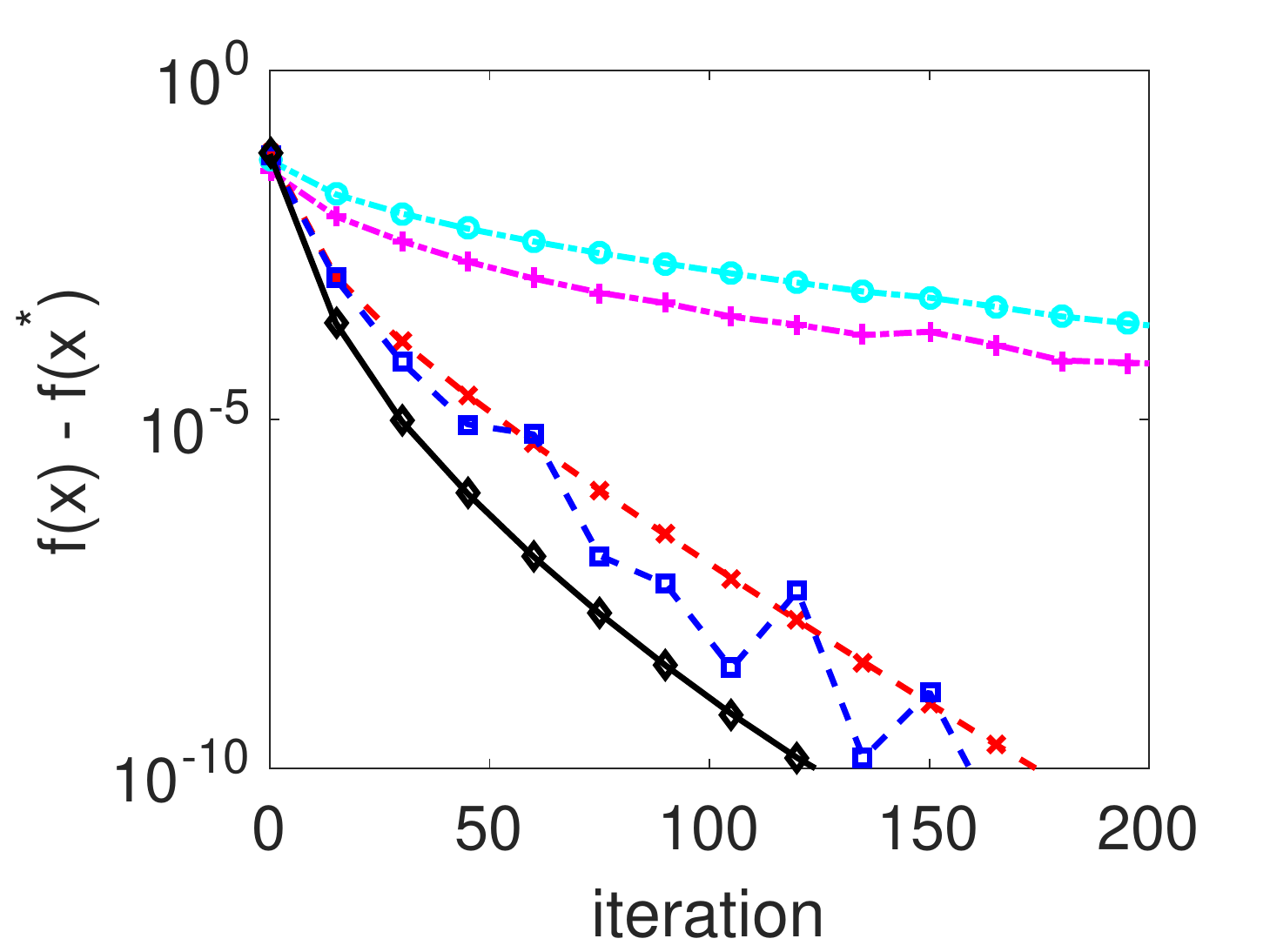}}
	\hspace{-0.15in}
	\subfigure[$\lambda = \frac{1}{100n}$] {\includegraphics[width=2.5in]{figure3/rcv1_primal_100n_time.pdf}}
	\hspace{-0.15in}
	\subfigure[$\lambda = \frac{1}{100n}$] {\includegraphics[width=2.5in]{figure3/rcv1_primal_100n_iter.pdf}}
	\hspace{-0.15in}
	\subfigure[$\lambda = \frac{1}{1000n}$]{\includegraphics[width=2.5in]{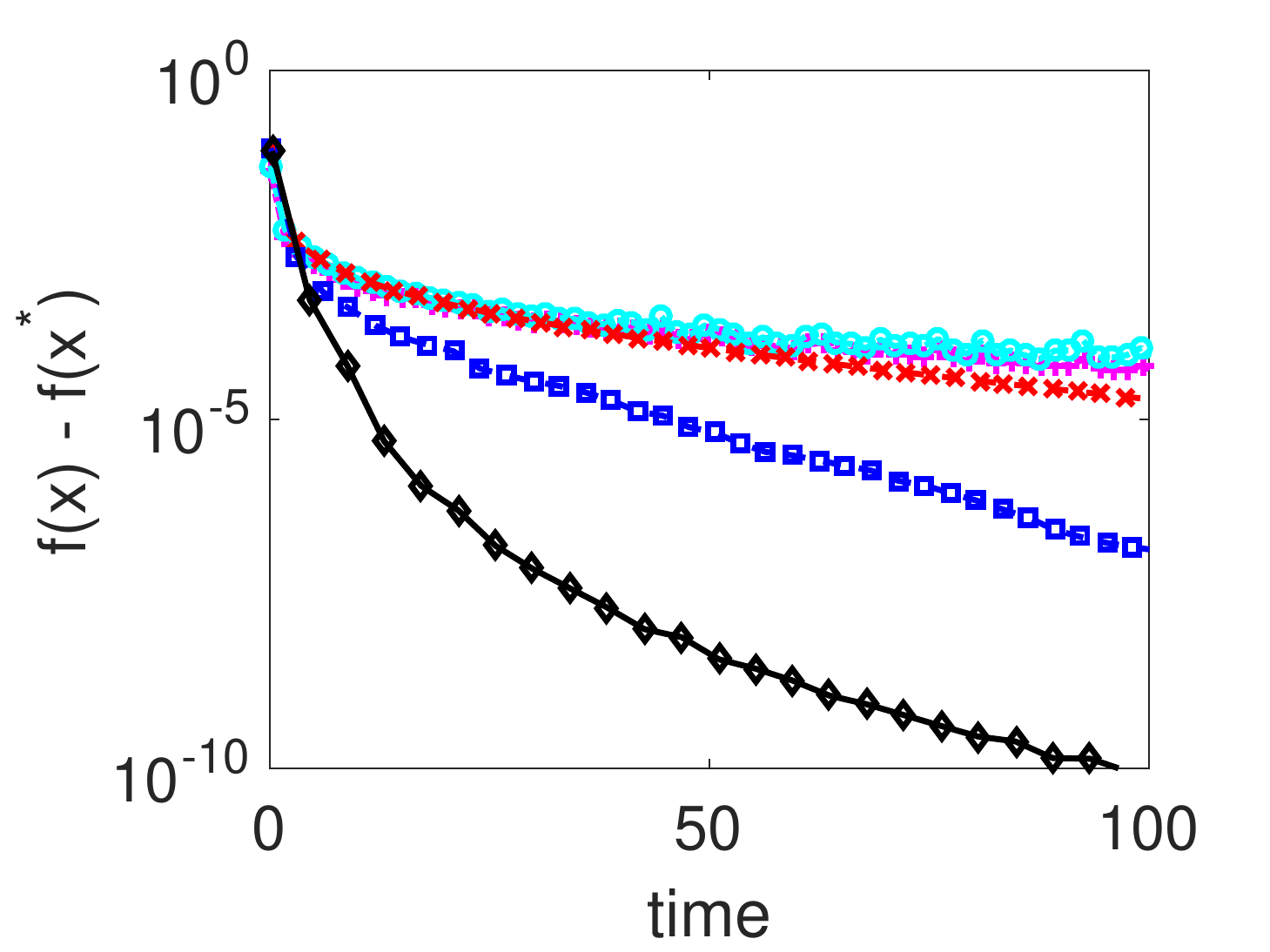}}
	\hspace{-0.15in}
	\subfigure[$\lambda = \frac{1}{1000n}$]{\includegraphics[width=2.5in]{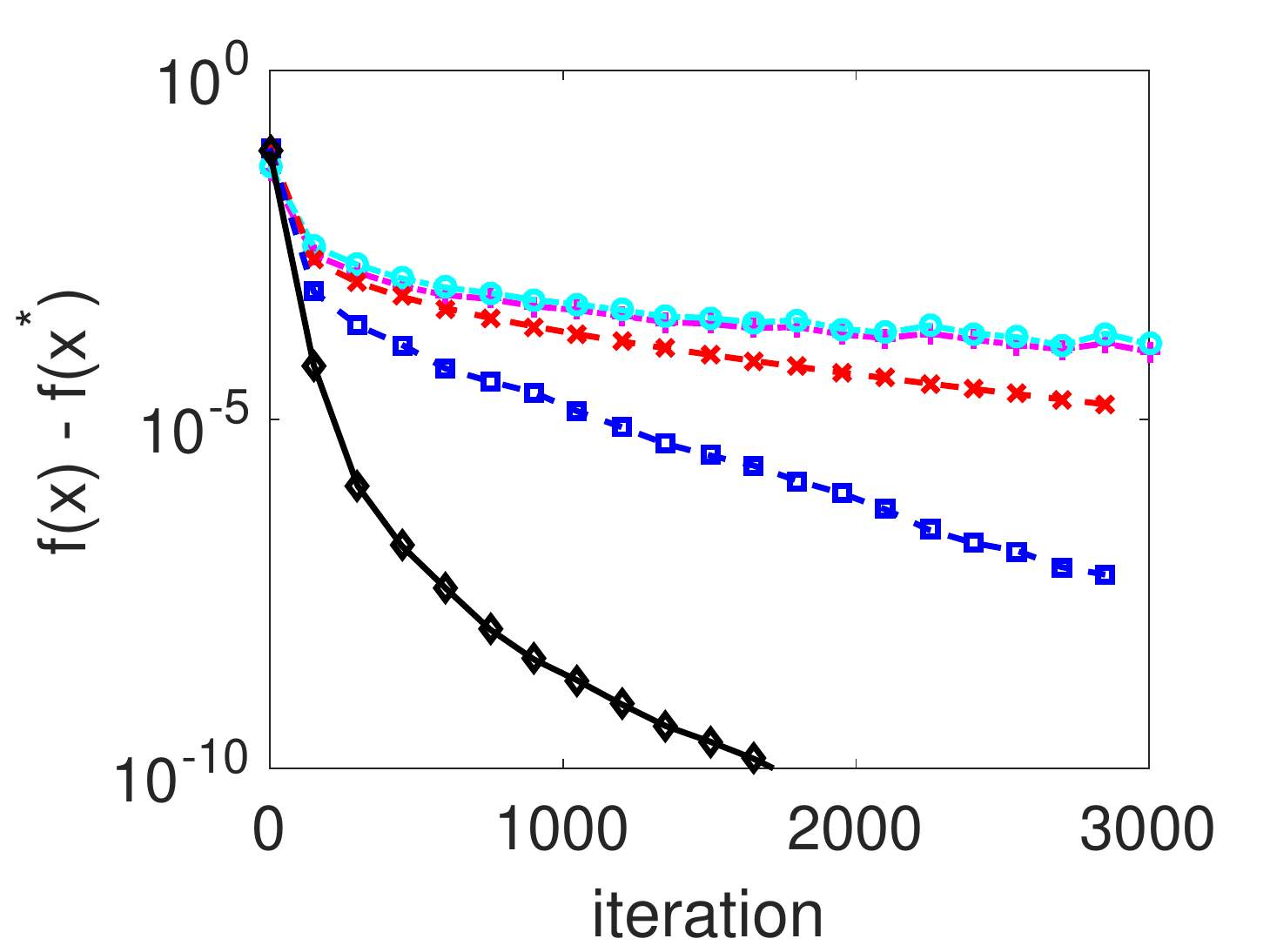}}
	\vspace{-0.1cm}
	
	\includegraphics[height=0.4in]{figure3/primal_legend_crop.pdf}
	\vspace{-0.3cm}
	\caption{Residuals vs CPU training time and iterations for solving Ridge Regression problem with different $\lambda$ on rcv1 datasets.}
	\label{fig:lambda}
	\vspace{-0.3cm}
\end{figure*}

Our algorithm has big advantages for ill-condition problem, i,e, when the regularization term $\lambda$ in Ridge Regression is small. As we can see from Fig.~\ref{fig:lambda}, when $\lambda$ is large, our algorithm has similar performance as other state-of-the-art algorithms. However, when $\lambda$ is small, we gains the huge advantages in terms of the convergence.
